
%
\def\conv{\mathop{\vrule height2,6pt depth-2,3pt 
    width 5pt\kern-1pt\rightharpoonup}}

\advance\vsize by 1 true cm
%
\def\dess #1 by #2 (#3){
  \vbox to #2{
    \hrule width #1 height 0pt depth 0pt
    \vfill
    \special{picture #3} 
    }
  }

\def\dessin #1 by #2 (#3 scaled #4){{
  \dimen0=#1 \dimen1=#2
  \divide\dimen0 by 1000 \multiply\dimen0 by #4
  \divide\dimen1 by 1000 \multiply\dimen1 by #4
  \dess \dimen0 by \dimen1 (#3 scaled #4)}
  }
%
\def \trait (#1) (#2) (#3){\vrule width #1pt height #2pt depth #3pt}
\def \fin{\hfill
	\trait (0.1) (5) (0)
	\trait (5) (0.1) (0)
	\kern-5pt
	\trait (5) (5) (-4.9)
	\trait (0.1) (5) (0)
\medskip}
%


\font\sevenbf=cmbx7

\baselineskip=15pt
\abovedisplayskip=15pt plus 4pt minus 9pt
\belowdisplayskip=15pt plus 4pt minus 9pt
\abovedisplayshortskip=3pt plus 4pt
\belowdisplayshortskip=9pt plus 4pt minus 4pt
\let\epsilon=\varepsilon

\def\biblio #1 #2\par{\parindent=30pt\item{}\kern -30pt\rlap{[#1]}\kern
30pt #2\smallskip}
 %
\catcode`\@=11
\def\@lign{\tabskip=0pt\everycr={}}
\def\equations#1{\vcenter{\openup1\jot\displ@y\halign{\hfill\hbox
{$\@lign\displaystyle##$}\hfill\crcr
#1\crcr}}}
\catcode`\@=12
%
\def\pmb#1{\setbox0=\hbox{#1}%
\hbox{\kern-.04em\copy0\kern-\wd0
\kern.08em\copy0\kern-\wd0
\kern-.02em\copy0\kern-\wd0
\kern-.02em\copy0\kern-\wd0
\kern-.02em\box0\kern-\wd0
\kern.02em}}
%
\def\undertilde#1{\setbox0=\hbox{$#1$}
\setbox1=\hbox to \wd0{$\hss\mathchar"0365\hss$}\ht1=0pt\dp1=0pt
\lower\dp0\vbox{\copy0\nointerlineskip\hbox{\lower8pt\copy1}}}
%

%

\def\maj#1#2,{\rm #1\sevenrm #2\rm{}}
\def\Maj#1#2,{\bf #1\sevenbf #2\rm{}}
\outer\def\lemme#1#2 #3. #4\par{\medbreak
\noindent\maj{#1}{#2},\ #3.\enspace{\sl#4}\par
\ifdim\lastskip<\medskipamount\removelastskip\penalty55\medskip\fi}

\def\Remark #1. {\noindent{\Maj REMARK,\ \bf #1. }}

\outer\def\Lemme#1#2 #3. #4\par{\medbreak
\noindent\Maj{#1}{#2},\ \bf #3.\rm\enspace{\sl#4}\par
\ifdim\lastskip<\medskipamount\removelastskip\penalty55\medskip\fi}



\def\Notation #1. {\noindent{\Maj NOTATION,\ \bf #1. }}

\def\Example #1. {\noindent{\Maj EXAMPLE,\ \bf #1. }}

\hfuzz=1cm


\catcode`\ˆ=\active     \def ˆ{\`a}
\catcode`\‰=\active     \def ‰{\^a}
\catcode`\=\active     \def {\c c}
\catcode`\Ž=\active    \def Ž{\'e} 

\catcode`\=\active   \def {\`e}
\catcode`\=\active   \def {\^e}
\catcode`\'=\active   \def '{\"e}
\catcode`\"=\active   \def "{\^\i}
\catcode`\•=\active   \def •{\"\i}
\catcode`\™=\active   \def ™{\^o}
\catcode`\š=\active   \defš{}
\catcode`\=\active   \def {\`u}
\catcode`\ž=\active   \def ž{\^u}
\catcode`\Ÿ=\active   \def Ÿ{\"u}
\catcode`\ =\active   \def  {\tau}
\catcode`\¡=\active   \def ¡{\circ}
\catcode`\¢=\active   \def ¢{\Gamma}
\catcode`\¤=\active   \def ¤{\S\kern 2pt}
\catcode`\¥=\active   \def ¥{\puce}
\catcode`\§=\active   \def §{\beta}
\catcode`\¨=\active   \def ¨{\rho}
\catcode`\©=\active   \def ©{\gamma}
\catcode`\­=\active   \def ­{\neq}
\catcode`\°=\active   \def °{\ifmmode\ldots\else\dots\fi}
\catcode`\±=\active   \def ±{\pm}
\catcode`\²=\active   \def ²{\le}
\catcode`\³=\active   \def ³{\ge}
\catcode`\µ=\active   \def µ{\mu}
\catcode`\¶=\active   \def ¶{\delta}
\catcode`\·=\active   \def ·{\Sigma}
\catcode`\¸=\active   \def ¸{\Pi}
\catcode`\¹=\active   \def ¹{\pi}
\catcode`\»=\active   \def »{\Upsilon}
\catcode`\¾=\active   \def ¾{\alpha}
\catcode`\À=\active   \def À{\cdots}
\catcode`\Â=\active   \def Â{\lambda}
\catcode`\Ã=\active   \def Ã{\sqrt}
\catcode`\Ä=\active   \def Ä{\varphi}
\catcode`\Å=\active   \def Å{\xi}
\catcode`\Æ=\active   \def Æ{\Delta}
\catcode`\Ç=\active   \def Ç{\cup}
\catcode`\È=\active   \def È{\cap}
\catcode`\Ï=\active   \def Ï{\oe}
\catcode`\Ñ=\active   \def Ñ{\to}
\catcode`\Ò=\active   \def Ò{\in}
\catcode`\Ô=\active   \def Ô{\subset}
\catcode`\Õ=\active   \def Õ{\superset}
\catcode`\Ö=\active   \def Ö{\over}
\catcode`\×=\active   \def ×{\nu}
\catcode`\Ù=\active   \def Ù{\Psi}
\catcode`\Ú=\active   \def Ú{\Xi}
\catcode`\Ü=\active   \def Ü{\omega}
\catcode`\Ý=\active   \def Ý{\Omega}
\catcode`\ß=\active   \def ß{\equiv}
\catcode`\à=\active   \def à{\chi}
\catcode`\á=\active   \def á{\Phi}
\catcode`\ä=\active   \def ä{\infty}
\catcode`\å=\active   \def å{\zeta}
\catcode`\æ=\active   \def æ{\varepsilon}
\catcode`\è=\active   \def è{\Lambda}  
\catcode`\é=\active   \def é{\kappa}
\catcode`\ë=\active   \defë{\Theta}
\catcode`\ì=\active   \defì{\eta}
\catcode`\í=\active   \defí{\theta}
\catcode`\î=\active   \defî{\times}
\catcode`\ñ=\active   \defñ{\sigma}
\catcode`\ò=\active   \defò{\psi}

\def\date{\number\day\
\ifcase\month \or janvier \or f\'evrier \or mars \or avril \or mai \or juin \or juillet \or ao\^ut  \or
septembre \or octobre \or novembre \or d\'ecembre \fi
\ \number\year}
\font \GGgras=cmb10 at 14pt
\font \Ggras=cmb10 at 12pt
\font \ggras=cmb10 at 11pt

\def\Ga{{\bf \Ggras a}}
\def\Gb{{\bf \Ggras b}}
\def\Gc{{\bf \Ggras c}}

\def\Ge{{\bf \Ggras e}}

\def\Gn{{\bf \Ggras n}}

\def\Gt{{\bf \Ggras t}}

\def\Gx{{\bf \Ggras x}}

\def\GA{{\bf \Ggras A}}
\def\GB{{\bf \Ggras B}}

\def\GE{{\bf \Ggras E}}
\def\GF{{\bf \Ggras F}}
\def\GG{{\bf \Ggras G}}

\def\GI{{\bf \Ggras I}}

\def\GM{{\bf \Ggras M}}

\def\GP{{\bf \Ggras P}}
\def\GQ{{\bf \Ggras Q}}

\def\GU{{\bf \Ggras U}}
\def\GV{{\bf \Ggras V}}
\def\GW{{\bf \Ggras W}}

\def\GR{{\bf \Ggras R}}

\def\liminf{\mathop{\underline{\rm lim}}}

\def\sym{\fam\comfam\com}
\font\tensym=msbm10
\font\sevensym=msbm7
\font\fivesym=msbm5
\newfam\symfam
\textfont\symfam=\tensym
\scriptfont\symfam=\sevensym
\scriptfont\symfam=\fivesym
\def\sym{\fam\symfam\relax}
\def\N{{\sym N}}
\def\R{{\sym R}}

\def\Z{{\sym Z}}
\let\ds\displaystyle
\medskip
\centerline{\GGgras Decomposition of  deformations of thin rods. }

\centerline{\GGgras Application to nonlinear elasticity}
\bigskip
\centerline{\it Dominique Blanchard ${ }^a$ and Georges Griso ${ }^b$}

\centerline{\it ${ }^a$Universit\'e de Rouen, UMR 6085,  76801   Saint Etienne du Rouvray Cedex, France,}

\centerline{\it e-mail: dominique.blanchard@univ-rouen.fr}

\centerline{\it ${ }^b$ Laboratoire J.-Louis Lions, Universit\'e P.
et M. Curie, Case Courrier 187,}

\centerline{\it 75252 Paris C\'edex 05 - France,
e-mail:  griso@ann.jussieu.fr}

\bigskip
\noindent {\ggras Keywords: } {nonlinear elasticity, linear elasticity, $\Gamma$-convergence, rods, unfolding methods.}

\noindent{\ggras  2000 MSC: } 74B20, 74B05, 74K10, 49J45, 35A15.
\bigskip
\centerline{\Ggras Abstract. }
\medskip
This paper deals with the introduction of a decomposition of the deformations of curved  thin beams, with section of order  $\delta$, which takes into account the specific geometry of such beams. A deformation $v$ is split into an elementary deformation and a warping. The elementary deformation is the analog of a Bernoulli-Navier's displacement for linearized deformations replacing the infinitesimal rotation by a rotation in $SO(3)$ in each cross section of the rod. Each part of the decomposition is estimated with respect to the $L^2$ norm of the distance from gradient $v$ to $SO(3)$. This result relies on revisiting the  rigidity theorem of  Friesecke-James-MŸller in which we estimate the constant for a bounded open set star-shaped with respect to a ball. Then we use the decomposition of the deformations to derive a few asymptotic  geometrical behavior: large deformations of extensional type,  inextensional deformations  and linearized deformations. To illustrate the use  of our decomposition in nonlinear elasticity, we consider  a St Venant-Kirchhoff  material and upon various scaling on the applied forces we obtain the $\Gamma$-limit of the rescaled elastic energy. We first analyze the case of bending forces of order $\delta^2$ which leads to a nonlinear inextensional model. Smaller pure bending  forces give  the classical linearized model. A coupled extensional-bending model is obtained for a class of   forces of order $\delta^2$ in traction and of order $\delta^3$ in bending.
\bigskip
\noindent {\GGgras I. Introduction}
\medskip
This paper pertains to the field of modeling the deformations of a thin  structure who has a curved rod-like geometry with a few applications to elastic rods. Let us consider a curved rod of fixed length and with cross sections of small diameter of order $\delta$. Let us denote  by  $s_3$ the arc length of the middle line of the rod, by  $\Gn_1(s_3), \Gn_2(s_3)$ two normal vectors of this line and  the corresponding coordinates   by  $s=(s_1,s_2,s_3)$. In this setting, the aim of this paper is twofold. In a first result, we show that a deformation $v$ of such a rod can be decomposed as the sum of an elementary deformation and of a residual one as follows (see  (II.2.1)):
$$v(s)={\cal V} (s_3)+\GR(s_3)\big(s_1\Gn_1(s_3)+s_2\Gn_2(s_3)\big)+\overline{v}(s) .\leqno (I.1)$$
In the above decomposition, the field ${\cal V} (s_3)$ is the mean of $v$ over each section  and $\GR(s_3)\big(s_1\Gn_1(s_3)+s_2\Gn_2(s_3)\big)$ is the rotation of the same section, meaning that $\GR(s_3)\in SO(3)$ (the special orthogonal group i.e. the set of orthogonal $3\times 3$-matrices with determinant equal to $1$). The residual field 
$\overline{v}(s)$ represents the warping of a section. The main interest of our decomposition is the fact that each term is estimated with respect to  $\delta$ and  the $L^2$-norm of the distance between $\nabla v$ to $SO(3)$. In order to obtain such decompositions, we first adapt the proof  of the so called "Rigidity Theorem" established by Friesecke-James-MŸller in [11]. Our improvement only consists in evaluating the  dependence of the quantity which measure the distance from the gradient of a deformation (defined on an open set $\Omega$)  to $SO(3)$ in terms of two geometrical parameters characterizing $\Omega$ (see Theorem II.1.1).  As far as thin structures are concerned, the main interest of this result  is the possibility to slice the considered structure into small pieces for which the two geometrical parameters are uniformly controlled. This point is particularly helpful for a  curved rod with a variable curvature which is the case investigated in the present paper. This allows to define the elementary deformation as a continuous field and to derive  estimates on ${\cal V}$, $\GR$, $\overline {v}$ and on the distance between $\nabla v$ and $\GR$. These estimates first permit to identify a few known critical orders for the quantity $||\hbox{dist}(\nabla v,SO(3))||_{L^2}$ with respect to $\delta$ (see [15], [19], [20]). Then, we explicitly investigate two cases namely where  $||\hbox{dist}(\nabla v,SO(3))||_{L^2}$ is of order $\delta^2$ and $\delta^\kappa$ where $\kappa$ is a {\it real} number strictly greater than $2$. Let us emphasize that the decomposition of $v$ together with the estimates on 
${\cal V}$, $\GR$, $\overline {v}$ allow to identify the limit of the Green-Lagrange strain tensor in terms of the limit of the components of the decomposition of $v$. Moreover this decomposition of a deformation is, in some sense, stable with respect to the limit process with respect to $\delta$, which can be seen as a justification of this splitting of $v$.

The second type of results concerns the asymptotic behavior of the deformations of elastic rods when $\delta$ goes to $0$, assuming that the elastic energy is comparable to  $||\hbox{dist}(\nabla v,SO(3))||^2_{L^2}$ and more precisely for a Saint Venant-Kirchhoff's material. We consider an elastic rod submitted to dead forces (which are assumed to be volume forces  to simplify the  computations but this is not essential). We strongly use the decomposition  (I.1) to choose the scaling for the applied forces. In order to obtain an elastic energy of order  $\delta^{2\kappa}$ with $\kappa \ge 2$, we are led to split the forces into two types: order $\delta^{\kappa-1}$ for the loads with mean equal to $0$ over each cross section  and order $\delta^{\kappa}$ for general loads. We mainly investigate the cases $\kappa=2$ and $\kappa>2$.
Then we also use our decomposition to identify the limit energy through a $\Gamma$-convergence argument in both cases. Let us briefly summarize the obtained results. 

In the case  $\kappa=2$, we obtain a minimization problem which depends only on  the fields ${\cal V}$ and  $\GR$ (and indeed on the forces and  the boundary conditions of the 3D problem) which corresponds to  the nonlinear energy  for inextensible rods obtained in [15] and [19].  Moreover if the rod is clamped  on one (and only one) of its extremities, we show that this minimization problem is equivalent to an integro-differential problem for $\GR$ and that for small enough forces there is uniqueness of the solution.

In the case  $\kappa>2$, the limit minimization problem corresponds to the standard linear bending-torsion energy which is also obtained in the case $\kappa=3$ in [15] and [20]. 

We also examine a situation where the forces satisfy a specific geometrical assumption (which corresponds to pure traction-compression for a straight rod) but are of order $\delta^{\kappa-1}$ ($\kappa \ge 3$) and  nevertheless which leads to an elastic energy of order  $\delta^{2\kappa}$. We obtain a linear limit model for extensional displacement in the elastic 1D rod (with an elastic limit energy  already derived in the case of a straight rod and a 3D energy of order   $\delta^6$ in [15] and [20]). 

As a general reference on elasticity, we refer to [7] and [3].  A general introduction to the mathematical modeling of elastic  rod  models can be found in [2], [24], see also e.g. [1], [16]. For the justification of rods or plates models in nonlinear elasticity we refer  [1], [8], [12], [14], [15], [18], [19], [20], [21],  [22], [23]. For a general introduction of $\Gamma$-convergence we refer to [9]. The  rigidity  theorem and its applications to thin structures using $\Gamma$-convergence arguments can be found in [11], [12], [19], [20]. For the decomposition of the deformations in thin structures, we refer to  [13], [14] and for a few applications the junctions of multi-structures and homogenization to [4], [5], [6].

The paper is organized as follow. Section II is devoted to introduce the decomposition (I.1) of the deformations in a thin curved rod and to establish the estimates on   ${\cal V}$, $\GR$, $\overline {v}$.  In Section III, after rescaling the rod and the various fields with respect to $\delta$, we investigate the limit of the Green-St Venant tensor in the two cases  $||\hbox{dist}(\nabla v,SO(3))||_{L^2}\sim \delta^2$ and  $||\hbox{dist}(\nabla v,SO(3))||_{L^2}\sim \delta^\kappa$ for $\kappa>2$. In Section IV, we consider an elastic curved rod made of a St Venant- Kirchhoff 's material (see IV.1.9). After rescaling the applied forces, we identify the limit energy (as $\delta$ goes to $0$) through a $\Gamma$-convergence technique. The $\Gamma$-limit is a functional of   ${\cal V}$ and $\GR$ if $\kappa =2$ and of the displacement field  ${\cal U}$ and of an infinitesimal rotation field  ${\cal R}\land \vec x$ if $\kappa>2$. Then,  a specific choice of applied forces leads to a linear extentional model. Section V is devoted to  give an equivalent formulation of the limit minimization problem obtained in the nonlinear case $\kappa=2$ which leads to a partial uniqueness result. At least an  appendix  at the end of the paper details a few technical points concerning the interpolation between two rotations and a density result.

\medskip
\noindent {\GGgras II. Decomposition  of a deformation in a thin curved rod}\medskip
\noindent In this section, we derive a decomposition of the type I.1 for a deformation $v$ of a curved rod together with the estimates given in Theorem II.2.2. In order to obtain these results,  we first adapt the proof of the  ''Theorem of Geometric Rigidity''  established  in [11]. As mentioned in the introduction we essentially 
evaluate the dependence of the quantity which measure the distance from $\nabla v$ to $SO(3)$ in terms of two geometrical parameters characterizing the domain. This is the object of Subsection II.1. Then Subsection II.2 is devoted to establish the estimates on the terms of the decomposition of $v$ with respect to $||\hbox{dist}(\nabla v,SO(3))||_{L^2}$ (see Theorem II.2.2). The techniques are similar to the ones developed for small displacements in [13] and [14]. At least, in Subsection II.3 where the rod is assumed to be clamped at least on one of its extremities, we deduce estimates of $v$ and $\nabla v$ in terms of  $||\hbox{dist}(\nabla v,SO(3))||_{L^2}$.
\medskip
\noindent  {\Ggras II.1. Estimating the constant in the Theorem of Geometric Rigidity}\medskip
We equip the vector space $\GM_n$ of $n\times n$ matrices  with the Frobenius norm defined by
$$\GA=\big(a_{ij}\big)_{1\le i,j\le n},\qquad\qquad |||\GA|||=\sqrt{\sum_{i=1}^n\sum_{j=1}^n|a_{ij}|^2}.$$
Recall that an open set $\Omega$ of $\R^n$ is said to be star-shaped with respect to a ball $B(O;R_1)$ if for any $x\in B(O;R_1)$ and any $y\in \Omega$ the segment $[x,y]$ in included in $\Omega$.
\medskip
\noindent{\ggras  Theorem II.1.1. }{\it Let $\Omega$ be an open set of  $\R^n$ contained in the ball $B\big(O ;R\big)$ and star-shaped with respect to the ball $ B\big(O; R_1\big)$, $(0<R_1\le R)$. For any $v\in \big(H^1(\Omega)\big)^n$, there exist  $\GR\in SO(n)$ and $\Ga\in \R^n$   such that 
$$\left\{\eqalign{
&||\nabla v-\GR||_{(L^2(\Omega))^ {n\times n}}\le C||\hbox{dist}(\nabla v;SO(n))||_{L^2(\Omega) }, \cr
&|| v-\Ga-\GR\, x||_{(L^2(\Omega ))^n }\le CR||dist\big(\nabla v; SO(n)\big)||_{ L^2(\Omega) }, \cr}\right.\leqno(II.1.1)$$ where the constant $C$ depends only on $n$ and $\ds {R\over R_1}$. }

\noindent{\ggras  Proof of Theorem II.1.1. } The proof of  the first inequality in Theorem II.1.1 is identical to the proof of Theorem 3.1 in [11] if we show that   the constants  which appear in the three main points of this proof only  depend upon $n$ and  $\ds {R\over R_1}$. These three main arguments are first an approximation lemma, then a specific covering of $\Omega$ and finally a Poincar\'e-Wirtinger's type inequality. In particular, we explicitly construct a covering of $\Omega$ which can be used in the proof of Theorem 3.1 of [11] and which only depends of $\ds {R\over R_1}$ and $n$.

We begin with the following lemma which just specifies the dependence of the constants in Proposition A.1 of [11].
\smallskip

\noindent{\ggras  Lemma II.1.2.  }{\it Let  $n\ge1$ be an integer  and  $1\le p<\infty$ be a real number.  Let $\Omega$ be an open set of  $\R^n$ contained in the ball $B\big(O ;R\big)$ and star-shaped with respect to the ball $ B\big(O; R_1\big)$, $(0<R_1\le R)$. There exists a constant  $C=C(n,p,R/R_1)$ such that for any function  $v\in \big(W^{1,p}(\Omega)\big)^n $ and for any real number $\lambda>0$ there exists a function  $w\in \big(W^{1,\infty}(\Omega)\big)^n $ such that
$$\eqalign{
(i)&\qquad ||\nabla w||_{(L^\infty(\Omega))^{n\times n}}\le C \lambda\cr
(ii)&\qquad \big|\big\{x\in \Omega\enskip;\enskip v(x)\not=w(x)\big\}\big|\le 
{C \over \lambda^p}\int_{\{x\in \Omega\,;\, |||\nabla v(x)||| >
\lambda\}}|||\nabla v(x)|||^pdx\cr 
(iii)&\qquad||\nabla v-\nabla w||_{(L^p(\Omega))^{n\times n}}\le C \int_{\{x\in \Omega\,;\,
|||\nabla v(x)|||> \lambda \}}|||\nabla v(x)|||^pdx}$$}
\noindent{\ggras  Proof of Lemma II.1.2.  } \noindent Let us denote by $B_n=B\big(O;1\big)$ the unit ball of $\R^n$ and set  $S_n=\partial B_n$. The proof of  Lemma II.1.2 is given in [11] except what concerns the dependence of the constant in the inequalities with respect to the geometrical parameter $R/R_1$ which will be extensively used in the sequel. We recall the following result proved in [11] ( Proposition A1; see also Evans and Gariepy [10], Section 6.6.2 and 6.6.3): there exists a constant $C_0$ which depends on  $n$ and   $p$ such that for any function $\widetilde{v}$ in $\big(W^{1,p}(B_n )\big)^n $ and for any real number $\tilde\lambda>0$ there exists a function  $\widetilde{w}\in \big(W^{1,\infty}(B_n )\big)^n $ such that 
$$\left\{\eqalign{
& ||\nabla_y\widetilde{w}||_{(L^\infty(B_n))^{n\times n}}\le C_0 \tilde \lambda\cr
& \big|\big\{y\in B_n\enskip;\enskip \widetilde{v}(y)\not=\widetilde{w}(y)\big\}\big|\le 
{C_0 \over\tilde\lambda^p}\int_{\{y\in B_n\,;\, |||\nabla_y\widetilde{v}(y)||| >\tilde
\lambda\}}|||\nabla_y\widetilde{v}(y)|||^pdy\cr 
&||\nabla_y\widetilde{v}-\nabla_y\widetilde{w}||_{(L^p(B_n))^{n\times n}}\le C_0 \int_{\{y\in B_n\,;\,
|||\nabla_y\widetilde{v}(y)|||>\tilde \lambda \}}|||\nabla_y\widetilde{v}(y)|||^pdy\cr }\right.\leqno (II.1.2)$$ 
Since  $\Omega$ is, in particular,  star-shaped with respect to the origin  $O$, for any  direction $s$ of $S_n$ the ray issued from $O$ and with direction $s$ meets the boundary $\partial \Omega$ on a unique point $P(s)$.  In order to transform the ball $B_n$ into the set $\Omega$, we first introduce the function $F$ from $S_n$ into $\R^+$ by
$$\forall s\in S_n,\qquad\qquad F(s)=\Vert\overrightarrow{OP}(s)\Vert_2,$$ 
\noindent where $\Vert.\Vert_2$ denotes the euclidian norm on $\R^n$

\noindent Now the function  $G$ from $\R^n$ into $\R^n$ is defined by 
$$G\; : \; y\in \R^n\longrightarrow \left\{\eqalign{ & yF\Bigl({y\over \Vert y\Vert_2}\Bigr) \; \hbox{if}\; y\not=O\cr
&O\qquad\qquad \hbox{if}\;  y=O.\cr}\right.$$ This function $G$ is one to one  from $\R^n$ onto $\R^n$ and maps $B_n$ onto  $\Omega$. Moreover, due to the geometrical assumptions on $\Omega$,  the function $G$ is Lipschitz-continuous and satisfies the following inequalities for  almost any $y\in \R^n$
$$ R_1C_1\le |||\nabla_y G(y)|||\le RC_2,\quad{C_1\over R}\le |||\nabla_x G^{-1}(x)|||\le {C_2\over R_1},\quad R_1^nC_1\le |\det\big(\nabla_y G(y)\big)|\le R^nC_2\leqno (II.1.3)$$  
\noindent where the constants $C_1$ and $C_2$ depend on $n$ and  $R/R_1$. The proof of the above estimates is left to the reader (see also [13] and [14]).

\noindent Let  $v\in \big(W^{1,p}(\Omega))^n$. We define the function
$\widetilde{v}=v\circ G$ which belongs to $ \big(W^{1,p}(B_n)\big)^n $ and we have for almost any  $y\in B_n$ 
$$\nabla_y\widetilde{v}(y)=\nabla_x v\big(G(y)\big)\nabla_yG(y).$$ Taking into account (II.1.3), we obtain
$$R_1C_3|||\nabla_x v\big(G(y)\big)|||\le |||\nabla_y\widetilde{v}(y)|||\le RC_4|||\nabla_x v\big(G(y)\big)|||\leqno(II.1.4)$$ 
\noindent where $C_3$ and $C_4$ depend on $n$ and $R/R_1$. Using the estimates on the jacobian given by (II.1.3), we deduce that 
$$\left\{\eqalign{
C_5  {R_1^p\over R^n} ||\nabla_x v||^p_{(L^p(\Omega))^{n\times n}} & \le ||\nabla_y\widetilde{v}||^p_{(L^p(B_n))^{n\times n}}\le C_6{R^p\over R_1^n} ||\nabla_x v||^p_{(L^p(\Omega))^{n\times n}}\qquad\hbox{for}\;\; 1\le p<\infty\cr  
C_5 R_1||\nabla_x v||_{(L^\infty(\Omega))^{n\times n}} &\le
||\nabla_y\widetilde{v}||_{(L^\infty(B_n))^{n\times n}}\le  C_6 R||\nabla_x v||_{(L^\infty(\Omega))^{n\times n}}\cr}\right.\leqno(II.1.5)$$  where $C_5$ and $C_6$ depend on $n$ and $R/R_1$. Now we apply the result recalled at the beginning of the proof  so that for any $\lambda>0$, setting $\tilde \lambda=C_4R \lambda$,  there exists a function $\widetilde{w}\in \big(W^{1,\infty}(B_n )\big)^n $ such that (II.1.2) holds true.
Let us set  $w=\widetilde{w}\circ G^{-1}$ which belongs to  $\big(W^{1,\infty}(\Omega)\big)^n $. Thanks to (II.1.2) and (II.1.5)    we have
$$||\nabla_x w||_{(L^\infty(\Omega))^{n\times n}}\le {C_0  \tilde \lambda\over C_5R_1}={C_0 C_4 \over C_5 }{ R\over  R_1}\lambda,$$ and   $i)$ is proved. We use   (II.1.3) and (II.1.4) to obtain
$$\eqalign{
&\big|\big\{x\in \Omega\enskip;\enskip v(x)\not=w(x)\big\}\big|\le C_2R^n\big|\big\{y\in
B_n\enskip;\enskip
\widetilde{v}(y)\not=\widetilde{w}(y)\big\}\big|\cr
\le &{C_0C_2R^n \over\tilde\lambda^p}\int_{\{y\in B_n\,;\,|||\nabla_y\widetilde{v}(y)|||>\tilde \lambda
\}}|||\nabla_y\widetilde{v}(y)|||^pdy\le {C_0 C_2\over C_1}{R^n\over R_1^n}{1\over \lambda^p}\int_{\{x\in \Omega\,;\,|||\nabla_x v(x)|||>  \lambda   \}}|||\nabla_x v(x)|||^pdx\cr
}$$  and  $ii)$ is established. Now we prove $iii)$. We have for $\lambda$ and $w$ satisfying i) and ii)
$$\eqalign{
\int_\Omega|||\nabla_x v(x)-\nabla_x w(x)|||^pdx&=\int_{v\not=w}|||\nabla_x v(x)-\nabla_x w(x)|||^pdx\le
2^p\int_{v\not=w}\big\{|||\nabla_x v(x)|||^p+|||\nabla_x w(x)|||^p\big\}dx\cr
&\le 2^p\int_{v\not=w}\big\{\lambda^pdx+|||\nabla_x w(x)|||^p\big\}dx+2^p\int_{|||\nabla v(x)|||>\lambda}|||\nabla_x v(x)|||^pdx\cr
&\le C\int_{v\not=w} \lambda^pdx+ 2^p\int_{|||\nabla_x v(x)|||>\lambda} |||\nabla_x v(x)|||^pdx\le  C\int_{|||\nabla_x v(x)|||>\lambda}|||\nabla_x v(x)|||^pdx\cr}$$ Finally we obtain
$$ ||\nabla_x v-\nabla_xw||_{(L^p(\Omega))^{n\times n}}\le C \int_{\{x\in \Omega\,;\,
|||\nabla_xv(x)|||> \lambda \}}|||\nabla_xv(x)|||^pdx,$$
where the constant  depends on $n$, $p$ and $R/R_1$. This concludes the proof of Lemma II.1.2.\fin\medskip
\medskip
We now turn to the second argument in the proof of Theorem II.1.1, namely the specific covering of $\Omega$. In the following we construct a covering ${\cal Q}$ of  $\Omega$ with cubes of the type $Q(a,r)=a+]-r,r[^n$, $r>0$ satisfying the following properties:

* for every $Q(a,r)\in {\cal Q}$, the cube $Q(a,2r)$ is included in $\Omega$,

* there exists a constant $\ds C\big(n,{R\over R_1}\big)$ such that
$$\forall x\in Q(a,r)\in{\cal Q},\qquad  r\le dist_\infty(x,\partial \Omega)\le C\big(n,{R\over R_1}\big) r,\leqno (II.1.6)$$

* there exists a finite integer $N(n)$ such that for every cube $Q(a,r)$ of the covering  ${\cal Q}$, the number of cubes of the type $Q(b,2r')$ which meet $Q(a,r)$, where $Q(b,r')$ belongs to  ${\cal Q}$, is at most  $N(n)$.

For $i\in\N$ and $k\in \N$, we set  $r_i =\ds{R_1\over 2^ i3\sqrt n}$  and  $${\cal R}_k=\Bigl\{a\in r_k\Z^n\; |\; a=r_k\,(i_1, i_2,\ldots, i_n),\quad i_p \in \Z \;
\hbox{and odd}\Big\}.$$
\noindent The covering ${\cal Q}$ of  $\Omega$ is constructed by induction  as follows :

\noindent * consider all the cubes  $Q(a, r_0 )$, $a\in {\cal R}_0 $, such that $Q(a,2r_0 )\subset \Omega$ and denote by  ${\cal Q}_0 $ the family of these cubes $Q(a, r_0)$ and by  ${\cal U}_0 =\ds\bigcup_{Q(a, r_0 )\in{\cal Q}_0 }Q(a, r_0 )$,

\noindent * in step $k\ge 1$, consider the cubes  $Q(a, r_k)$, $a\in {\cal R}_k$, such that
$Q(a, r_k)\subset \Omega\setminus \overline{{\cal U}_0 \cup\ldots \cup{\cal U}_{k-1}}$ and such that $Q(a, 2r_k)\subset \Omega$, and denote by  ${\cal
Q}_k$   the family of these cubes  $Q(a, r_k)$ and by ${\cal U}_k=\ds\bigcup_{Q(a,r_k)\in {\cal Q}_k}Q(a, r_k)$. 

\noindent We denote by  ${\cal Q}$ the countable family of all the cubes constructed through the above process. 
\smallskip
\noindent 
The above explicit construction  permits to show that the covering  ${\cal Q}$ verifies the required properties (as an example we can take $\ds C\big(n,{R\over R_1}\big)=5\sqrt n{R\over R_1} $ and $N(n)=2^{n+3}$).

\medskip
As far as the third argument in the proof of Theorem I.2.1 is concerned, we now recall  the following Poincar\'e-Wirtinger's  inequality (see [13] for a proof and various applications). Since $\Omega$ is contained in the ball $B\big(O ;R\big)$ and is star-shaped with respect to the ball $ B\big(O; R_1\big)$, there exists a constant $C$ which depends on $n$ and $\ds {R\over R_1}$ such that for any $\phi\in H^1(\Omega)$ (see [13])
$$|| \phi-{\cal M}(\phi)||_{L^2(\Omega)}\le C ||\rho \nabla \phi||_{(L^2(\Omega))^n}\leqno (II.1.7)$$
where ${\cal M}(\phi)$ is the mean of $\phi $ over $\Omega$ and $\rho(x)=\hbox{dist}(x,\partial \Omega)$. 
Using Lemma II.1.2, the specific covering of $\Omega$ described above and the Poincar\'e-Wirtinger's inequality  (II.1.7) permit to reproduce the proof of Theorem 3.1 in [11] in order to obain the first estimate of  our  Theorem II.1.1 with a constant which depends only on $n$ and $\ds {R\over R_1}$. To end the proof, we apply inequality (II.1.7) to the field $ v(x)-\GR x$ and we use the first estimate in Theorem II.1.1.
\fin
\smallskip
\noindent  {\Ggras II.2. Decomposition of the deformation in a curved rod. Estimates
 }\medskip
\noindent {\ggras II.2.1. The geometry}
\smallskip
Let us introduce a few  notations and definitions concerning the geometry of a curved rod (see [13], [14] for a detailed presentation).

Let $\zeta$ be a  curve in the euclidian space
$\R^3$ parametrized by  its  arc length $s_3$. The current point of the curve  is denoted $M(s_3)$.

\noindent We suppose that the mapping $s_3\longrightarrow M(s_3)$  belongs to $\big({\cal C}^2([0,L])\big)^3 $ and that it is  one to one.  We have
$${ {dM}\over ds_3}={\Gt},\qquad\qquad \|{\Gt}\|_2=1,$$ where $\|\cdot\|_2$ is the euclidian norm in $\R^3$. 

Let $\Gn_1$ be a function belonging to $\big({\cal C}^1([0,L])\big)^3 $ and such that
$$\forall s_3\in [0,L],\qquad \|\Gn_1(s_3)\|_2=1\qquad \hbox{and}\qquad \Gt(s_3)\cdot\Gn_1(s_3)=0.$$
We set $${\Gn_2}={\Gt}\land {\Gn_1}.$$

In the sequel, $\omega$ denotes a bounded domain in $\R^2$  with lipschitzian boundary  (while obviously, 
$\overline{\omega}$  denotes the closure of $\omega$).  We choose the origin O of coordinates at the center of mass of
$\omega$ and we choose the coordinates axes  $(O;\Ge_1)$ and $(O;\Ge_2)$ as the principal axes of inertia of $\omega$, so that $\ds\int_\omega x_1dx_1dx_2=\int_\omega x_2dx_1dx_2=\int_\omega x_1x_2dx_1dx_2=0$.
The reference cross-section $\omega_\delta$ of the rod  is obtained by transforming $\omega$ with a dilatation of ratio $\delta>0$ and we set 
$$ \Omega_\delta= \omega_\delta\times(0,L).$$
\medskip
 Introduce  now  the mapping  $\Phi:  \R^2\times [0,L]\longrightarrow \R^3$ defined by 
$$\Phi\;\; :
\;(s_1,s_2,s_3)\longmapsto M(s_3)+s_1 {\Gn_1}(s_3)+s_2 {\Gn_2}(s_3)$$  \noindent There exists  $\delta_0 >0$  depending only on $\zeta$, such that  the restriction of  $\Phi$ to the compact set
$\overline{\Omega}_{\delta }$ is a ${\cal C}^1-$ diffeomorphism between $\overline{\Omega}_{\delta }$ and $\Phi(\overline{\Omega}_{\delta })$. Moreover, there exists two  positive constants $c $ and $c_1$    such that
$$\forall \delta \in (0,\delta_0 ],\qquad \forall s\in \overline{\Omega}_{\delta},\qquad\qquad c \le |||\nabla\Phi(s)|||\le c_1.$$ 
\noindent{\ggras Definition II.2.1. }{\it For $\delta\in
(0,\delta_0 ]$, the curved rod ${\cal P}_\delta$ is defined by 
$${\cal P}_\delta=\Phi \big(\Omega_\delta\big).$$} The cross-section of the curved rod is
isometric to $\omega_\delta$. In ${\cal P}_\delta$, the point $M(s_3)$ is the center of gravity of the cross-section  $\Phi(\omega_\delta\times\{s_3\})$ and the axes of direction  $\Gn_1(s_3)$ and
$\Gn_2(s_3)$ are the principal axes of this cross-section.
\smallskip
\noindent{\ggras Notation. }  {\sl  Reference domains and running points}.   We denote $x$ and $s$ respectively the running point of
${\cal P}_\delta$ and of  $\Omega_\delta$ so that  $x=\Phi(s)$. 

A deformation $v$ defined on ${\cal P}_\delta$ can be also considered as  a deformation defined on $\Omega_\delta$ which  we will also denote by $v$, as a convention. As far as the gradients of $v$ are concerned we have  $\nabla_s v=\nabla_x v.\nabla \Phi$ 
$$\forall \delta \in (0,\delta_0 ],\qquad c|||\nabla_x v(x) |||\le |||\nabla_s v(s)|||\le C|||\nabla_x v(x) |||$$
where the constants are  positive and do not depend on $\delta$.
\medskip
\noindent {\ggras II.2.2. The elementary deformation}
\smallskip

In this subsection, we show that any deformation $v\in \big(H^1({\cal P}_\delta)\big)^3$ of the rod ${\cal P}_\delta$ can be decomposed as (using the above convention)
$$ v(s)={\cal V} (s_3)+\GR(s_3)\big(s_1\Gn_1(s_3)+s_2\Gn_2(s_3)\big)+\overline{v}(s),\qquad s\in \Omega_\delta,\leqno (II.2.1)$$
where ${\cal V} $ belongs to $\big(H^1(0,L)\big)^3$, $\GR$  belongs to $\big(H^1(0,L)\big)^{3\times 3}$ and satisfies for any $s_3\in [0,L]$: $\GR(s_3)\in SO(3)$ and $\overline{v}$ belongs to $\big(H^1({\cal P}_\delta)\big)^3 $ (or $\big(H^1(\Omega_\delta)\big)^3 $  using  again the same convention as for $v$). Let us give a few comments on the above decomposition. The term ${\cal V}$ gives the deformation of the center line of the rod and it is indeed a function of  the arc length $s_3$. The second term $\GR(s_3)\big(s_1\Gn_1(s_3)+s_2\Gn_2(s_3)\big)$ describes the rotation of the cross section (of the curved rod) which contains the point $M(s_3)$. The sum of the  two first terms ${\cal V} (s_3)+\GR(s_3)\big(s_1\Gn_1(s_3)+s_2\Gn_2(s_3)\big)$ is called an elementary deformation of the rod.
\medskip
\noindent  {\ggras II.2.3. The main theorem}
\medskip
The following theorem gives a decomposition (II.2.1) of a deformation and estimates on the terms of this decomposition.
\smallskip
\noindent{\ggras  Theorem II.2.2. }{\it Let $v\in \big(H^1({\cal P}_\delta)\big)^3 $, there exists an elementary deformation ${\cal V} +\GR \big(s_1\Gn_1 +s_2\Gn_2 \big)$ and a warping $\overline{v}$ satisfying (II.2.1) and such that 
$$\left\{\eqalign{
&||\overline{v}||_{(L^2(\Omega_\delta))^3}\le C\delta ||\hbox{dist}(\nabla_x v,SO(3))||_{L^2({\cal P}_\delta)}\cr
&||\nabla_s\overline{v}||_{(L^2(\Omega_\delta))^{3\times 3}}\le C ||\hbox{dist}(\nabla_x v,SO(3))||_{L^2({\cal P}_\delta)}\cr
&\Bigl\|{d\GR\over ds_3}\Big\|_{(L^2(0,L))^{3\times 3}}\le {C\over \delta^2} ||\hbox{dist}(\nabla_x v,SO(3))||_{L^2({\cal P}_\delta)}\cr
& \Bigl\|{d{\cal V}\over ds_3}-\GR \Gt\Big\|_{(L^2(0,L))^3}\le {C\over \delta}||\hbox{dist}(\nabla_x v,SO(3))||_{L^2({\cal P}_\delta)}\cr
& \bigl\|\nabla_x v-\GR \big\|_{(L^2(\Omega_\delta))^{3\times 3}}\le C||\hbox{dist}(\nabla_x v,SO(3))||_{L^2({\cal P}_\delta)}\cr}\right.\leqno(II.2.2)$$ where the constant $C$ does not depend on $\delta$.}
\medskip
\noindent{\ggras Proof. }  Let $N$ be an integer belonging to  $\displaystyle  \Big[ {2L\over 3\delta}, {L\over \delta}\Big]$  and  let 
$\ds 0 \le \alpha \le L- {L\over N}$. 

\noindent We have $\displaystyle \delta\le{L\over N}\le{3\over 2}\delta$. Let $R>1$  be such that the reference cross-section  $\omega$ is contained in the ball $B(O;R)$. Then the domain $\ds\Omega_{\delta, \alpha}=\omega_\delta\times ]\alpha, \alpha+ {L\over N}[$ has a diameter less than $3R\delta$. In the sequel we will work with the portions ${\cal P}_{\delta, \alpha}$ of the rod ${\cal P}_\delta$ defined by $${\cal P}_{\delta,\alpha}=\Phi\bigl(\Omega_{\delta, \alpha} \bigr).$$ As in [13] we distinguish two cases.   
 
 \noindent{\ggras  First case. } If $\omega$ is star-shaped with respect to a ball of radius  $R_1\le1/2$, it is shown in [13] that each ${\cal P}_{\delta,\alpha}$ is star-shaped with respect to a ball of radius  $\ds{R_1\delta \over 8}$ and we are in a position to apply Theorem II.1.1 to the function $v$ into each part ${\cal P}_{\delta,\alpha}$ for which the ratio $\ds{R\over R_1}$ is independent of $\delta$.  As a consequence, there exist   $\GR_\alpha\in SO(3)$ and $\Ga_\alpha\in \R^3$   such that 
$$\left\{\eqalign{
&||\nabla_x v-\GR_\alpha||_{(L^2({\cal P}_{\delta,\alpha}))^ {3\times 3}}\le C||\hbox{dist}(\nabla_x v;SO(3))||_{L^2({\cal P}_{\delta,\alpha})}  \cr
&|| v-\Ga_\alpha-\GR_\alpha\, \big(x-M(\alpha)\big)||_{(L^2({\cal P}_{\delta,\alpha}))^3 }\le C\delta||dist\big(\nabla_x v; SO(3)\big)||_{L^2({\cal P}_{\delta,\alpha})}. \cr}\right.\leqno (II.2.3)$$ The constant $C$ does not depend on $\alpha$ and $\delta$.

 \noindent{\ggras  Second Case. }  Let us consider the general case, where the cross-section is a bounded domain in $\R^2$ with lipschitzian boundary. There
exists a finite sequence of  open sets $\omega^{(1)}, \ldots, \omega^{(K)}$ such that
$$\omega=\bigcup_{1\le l\le K}\omega^{(l)},\qquad\qquad \omega_\delta=\bigcup_{1\le l\le K}\omega^{(l)}_\delta, $$
and such that every $\omega^{(l)}$  is star-shaped with respect to a disc of radius $R_1$, $0<R_1\le 1/2$.  Moreover, the
open set $\omega$ is connected, then there exists $R_2\in ]0,R_1]$ such that  $\omega^{(r)}\cap \omega^{(s)}$ contains  a disc of radius $R_2$ if the intersection is not empty.

The domain $\ds \Omega^{(l)}_{\delta,\alpha}=\omega^{(l)}_\delta\times ]\alpha,\alpha+{L\over N}[$ is star-shaped with respect to  a
ball  of radius $R_1\delta$. As in the first case, there exist   $\GR^{(l)}_\alpha\in SO(3)$ and $\Ga^{(l)}_\alpha\in \R^3$   such that 
$$\eqalign{
&||\nabla_x v-\GR^{(l)}_\alpha||_{(L^2({\cal P}^{(l)}_{\delta,\alpha}))^ {3\times 3}}\le C||\hbox{dist}(\nabla_x v;SO(3))||_{L^2({\cal P}^{(l)}_{\delta,\alpha})}  \cr
&|| v-\Ga^{(l)}_\alpha-\GR^{(l)}_\alpha\, \big(x-M(\alpha)\big)||_{(L^2({\cal P}^{(l)}_{\delta,\alpha}))^3 }\le C\delta||dist\big(\nabla_x v; SO(3)\big)||_{L^2({\cal P}^{(l)}_{\delta,\alpha})}. \cr}$$ The constant $C$ does not depend on $\alpha$ , $\delta$ and $l$.
 
If $\omega^{(r)}\cap \omega^{(s)}\not =\emptyset$
the portion ${\cal P}^{(r)}_{\delta,\alpha}\cap {\cal P}^{(s)}_{\delta,\alpha}$ contains a ball of radius $R_2\delta/8$. This allows us
to compare the elementary deformations  $\Ga^{(r)}_\alpha+\GR^{(r)}_\alpha\,  \big(x-M(\alpha)\big)$ and $\Ga^{(s)}_\alpha+\GR^{(s)}_\alpha\,  \big(x-M(\alpha)\big)$ in this ball.  We obtain
$$|||\GR^{(r)}_\alpha-\GR^{(s)}_\alpha|||^2\le {C\over \delta^3}||dist\big(\nabla_x v; SO(3)\big)||^2_{L^2(
{\cal P}^{(r)}_{\delta,\alpha}\cup{\cal P}^{(s)}_{\delta,\alpha})}\qquad||\Ga^{(r)}_\alpha-\Ga^{(s)}_\alpha||^2\le {C\over \delta}||dist\big(\nabla_x v; SO(3)\big)||^2_{L^2(
{\cal P}^{(r)}_{\delta,\alpha}\cup{\cal P}^{(s)}_{\delta,\alpha})}$$ where the constant only depends on $R$ , $R_1$ and $R_2$.

Setting $\GR_\alpha=\GR^{(1)}_\alpha$ and $\Ga_\alpha=\Ga^{(1)}_\alpha$ and proceeding step by step with respect to $l$, we finally deduce that (II.2.3) holds true with  a constant $C$ which does not depend on $\alpha$ and $\delta$ as in the first case.

Now we consider two splittings of ${\cal P}_\delta$ by considering two sets of arc length 
$$\alpha_k=k{L\over N}, \quad k=0,\ldots,N,\qquad \beta _k=\alpha_k+{L\over2N}, \quad k=0,\ldots,N-1.$$

We  consider the elementary  deformations $\Ga_{\alpha_k}+\GR_{\alpha_k}\, \big(x-M(\alpha_k)\big)$ and $\Ga_{\beta_k}+\GR_{\beta_k}\, \big(x-M(\beta_k)\big)$ of the portions ${\cal P}_{\delta,\alpha_k}$ and ${\cal P}_{\delta,\beta_k}$ which satisfies estimates (II.2.3) with a constant independent of $k$.  Considering ${\cal P}_{\delta,\alpha_k}\cap {\cal P}_{\delta,\beta_k}$ and ${\cal P}_{\delta,\alpha_{k+1}}\cap {\cal P}_{\delta,\beta_k}$, we can compare $\GR_{\alpha_k}$ and $\GR_{\alpha_{k+1}}$. We obtain 
$$\sum_{k=0}^{N-1}{L\over N}\Big|\Big|\Big|{\GR_{\alpha_{k+1}}-\GR_{\alpha_k}\over L/N}\Big|\Big|\Big|^2\le {C\over \delta^4}||dist\big(\nabla_x v; SO(3)\big)||^2_{L^2({\cal P}_{\delta})},\leqno (II.2.4)$$
where $\GR_{\alpha_N}=\GR_{\alpha_{N-1}}$ and where the constant $C$ is independent of $\delta$.

Now we are in a position to define the elementary deformation associated to $v$. We set for $s_3\in [0,L]$ 
$${\cal V}(s_3)={1\over |\omega_\delta|}\int_{\omega_\delta}v(s_1,s_2,s_3) ds_1ds_2.\leqno (II.2.5)$$
In order to define $\GR$ 
 we use the following  argument whose proof is postponed to the appendix. There exists a field of  matrices  $\GR$ belonging to $\big(H^1(0,L)\big)^{3\times 3}$, with $\GR(s_3)\in SO(3)$ for all $s_3\in [0,L]$, such that $\GR(\alpha_k)=\GR_{\alpha_k}$ for  $k=0,\ldots,N$ and $$\Big\|{d\GR\over ds_3}\Big\|^2_{(L^2(0,L))^{3\times3}}\le 4 \sum_{k=0}^{N-1}{L\over N}\Big|\Big|\Big|{\GR_{\alpha_{k+1}}-\GR_{\alpha_k}\over L/N}\Big|\Big|\Big|^2. \leqno(II.2.6)$$
Indeed the field $\overline{v}$ is defined by $$ \overline{v}(s)=v(s)-{\cal V}(s_3)-\GR(s_3)\big(s_1\Gn_1+s_2\Gn_2\big)\quad \hbox{for a. e. }  s\in \Omega_\delta.\leqno(II.2.7)$$ The third estimate of (II.2.2) follows directly from (II.2.4)  and (II.2.6).
 \smallskip
 From (II.2.3) we obtain
$$\left\{\eqalign{
&\sum_{k=0}^{N-1}||\nabla_x v-\GR_{\alpha_k}||^2_{(L^2({\cal P}_{\delta,\alpha_k}))^ {3\times 3}}\le C||\hbox{dist}(\nabla_x v;SO(3))||^2_{L^2({\cal P}_{\delta})}  \cr
&\sum_{k=0}^{N-1}|| v-\Ga_{\alpha_k}-\GR_{\alpha_k}\, \big(x-M(\alpha_k)\big)||^2_{(L^2({\cal P}_{\delta,\alpha_k}))^3 }\le C\delta^2||dist\big(\nabla_x v; SO(3)\big)||^2_{L^2({\cal P}_{\delta})}. \cr}\right.\leqno(II.2.8)$$
Now, we take  the mean  value over the cross-sections of $\Omega_\delta$, then using the definition  of ${\cal V}$  we deduce  
$$\sum_{k=0}^{N-1}\bigl\|{\cal V}-\Ga_{\alpha_k}-\GR_{\alpha_k} \big(M-M(\alpha_k)\big)\bigr\|^2_{(L^2(\alpha_k,\alpha_{k+1}))^3}\le C||dist\big(\nabla_x v; SO(3)\big)||^2_{L^2({\cal P}_{\delta})}.\leqno (II.2.9)$$ Thanks to the definition of $\GR$ and the third estimate in (II.2.2) we get
$$\sum_{k=0}^{N-1}\bigl\|\GR-\GR_{\alpha_k} \bigr\|^2_{(L^2(\alpha_k,\alpha_{k+1}))^{3\times 3}}\le {C\over \delta^2}||dist\big(\nabla_x v; SO(3)\big)||^2_{L^2({\cal P}_{\delta})}.\leqno (II.2.10)$$ 
Consequently (II.2.9) and  (II.2.10)  give the first estimate in (II.2.2) while  (II.2.10) leads to the last estimate in (II.2.2).

Due to the definition of $\Phi$, we have $\ds {\partial \Phi\over \partial s_\alpha}=\Gn_\alpha$ and $\ds {\partial \Phi\over \partial s_3}=\Gt+s_1{d\Gn_1\over ds_3}+s_2{d \Gn_2\over ds_3}$ so that the relation  $\nabla_s v=\nabla_x v .\nabla \Phi$ leads to 
$${\partial v\over \partial s_1}=\nabla_x v\;\Gn_1\qquad {\partial v\over \partial s_2}=\nabla_x v\;\Gn_2\qquad {\partial v\over \partial s_3}=\nabla_x v\Big(\Gt+s_1{d \Gn_1\over ds_3}+s_2{d\Gn_2\over ds_3}\Big).\leqno(II.2.11)$$
Then inserting (II.2.10) into (II.2.8)  gives in particular
$$ \Big\|{\partial v\over \partial s_3}-\GR\big(\Gt+s_1{d \Gn_1\over ds_3}+s_2{d\Gn_2\over ds_3}\big)\Big\|^2_{(L^2(\Omega_\delta))^3}\le  C||dist\big(\nabla_x v; SO(3)\big)||^2_{L^2({\cal P}_{\delta})}.\leqno (II.2.12)$$ 

Integrating first over  $\omega_\delta\times \{s_3\}$ in  (II.2.12) leads to the forth estimate of (II.2.2) (recall that $\ds \int_\omega s_1 ds_1ds_2=\int_\omega s_2 ds_1ds_2=0$). It remains to show the estimate on $\nabla_s \overline{v}$. From (II.2.11) we get
$${\partial \overline{v}\over \partial s_1}=(\nabla_x v-\GR)\;\Gn_1\qquad {\partial \overline{v}\over \partial s_2}=(\nabla_x v-\GR)\;\Gn_2\leqno(II.2.13)$$
and then with (II.2.8) and (II.2.10) 
$$ \Big\|{\partial \overline{v}\over \partial s_1}\Big\|^2_{(L^2(\Omega_\delta))^3}+ \Big\|{\partial \overline{v}\over \partial s_2}\Big\|^2_{(L^2(\Omega_\delta))^3}\le  C||dist\big(\nabla_x v; SO(3)\big)||^2_{L^2({\cal P}_{\delta})}.\leqno (II.2.14)$$ 
Now we estimate $\ds\Big\|{\partial \overline{v}\over \partial s_3}\Big\|_{(L^2(\Omega_\delta))^3}$. We have from (II.2.7) and (II.2.11)
$${\partial \overline{v}\over \partial s_3}=\big(\nabla_x v-\GR\big)\big(\Gt+s_1{d \Gn_1\over ds_3}+s_2{d\Gn_2\over ds_3}\big)-{d\GR\over ds_3}\big(s_1 \Gn_1 +s_2 \Gn_2\big)-{d{\cal V}\over ds_3}+\GR\;\Gt.\leqno(II.2.15)$$
Proceeding as above to bound the first term and using the third, the forth  and the last estimates of (II.2.2) to control the last two terms of the above inequality permit to obtain the estimate on $\ds{\partial \overline{v}\over \partial s_3}$ given in (II.2.2). \fin
\medskip
In order  to split the bending and the stretching of the rod, we now introduce the following quantities:
$$\forall s_3\in [0,L],\qquad {\cal V}_B(s_3)={\cal V}(0)+\int_0^{s_3}\GR(z)\Gt(z)dz,\qquad {\cal V}_S(s_3)={\cal V}(s_3)-{\cal V}_B(s_3).\leqno(II.2.16)$$ Let us notice that ${\cal V}_B$ is the bending deformation of the middle line while ${\cal V}_S$ is the stretching deformation. Due to the third estimate in (II.2.2) we have
$$ \bigl\|{\cal V}_S\big\|_{(H^1(0,L))^3}\le {C\over \delta}||\hbox{dist}(\nabla_x v,SO(3))||_{L^2({\cal P}_\delta)}.\leqno(II.2.17)$$ Inserting the definition (II.2.16) into the decomposition (II.2.1) gives
$$ v(s)={\cal V}_B (s_3)+\GR(s_3)\big(s_1\Gn_1(s_3)+s_2\Gn_2(s_3)\big)+{\cal V}_S(s_3)+\overline{v}(s),\qquad s\in \Omega_\delta.\leqno(II.2.18)$$

Estimates (II.2.2) and  (II.2.17) permit to interpret the part ${\cal V}_B (s_3)+\GR(s_3)\big(s_1\Gn_1(s_3)+s_2\Gn_2(s_3)\big)$ of the decomposition (II.2.18) as an approximation of  the parametrization of the deformed   rod at least if 
$||\hbox{dist}(\nabla_x v,SO(3))||_{L^2({\cal P}_\delta)}<<\delta$.
\medskip
\noindent {\ggras II.2.4. The boundary condition}
 \smallskip
Let us denote by $\GI_3$ the unit $3\times 3$ matrix and by $I_d$ the identity map of $\R^3$.

 In this subsection, we derive  boundary conditions  on the terms of the elementary deformation given by Theorem II.2.2. Indeed these conditions depend on the boundary conditions on the field $v$. We discuss essentially  the usual  case of a clamped condition on the  extremity of the rod defined by
 $$\Gamma_{0,\delta}=\Phi(\omega_\delta\times \{0\}).$$
 Then we assume that 
 $$v(x)=x \qquad \hbox{on} \quad \Gamma_{0,\delta}.$$
 In the following we show that the elementary deformation ${\cal V}(s_3)+\GR(s_3)\big(s_1\Gn_1+s_2\Gn_2\big)$ given by Theorem II.2.2 can be chosen such that $\overline{v}=0$ on the corresponding boundary $\omega_\delta\times \{0\}$.
Due to the definition (II.2.5) of ${\cal V}$, we first have 
$$ {\cal V}(0)=M(0),$$
 (recall that the point $M(0)$ is the beginning of the middle line of the rod ${\cal P}_\delta$; see the notations in  Subsection 3.1).
Now we prove that in Theorem II.2.2  the choice  $\GR(0)=\GI_3$ is licit as a boundary condition for the matrix $\GR(s_3)$.
Estimates (II.2.3) for the first portion ${\cal P}_{\delta,\alpha_0}$ imply 
$$||v(.,., 0)-\Ga_{\alpha_1}-\GR_{\alpha_0}\big(s_1\Gn_1(0) +s_2\Gn_2(0)\big)||^2_{(L^2(\omega_\delta))^3}\le C\delta dist\big(\nabla_x v; SO(3)\big)||^2_{L^2({\cal P}_{\delta,\alpha_0})}.$$
Using now the boundary  condition written in the equivalent form $v(s_1,s_2,0)=M(0)+s_1\Gn_1(0) +s_2\Gn_2(0)$ in the above estimate and using again $\ds \int_\omega s_1 ds_1ds_2=\int_\omega s_2 ds_1ds_2=\int_\omega s_1s_2 ds_1ds_2=0$ lead to 

$$|||\GR_{\alpha_0}-\GI_3|||^2\le {C\over \delta^3}||dist\big(\nabla_x v; SO(3)\big)||^2_{L^2(
{\cal P}_{\delta,\alpha_0})}.$$
As a consequence we can substitute  $\GR_{\alpha_0}$ by $\GI_3$ in the construction of the function $\GR(s_3)$  without altering  estimates (II.2.2)  so that $\GR(0)=\GI_3$. Indeed this leads to $\overline{v}=0$ on the boundary $\omega_\delta\times \{0\}$. 
The above arguments can be easily adapted if the imposed deformation  on $ \Gamma_{0,\delta}$ is of the form $v(x)=A+\GP(x-M(0))$ where $A$ is the deformation of the point $M(0)$ and $\GP$ is $3\times 3$ matrix. This leads to two boundary conditions of the type $ {\cal V}(0)=A$ and $\GR(0)=\GQ$ where, as an example,  the rotation $\GQ$ minimizes the distance from $\GP$ to $SO(3)$. In the same spirit, if the rod is  clamped on the two extremities of ${\cal P}_{\delta}$, one can modify the construction of the function $\GR(s_3)$ in such a way that $\GR(0)=\GR(L)=\GI_3$ keeping (II.2.2) true.

 \medskip
\noindent{\Ggras II.3.    $ H^1$- Estimates }
  \medskip
Throughout the paper we now assume that the boundary $\Gamma_{0,\delta}$ is clamped so that 
$${\cal V}(0)=M(0)={\cal V}_B(0) ,\qquad{\cal V}_S(0)=0\qquad  \hbox{and} \qquad \GR(0)=\GI_3 \leqno (II.3.1)$$
and we denote by $C$ a generic constant independent of $\delta$.  
\par
We derive a first $ H^1$- estimates using directly (II.2.2) and the fact that $\|\GR\Gt\|_2=1$. It gives
$$\Bigl\|{d{\cal V}\over ds_3}\Big\|_{(L^2(0,L))^3}\le C\Big(1+{1\over \delta}||\hbox{dist}(\nabla_x v,SO(3))||_{L^2({\cal P}_\delta)}\Big).\leqno (II.3.2)$$
Using the boundary condition (II.3.1), it leads to 
$$ \bigl\|{\cal V}\big\|_{(L^2(0,L))^3}\le  C\Big(1+{1\over \delta}||\hbox{dist}(\nabla_x v,SO(3))||_{L^2({\cal P}_\delta)}\Big).\leqno (II.3.3)$$
Inserting (II.3.2), (II.3.3) into (II.2.1) and using the estimates of (II.2.2), we deduce that 
$$\left \{\eqalign{
&||v||_{(L^2({\cal P}_\delta))^3}+||\nabla_x v||_{(L^2({\cal P}_\delta))^{3\times 3}}\le  C\Big(\delta+||\hbox{dist}(\nabla_x v,SO(3))||_{L^2({\cal P}_\delta)}\Big),\cr
&||v-{\cal V}||_{(L^2({\cal P}_\delta))^3}\le  C\delta\Big(\delta+||\hbox{dist}(\nabla_x v,SO(3))||_{L^2({\cal P}_\delta)}\Big).}\right.
\leqno (II.3.4)$$
Notice that using (II.2.2) also leads to the following estimates
$$\left\{\eqalign{
&\bigl\|{\GR}-\GI_3\big\|_{(L^2(0,L))^{3\times 3}}\le {C\over \delta^2} ||\hbox{dist}(\nabla_x v,SO(3))||_{L^2({\cal P}_\delta)}\cr
& \Bigl\|{d{\cal V}\over ds_3}- \Gt\Big\|_{(L^2(0,L))^3}\le {C\over \delta^2}||\hbox{dist}(\nabla_x v,SO(3))||_{L^2({\cal P}_\delta)}.\cr}\right.\leqno(II.3.5)$$
Since $\ds \Gt={ dM\over ds_3}$ the last inequality of (II.3.5) together with the boundary condition (II.3.1) give 
$$ \bigl\|{\cal V}- M\big\|_{(L^2(0,L))^3}\le {C\over \delta^2}||\hbox{dist}(\nabla_x v,SO(3))||_{L^2({\cal P}_\delta)}.\leqno (II.3.6)$$
From  (II.2.1), (II.3.5), (II.3.6)  and  estimates (II.2.2), we deduce that 
$$||v-I_d||_{(L^2({\cal P}_\delta))^3}+||\nabla_x v-\GI_3||_{(L^2({\cal P}_\delta))^{3\times 3}}\le  {C\over \delta}||\hbox{dist}(\nabla_x v,SO(3))||_{L^2({\cal P}_\delta)}.\leqno (II.3.7)$$ From (II.2.2) and (II.3.5) we also have $$\bigl\|{\GR}-\GI_3\big\|_{(H^1(0,L))^{3\times 3}}\le {C\over \delta^2} ||\hbox{dist}(\nabla_x v,SO(3))||_{L^2({\cal P}_\delta)}.\leqno(II.3.8)$$
From (II.2.2) and   (II.3.8) and the fact that $(v-I_d)-({\cal V}-M)=(\GR-\GI_3)(s_1\Gn_1+s_2\Gn_2)+\overline{v}$, we obtain
$$||(v-I_d)-({\cal V}-M)||_{(L^2({\cal P}_\delta))^3}\le  C ||\hbox{dist}(\nabla_x v,SO(3))||_{L^2({\cal P}_\delta)}\leqno (II.3.9)$$
 and
$$\bigl\|\nabla_x v+(\nabla_x v)^T-2\GI_3\big\|_{(L^2(\Omega_\delta))^{3\times 3}}\le C||\hbox{dist}(\nabla_x v,SO(3))||_{L^2({\cal P}_\delta)}+C\delta\bigl\|{\GR}+\GR^T-2\GI_3\big\|_{(L^2(0,L))^{3\times 3}}.$$ Due to (II.3.8), the continuous embedding of $H^1(0,L)$ in $C^0([0,L])$ and the equality ${\GR}+\GR^T-2\GI_3=\GR^T\bigl(\GR-\GI_3)^2$ we finally obtain
$$\left\{\eqalign{
\bigl\|\nabla_x v+(\nabla_x v)^T-2\GI_3\big\|_{(L^2({\cal P}_\delta))^{3\times 3}}\le& C||\hbox{dist}(\nabla_x v,SO(3))||_{L^2({\cal P}_\delta)}\cr
+&{C\over \delta^3}||\hbox{dist}(\nabla_x v,SO(3))||^2_{L^2({\cal P}_\delta)}.\cr}\right.\leqno(II.3.10)$$

\medskip
\noindent {\GGgras III. Asymptotic behavior of a sequence of deformations}\medskip
In view of the first estimate in (II.3.5) we can distinguish three main cases for the behavior of the quantity  $||\hbox{dist}(\nabla_x v,SO(3))||_{L^2({\cal P}_\delta)}$ (which will be a bound from below of the elastic energy)
$$||\hbox{dist}(\nabla_x v,SO(3))||_{L^2({\cal P}_\delta)}=\left\{\eqalign{
&O(\delta^{\kappa}),\qquad 1\le \kappa<2,\cr
&O(\delta^2),\cr
&O(\delta^{\kappa}),\qquad \kappa>2.\cr}\right.$$
This hierarchy of behavior for $||\hbox{dist}(\nabla_x v,SO(3))||_{L^2({\cal P}_\delta)}$ has already be observed in terms of elastic energy in [15].  

In this section we investigate the behavior of a sequence of deformations of ${\cal P}_\delta$ for  $||\hbox{dist}(\nabla_x v,SO(3))||_{L^2({\cal P}_\delta)}=O(\delta^{\kappa}),$ for $\kappa\ge 2.$ Actually, the first estimate in (II.3.5)  is useless if  $\kappa<2$ then we can not analyze this case using the decomposition (II.2.1). As usual we first rescale $\Omega_\delta$ in order to work over a fix domain in Subsection III.1. In Subsection III.2 we investigate the case $\kappa=2$ while in Subsection III.3 we deal with $\kappa>2$. Let us emphasize that we explicit the limit of the Green-St Venant's tensor in both cases. Subsection III.4 gives a few comparisons with the linearized deformations.
\medskip
\noindent {\Ggras III.1.    Rescaling of $\Omega_\delta$ }
  \medskip
We set $\Omega=\omega\times (0,L)$ and, we rescale $\Omega_\delta$ using the operator 
$$(\Pi_\delta \phi)(S_1,S_2,s_3)=\phi(\delta S_1,\delta S_2,s_3)\hbox{ for any}\quad (S_1,S_2,s_3) \in \Omega$$
defined for any function $\phi$ defined over $\Omega_\delta$. Indeed, if $\phi\in L^2(\Omega_\delta)$ then $ (\Pi_\delta \phi)\in  L^2(\Omega)$.
The estimates (II.2.2) of $\overline{v}$   transposed over $\Omega$ are 
$$\left\{\eqalign{
&||\Pi_\delta \overline{v}||_{(L^2(\Omega))^3}\le C ||\hbox{dist}(\nabla_x v,SO(3))||_{L^2({\cal P}_\delta)}\cr
&||{\partial \Pi _\delta\overline {v}\over \partial S_1}||_{(L^2(\Omega))^3}\le C ||\hbox{dist}(\nabla_x v,SO(3))||_{L^2({\cal P}_\delta)}\cr 
&||{\partial \Pi _\delta\overline{v}\over \partial S_2}||_{(L^2(\Omega))^3}\le C ||\hbox{dist}(\nabla_x v,SO(3))||_{L^2({\cal P}_\delta)} \cr 
&||{\partial \Pi_\delta \overline{v}\over \partial s_3}||_{(L^2(\Omega))^3}\le {C \over \delta} ||\hbox{dist}(\nabla_x v,SO(3))||_{L^2({\cal P}_\delta)} \cr}\right.\leqno(III.1.1)$$
while estimates (II.3.4) and (II.3.7) on  $v$ become
$$\left\{\eqalign{
&||\Pi_\delta v||_{(L^2(\Omega))^3}\le  C\Big(1+{1\over \delta}||\hbox{dist}(\nabla_x v,SO(3))||_{L^2({\cal P}_\delta)}\Big)\cr
&||{\partial \Pi_\delta v \over \partial S_1}||_{(L^2(\Omega))^3}\le C\Big(\delta+||\hbox{dist}(\nabla_x v,SO(3))||_{L^2({\cal P}_\delta)}\Big)\cr
&||{\partial \Pi_\delta v \over \partial S_2}||_{(L^2(\Omega))^3}\le  C\Big(\delta+||\hbox{dist}(\nabla_x v,SO(3))||_{L^2({\cal P}_\delta)}\Big)\cr
&||{\partial \Pi_\delta v \over \partial s_3}||_{(L^2(\Omega))^3}\le   C\Big(1+{1\over \delta}||\hbox{dist}(\nabla_x v,SO(3))||_{L^2({\cal P}_\delta)}\Big)\cr}\right. \leqno (III.1.2)$$ and 
$$\left\{\eqalign{
&||\Pi_\delta (v-I_d)||_{(L^2(\Omega))^3}\le  {C\over \delta^2}||\hbox{dist}(\nabla_x v,SO(3))||_{L^2({\cal P}_\delta)}\cr
&||{\partial \Pi_\delta (v-I_d)\over \partial S_1}||_{(L^2(\Omega))^3}\le  {C\over \delta}||\hbox{dist}(\nabla_x v,SO(3))||_{L^2({\cal P}_\delta)}\cr
&||{\partial \Pi_\delta (v-I_d)\over \partial S_2}||_{(L^2(\Omega))^3}\le  {C\over \delta}||\hbox{dist}(\nabla_x v,SO(3))||_{L^2({\cal P}_\delta)}\cr
&||{\partial \Pi_\delta (v-I_d)\over \partial s_3}||_{(L^2(\Omega))^3}\le  {C\over \delta^2}||\hbox{dist}(\nabla_x v,SO(3))||_{L^2({\cal P}_\delta).}\cr}\right. \leqno (III.1.3)$$
All the estimates given in Section II.2.3 over $\Omega_\delta$ can be  easily transposed over $\Omega$.
\medskip
\noindent{\Ggras III.2.  Limit behavior  for  a sequence such that $||\hbox{dist}(\nabla_x v_\delta,SO(3))||_{L^2({\cal P}_\delta)} \sim \delta^2$ }
  \medskip
Let us consider a sequence of deformations $v_\delta$ of $\big(H^1({\cal P}_\delta)\big)^3 $  such that $v_\delta =I_d$ on $\Gamma_{0,\delta}$ and
$$||\hbox{dist}(\nabla_x v_\delta,SO(3))||_{L^2({\cal P}_\delta)} \le C \delta^2.\leqno(III.2.1)$$
Indeed different boundary conditions on $v_\delta$ can be considered on both the extremities of  ${\cal P}_\delta$ (see Subsection II.2.4).
We denote by ${\cal V}_\delta$, $\GR_\delta$ and $\overline{v}_\delta$ the three terms of the decomposition of $v_\delta$ given by Theorem II.2.2 and by ${\cal V}_{B,\delta}$ and ${\cal V}_{S,\delta}$  the two terms given by (II.2.16).  The estimates (II.2.2), (II.2.17), (II.3.5), (II.3.6), (III.1.1), (III.1.3) lead to the following lemma:

\noindent {\ggras  Lemma III.2.1. }
{\it There exists a subsequence still indexed by $\delta$ such that 
$$\left\{\eqalign{
&\GR_\delta \rightharpoonup \GR \quad \hbox{weakly in}\quad \big(H^1(0,L)\big)^{3\times 3}\cr
&{\cal V}_\delta\longrightarrow   {\cal V} \quad \hbox{strongly in}\quad \big(H^1(0,L)\big)^3\cr
&{\cal V}_{B,\delta}\longrightarrow   {\cal V} \quad \hbox{strongly in}\quad \big(H^1(0,L)\big)^3\cr
&{1\over \delta}{\cal V}_{S,\delta}\rightharpoonup {\cal V}_S\quad \hbox{weakly in}\quad \big(H^1(0,L)\big)^3\cr
&{1\over \delta^2}\Pi _\delta\overline {v}_\delta \rightharpoonup  \overline {v} \quad \hbox{weakly in}\quad \big(L^2(0,L;H^1(\omega))\big)^3\cr}\right.\leqno(III.2.2)$$ Moreover $\GR(s_3)$ belongs to $SO(3)$ for any $s_3 \in [0,L]$, ${\cal V}\in \big(H^2(0,L)\big)^3 $ and they satisfy 
$${\cal V}(0)=M(0), \qquad  \GR(0)=\GI_3,\qquad {\cal V}_S(0)=0,\qquad\hbox{and} \qquad{d{\cal V}\over ds_3}=\GR \Gt .\leqno (III.2.3)$$ 
\noindent Furthermore, we also have
$$\left\{\eqalign{
\Pi _\delta v_\delta&\longrightarrow    {\cal V} \quad \hbox{strongly in}\quad \big(H^1(\Omega)\big)^3,\cr
\Pi _\delta (\nabla_x v_\delta)& \longrightarrow  \GR\quad \hbox{strongly in}\quad \big(L^2(\Omega)\big)^{3\times 3}.\cr}\right.\leqno(III.2.4)$$}

The fields ${\cal V}$ and $\GR$ which are defined in Lemma III.2.1 describe the deformation of the limit 1D curved  rod as a deformation of the middle line ${\cal V}$ and a rotation of each cross section $\GR$. The convergences (III.2.4) show that $({\cal V},\GR)$ is the limit of the deformation of the (rescaled) 3D curved rod and represents the elementary deformation of this rod. \par
The last relation in (III.2.3) shows that $\ds \big \|{d{\cal V}\over ds_3}\big \|_2 =1$ and then the variable $s_3$ remains the arc length of the middle line of the deformed configuration. As a consequence, in the present case, there is no extension-compression of order 1 the rod. Then the limit behavior is essentially a bending model. At least, the forth convergence (III.2.2) means that the quantity ${\cal V}_S$  which describes the stretching  is of order $\delta$.\par
The following corollary gives a corrector result. 
\smallskip
\noindent{\ggras  Corollary III.2.2. } {\it For the same subsequence as in Lemma III.2.1, we  have
$$\left\{\eqalign{
{1\over \delta}\big(\Pi _\delta (\nabla_x v_\delta)-\GR_\delta\big)\Gn_\alpha& \rightharpoonup  {\partial\overline{v}\over \partial S_\alpha}\quad \hbox{weakly in}\quad (L^2(\Omega))^3,\cr
{1\over \delta}\big(\Pi _\delta (\nabla_x v_\delta)-\GR_\delta\big)\Gt &\rightharpoonup  {d\GR\over ds_3}\big(S_1\Gn_1+S_2\Gn_2)+{d{\cal V}_S\over ds_3} \quad \hbox{weakly in}\quad (L^2(\Omega))^3.\cr}\right.\leqno(III.2.5)$$}
\noindent{\ggras  Proof. } The first convergence in (III.2.5) is a direct consequence of (II.2.13) and (III.2.2). In order to obtain the second convergence, remark first that, thanks to  estimates (III.1.1)  and (III.2.1) the sequence $\ds{1\over \delta}\Pi_\delta\overline{v}_\delta$  is bounded in $H^1(\Omega)$. Due to (III.2.2), its weak limit must be equal to $0$. Using now (II.2.15) and the convergences  (III.2.2) leads to the result.\fin
\noindent 
To end this section, let us notice that the strong convergences in (III.2.4) together with the relation 
$(\nabla_xv_\delta)^T\nabla_x v_\delta-\GI_3=(\nabla_xv_\delta-\GR_\delta)^T\nabla_x v_\delta+(\GR_\delta)^T(\nabla_x v_\delta-\GR_\delta)$ permit to obtain the limit of the Green-St Venant's tensor in the rescaled domain $\Omega$
$${1\over 2\delta}\Pi _\delta \big((\nabla_xv_\delta)^T\nabla_x v_\delta-\GI_3\big)\rightharpoonup   \GE \qquad\hbox{weakly in}\quad (L^1(\Omega))^{3\times 3},\leqno(III.2.6)$$ where 
$$\eqalign{
\GE=&{1\over 2}\Big\{(\Gn_1\,|\, \Gn_2\,|\, \Gt) \Bigl({\partial\overline{v}\over \partial S_1}\,|\, {\partial\overline{v}\over \partial S_2}\,|\,{d\GR\over ds_3}\big(S_1\Gn_1+S_2\Gn_2)+{d{\cal V}_S\over ds_3}\Big)^T\GR\cr
+& \GR^T \Bigl({\partial\overline{v}\over \partial S_1}\,|\, {\partial\overline{v}\over \partial S_2}\,|\,{d\GR\over ds_3}\big(S_1\Gn_1+S_2\Gn_2)+{d{\cal V}_S\over ds_3}\Big) (\Gn_1\,|\, \Gn_2\,|\, \Gt)^T\Big\}\cr}\leqno(III.2.7)$$ and where  $\big(\Ga\, |\, \Gb\, |\, \Gc\big)$ denotes the $3\times 3$ matrix with  columns   $\Ga$,  $\Gb$ and  $\Gc$.
\medskip
Setting $\overline{w}=\GR^T\overline{v}$ and using the fact that the matrix $\ds \GR^T{d\GR\over ds_3}$ is antisymmetric, we can write $\GE$ as
$$\GE=(\Gn_1\,|\, \Gn_2\,|\, \Gt)\, \widehat{\GE}\,(\Gn_1\,|\, \Gn_2\,|\, \Gt)^T,\leqno (III.2.8)$$ where the symmetric matrix $\widehat{\GE}$ is defined by  
$$\widehat{\GE}= \pmatrix{
\ds{\partial\overline{w}\over \partial S_1}\cdot\Gn_1 & \ds{1\over 2}\Big\{{\partial\overline{w}\over \partial S_1}\cdot\Gn_2+{\partial\overline{w}\over \partial S_2}\cdot\Gn_1\Big\} & \ds {1\over 2}\Big\{{\partial\overline{w}\over \partial S_1}\cdot\Gt-S_2{d\GR\over ds_3}\Gn_1\cdot \GR\Gn_2+{d{\cal V}_S\over ds_3}\cdot\GR\Gn_1 \Big\}\cr
* & \ds{\partial\overline{w}\over \partial S_2}\cdot\Gn_2   &  \ds {1\over 2}\Big\{{\partial\overline{w}\over \partial S_2}\cdot\Gt+S_1{d\GR\over ds_3}\Gn_1\cdot \GR\Gn_2+{d{\cal V}_S\over ds_3}\cdot\GR\Gn_2\Big\}\cr
 * & *  & \ds  -S_1{d\GR\over ds_3}\Gt\cdot \GR\Gn_1-S_2{d\GR\over ds_3}\Gt\cdot \GR\Gn_2+{d{\cal V}_S\over ds_3}\cdot\GR\Gt}.\leqno(III.2.9)$$ 
 \smallskip
\noindent{\ggras  Corollary III.2.3. } {\it Assume that
$$\forall\delta\in]0,\delta_0],\qquad \det\big(\nabla v_\delta(x)\big)>0\qquad \hbox{for a.e. $x\in {\cal P}_\delta$}$$ then, for the same subsequence as in Lemma III.2.1, we  have 
$$||\GE||_{(L^2(\Omega))^{3\times 3}}=||\widehat{\GE}||_{(L^2(\Omega))^{3\times 3}}\le
\liminf_{\delta\to 0}{1\over \delta^2} ||\hbox{dist}(\nabla_x v_\delta,SO(3))||_{L^2({\cal P}_\delta)} . \leqno(III.2.10)$$}
\noindent{\ggras  Proof. } The map $A\to \sqrt{A^TA}$ is continuous from the space of the $3\times 3$ matrices  into the set of all symmetric  matrices and we have $|||\sqrt{A^TA}|||=|||A|||$,  where $|||\cdot|||$ is the Frobenius norm. Then, the second  strong convergence in (III.2.4) gives $$\Pi _\delta \Big(\sqrt{(\nabla_xv_\delta)^T\nabla_x v_\delta}\Big)\longrightarrow  \GI_3\quad \hbox{strongly in}\quad \big(L^2(\Omega)\big)^{3\times 3}.$$  Estimate (III.2.1) implies that the sequence $\ds {1\over \delta}\Pi _\delta \Big(\sqrt{(\nabla_xv_\delta)^T\nabla_x v_\delta}-\GI_3\Big)$ is bounded in $(L^2(\Omega))^{3\times 3}$.  The  identity $(\nabla_xv_\delta)^T\nabla_x v_\delta-\GI_3=\big(\sqrt{(\nabla_xv_\delta)^T\nabla_x v_\delta}-\GI_3\big)\big(\sqrt{(\nabla_xv_\delta)^T\nabla_x v_\delta}+\GI_3\big)$, the weak convergence (III.2.6) and  the above strong convergence give 
$${1\over \delta}\Pi _\delta \Big(\sqrt{(\nabla_xv_\delta)^T\nabla_x v_\delta}-\GI_3\Big)\rightharpoonup     \GE \qquad\hbox{weakly in}\quad (L^2(\Omega))^{3\times 3}.$$ We recall that for any $3\times 3$ matrix $\GA$ such that $\det(\GA)>0$, we have $\hbox{dist}(\GA,SO(3))=|||\sqrt{A^TA}-\GI_3|||$. 
By weak lower semi-continuity of the norm, we obtain  the result.
 \fin
 \medskip
\noindent{\Ggras III.3.  Limit behavior  for  a sequence such that $||\hbox{dist}(\nabla_x v_\delta,SO(3))||_{L^2({\cal P}_\delta)} \sim \delta^\kappa$ for $\kappa>2$. }
  \medskip
Let us consider a sequence of deformations $v_\delta$ of $\big(H^1({\cal P}_\delta)\big)^3 $  such that $v_\delta =I_d$ on $\Gamma_{0,\delta}$ and
$$||\hbox{dist}(\nabla_x v_\delta,SO(3))||_{L^2({\cal P}_\delta)} \le C \delta^\kappa.$$
The estimates   (II.3.7) and (III.1.3) lead to the following convergences:
$$\left\{\eqalign{
\GR_\delta&\longrightarrow    \GI_3 \quad \hbox{strongly in}\quad \big(H^1(0,L)\big)^{3\times 3},\cr
\Pi _\delta v_\delta&\longrightarrow    M \quad \hbox{strongly in}\quad \big(H^1(\Omega)\big)^3,\cr
\Pi _\delta (\nabla_x v_\delta)& \longrightarrow  \GI_3\quad \hbox{strongly in}\quad \big(L^2(\Omega)\big)^{3\times 3}.\cr}\right.\leqno(III.3.1)$$We now study the asymptotic behavior of the sequence of displacements
$$u_\delta=v_\delta-I_d.$$  Due to decomposition (II.2.1) we write
$$u_\delta(s)={\cal U}_\delta(s_3)+(\GR_\delta -\GI_3)(s_3)\big(s_1\Gn_1(s_3)+s_2\Gn_2(s_3)\big)
+\overline{v}_\delta(s),\qquad s\in \Omega_\delta,\leqno (III.3.2)$$ where ${\cal U}_\delta(s_3)={\cal V}_\delta(s_3)-M(s_3)=\big({\cal V}_{B,\delta}(s_3)-M(s_3)\big)+{\cal V}_{S,\delta}(s_3)$ and we have the following Lemma:
\medskip
\noindent {\ggras  Lemma III.3.1. }
{\it There exists a subsequence still indexed by $\delta$ such that 
$$\left\{\eqalign{
{1\over \delta^{\kappa-2}}\bigl(\GR_\delta-\GI_3\big)&\rightharpoonup \GA \quad \hbox{weakly in}\quad \big(H^1(0,L)\big)^{3\times 3}\cr
{1\over \delta^{\kappa-2}} {\cal U}_\delta &\longrightarrow   {\cal U} \quad \hbox{strongly in}\quad \big(H^1(0,L)\big)^3\cr
{1\over \delta^{\kappa-1}}{\cal V}_{S,\delta}& \rightharpoonup {\cal V}_S\quad \hbox{weakly in}\quad \big(H^1(0,L)\big)^3\cr
{1\over \delta^\kappa}\Pi _\delta\overline {v}_\delta& \rightharpoonup  \overline {v} \quad \hbox{weakly in}\quad \big(L^2(0,L;H^1(\omega))\big)^3\cr}\right.\leqno(III.3.3)$$ The function ${\cal U}$ belongs to $\big(H^2(0,L)\big)^3$, for any $s_3\in [0,L]$ the matrix $\GA(s_3)$ is antisymmetric and the following relations hold true:
$${\cal U}(0)={\cal V}_S(0)=0, \qquad  \GA(0)=0\qquad\hbox{and} \qquad{d{\cal U}\over ds_3}=\GA \Gt ,\leqno (III.3.4)$$ Moreover we have $$\left\{\eqalign{
{1\over \delta^{\kappa-2}}\Pi _\delta u_\delta&\longrightarrow    {\cal U} \quad \hbox{strongly in}\quad \big(H^1(\Omega)\big)^3,\cr
{1\over \delta^{\kappa-2}}\Pi _\delta (\nabla_x u_\delta)& \longrightarrow  \GA\quad \hbox{strongly in}\quad \big(L^2(\Omega)\big)^{3\times 3}.\cr}\right.\leqno(III.3.5)$$ }
\smallskip
\noindent  {\ggras Proof.}   The convergences (III.3.3), (III.3.5) and the relations (III.3.4) follow  directly from estimates (II.2.2), (II.3.6), (II.3.7) and (III.1.3). It remains to prove that $\GA(s_3)$ is antisymmetric.
Using  the first convergence in (III.3.1) and  the first convergence in   (III.3.3) we get
$${1\over \delta^{\kappa-2}}\GR^T_\delta\bigl(\GR_\delta-\GI_3\big)\longrightarrow \GA \quad \hbox{strongly in}\quad \big(L^2(0,L)\big)^{3\times 3}.$$ The matrix $\GR_\delta$ belongs to $SO(3)$, hence $\GR^T_\delta\bigl(\GR_\delta-\GI_3\big)=-(\GR_\delta-\GI_3)^T$. Then, from the first convergence in (III.3.3), we deduce that  the   matrix $\GA(s_3)$ is  antisymmetric. \fin
Since $\GA$ is antisymmetric, there exists a  field ${\cal R}\in \big(H^1(0,L)\big)^3$ (with ${\cal R}(0)=0$) such that for all $x\in \R^3$ 
$$\GA x={\cal R}\land x.\leqno(III.3.6)$$ From (III.3.4) and the above equality we obtain 
$${d{\cal U}\over ds_3}={\cal R}\land\Gt.\leqno(III.3.7)$$
The relation (III.3.6) means that, at the order $\delta^{\kappa-2}$, the cross sections of ${\cal P}_\delta$ are submitted to small rotations and (III.3.7) shows that the limit displacement is of Bernoulli-Navier's type. 
\smallskip
\noindent{\ggras  Corollary III.3.2. } {\it For the same subsequence as in Lemma III.3.1, we have
$$\left\{\eqalign{
{1\over \delta^{\kappa-1}}\big(\Pi _\delta (\nabla_x v_\delta)-\GR_\delta\big)\Gn_\alpha& \rightharpoonup  {\partial\overline{v}\over \partial S_\alpha}\quad \hbox{weakly in}\quad (L^2(\Omega))^3,\cr
{1\over \delta^{\kappa-1}}\big(\Pi _\delta (\nabla_x v_\delta)-\GR_\delta\big)\Gt &\rightharpoonup  {d{\cal R}\over ds_3}\land\big(S_1\Gn_1+S_2\Gn_2)+{d{\cal V}_S\over ds_3}\quad \hbox{weakly in}\quad (L^2(\Omega))^3.\cr}\right.\leqno(III.3.8)$$}
\noindent{\ggras  Proof. } The proof of Corollary III.3.2 is similar to that of Corollary III.2.2, but using  now (II.2.13), (II.2.15) and  the convergences of Lemma III.3.1.\fin
\noindent From Lemma III.3.1 and Corollary  III.3.2 we deduce the limit of the Green-St Venant's tensor in the rescaled domain $\Omega$
$${1\over 2\delta^{\kappa-1}}\Pi _\delta \big((\nabla_xv_\delta)^T\nabla_x v_\delta-\GI_3\big)\rightharpoonup  \GE \qquad\hbox{weakly in}\quad (L^1(\Omega))^{3\times 3},\leqno(III.3.9)$$   where the symmetric matrix $\GE$ is defined by 
$$\eqalign{
\GE=&{1\over 2}\Big\{\Bigl({\partial\overline{v}\over \partial S_1}\,|\, {\partial\overline{v}\over \partial S_2}\,|\,{d{\cal R}\over ds_3}\land\big(S_1\Gn_1+S_2\Gn_2)+{d{\cal V}_S\over ds_3}\Big)^T (\Gn_1\,|\, \Gn_2\,|\, \Gt) \cr
&+ (\Gn_1\,|\, \Gn_2\,|\, \Gt) ^T\Bigl({\partial\overline{v}\over \partial S_1}\,|\, {\partial\overline{v}\over \partial S_2}\,|\,{d{\cal R}\over ds_3}\land\big(S_1\Gn_1+S_2\Gn_2)+{d{\cal V}_S\over ds_3}\Big)\Big\}\cr}$$ We can write $\GE$ in the form
$$\GE=(\Gn_1\,|\, \Gn_2\,|\, \Gt)\widehat{\GE} (\Gn_1\,|\, \Gn_2\,|\, \Gt)^T $$ where the symmetric matrix $\widehat\GE$ is defined by 
$$\widehat{\GE}=\pmatrix{
\ds{\partial\overline{v}\over \partial S_1}\cdot\Gn_1 & \ds{1\over 2}\Big\{{\partial\overline{v}\over \partial S_1}\cdot\Gn_2+{\partial\overline{v}\over \partial S_2}\cdot\Gn_1\Big\} & \ds {1\over 2}\Big\{{\partial\overline{v}\over \partial S_1}\cdot\Gt-S_2{d{\cal R}\over ds_3}\cdot \Gt+{d{\cal V}_S\over ds_3}\cdot \Gn_1 \Big\}\cr
* & \ds{\partial\overline{v}\over \partial S_2}\cdot\Gn_2   &  \ds {1\over 2}\Big\{{\partial\overline{v}\over \partial S_2}\cdot\Gt+S_1{d{\cal R}\over ds_3} \cdot \Gt+{d{\cal V}_S\over ds_3}\cdot \Gn_2\Big\}\cr
 * & *  & \ds  -S_1{d{\cal R}\over ds_3} \cdot\Gn_2+S_2{d{\cal R}\over ds_3} \cdot\Gn_1+{d{\cal V}_S\over ds_3}\cdot \Gt}.\leqno(III.3.10)$$
From Lemma III.3.1 and the above convergences, we deduce the analog of Corollary II.2.3.

 \noindent{\ggras  Corollary III.3.3. } {\it Assume that
$$\forall\delta\in]0,\delta_0],\qquad \det\big(\nabla v_\delta(x)\big)>0\qquad \hbox{for a.e. $x\in {\cal P}_\delta$}$$ then, for the same subsequence as in Lemma III.3.1, we  have 
$$||\GE||_{(L^2(\Omega))^{3\times 3}}=||\widehat{\GE}||_{(L^2(\Omega))^{3\times 3}}\le
\liminf_{\delta\to 0}{1\over \delta^{\kappa}} ||\hbox{dist}(\nabla_x v_\delta,SO(3))||_{L^2({\cal P}_\delta)} . \leqno(III.3.11)$$}

 \medskip
\noindent {\Ggras III.4.    Comparison with linearized deformations}
  \medskip
In this subsection, we always consider  a sequence of deformations $v_\delta$ of $\big(H^1({\cal P}_\delta)\big)^3 $  satisfying  $v_\delta =I_d$ on $\Gamma_{0,\delta}$ and
$$||\hbox{dist}(\nabla_x v_\delta,SO(3))||_{L^2({\cal P}_\delta)} \le C \delta^\kappa.$$
We recall the decomposition (III.3.2) of the displacement $u_\delta=v_\delta-I_d$.
 
Let us notice that  (II.3.5) shows that for $\kappa>2$, both the displacement and its gradient are small (with respect to $\delta$). One can then address the problem of comparing the limit displacement ${\cal U}_\delta$ and the limit displacement in  the framework of small deformations. To this end let us first recall the decomposition of displacement for small deformations.
 
We define the strain semi-norm $|\cdot|_{\cal E}$ by setting 
$$\forall w\in \big(H^1({\cal P}_\delta)\big)^3 \qquad |w|_{\cal E}=\big\|{1\over 2}\big(\nabla_x w+(\nabla_x w)^T\big)\big\|_{(L^2({\cal P}_\delta))^{3\times 3}}$$  Now, using the results obtained in [13], we decompose $u_\delta$ in the sum of an elementary displacement and a warping
$$u_\delta(s)=U_{e,\delta} (s)+\overline{u}_\delta(s)={\cal U}_\delta(s_3)+{\cal R}_\delta(s_3)\land\big(s_1\Gn_1(s_3)+s_2\Gn_2(s_3)\big)+\overline{u}_\delta(s)\qquad \hbox{for a.e.}\enskip s\in \Omega_\delta.$$
The warping $\overline{u}_\delta$ satisfies the following equalities
$$\int_{\omega_\delta}\overline{u}_\delta(s_1,s_2,s_3)ds_1ds_2=0\qquad \int_{\omega_\delta}\overline{u}_\delta(s_1,s_2,s_3)\land\big(s_1\Gn_1(s_3)+s_2\Gn_2(s_3)\big)ds_1ds_2=0\qquad \hbox{for a.e.}\enskip s_3\in (0,L).$$  Notice that the first term ${\cal U}_\delta$ of the elementary displacement $U_{e,\delta} $ is the mean value of $u_\delta$ over the cross-section $\Phi(\omega_\delta\times\{s_3\})$ and then is the same as in (III.3.2).  Theorem 2.1 in [13] gives
$$\left\{\eqalign{
&\|\nabla_s \overline{u}_\delta\|_{(L^2(\Omega_\delta))^{3\times 3}}\le C|u_\delta|_{\cal E}\qquad \|\overline{u}_\delta\|_{(L^2(\Omega_\delta))^3}\le C\delta|u_\delta|_{\cal E}\cr
&\Big\|{d{\cal R}_\delta\over ds_3}\Big\|_{(L^2(0,L))^3}\le C{|u_\delta|_{\cal E}\over\delta^2}\qquad \Big\|{d{\cal U}_\delta\over ds_3}-{\cal R}_\delta\land\Gt\Big\|_{(L^2(0,L))^3}\le C{|u_\delta|_{\cal E}\over \delta}.\cr}\right.\leqno(III.4.1)$$  
 We recall (see [14]) the definitions of the inextensional displacements and   extensional displacements sets of the middle-line of the curved rod.
$$\left\{\eqalign{
&D_{In}=\Bigl\{U\in \big(H^1(0,L)\big)^3 \;|\; {dU\over ds_3}\cdot\Gt=0,\enskip U(0)=0\Big\}\cr
&D_{Ex}=\Bigl\{U\in \big(H^1(0,L)\big)^3 \;|\; {dU\over ds_3}\cdot\Gn_1={dU\over ds_3}\cdot\Gn_2=0,\enskip U(0)=0\Big\}\cr}\right.\leqno(III.4.2)$$ An element of $D_{In}$ is an inextensional displacement while an element of $D_{Ex}$ is an  extensional one. We recall (see [14]) that ${\cal U}_\delta$ can be written as the sum of an inextensional displacement and an extensional one
$${\cal U}_\delta=U_{I,\delta}+U_{E,\delta}\qquad U_{I,\delta}\in D_{In},\quad U_{E,\delta}\in D_{Ex},\leqno(III.4.3)$$ and we have (see again [14])
$$\Big\|{d{\cal U}_{\delta}\over ds_3}\Big\|_{(L^2(0,L))^3}+\Big\|{dU_{I,\delta}\over ds_3}\Big\|_{(L^2(0,L))^3}\le C{|u_\delta|_{\cal E}\over \delta^2}\qquad \Big\|{dU_{E,\delta}\over ds_3}\Big\|_{(L^2(0,L))^3}\le C{|u_\delta|_{\cal E}\over \delta}.\leqno(III.4.4)$$ 
In order to be obtain the same estimate on ${\cal U}_{\delta}$ that the one given by Lemma III.3.1.  in the previous section, we are lead to assume that $|u_\delta|_{\cal E}\le C\delta^\kappa$. Comparing with estimate (II.3.10) which gives
$$|u_\delta|_{\cal E}\le C(\delta^\kappa+\delta^{2\kappa-3})$$ since $||\hbox{dist}(\nabla_x v_\delta,SO(3))||_{L^2({\cal P}_\delta)} \le C \delta^\kappa$ we are led to  assume in the following that $\kappa \ge 3$.

The estimates (III.3.11), (III.4.1), (III.4.2) and (III.4.4) lead to the following lemma:
\medskip
\noindent{\ggras  Lemma III.4.1.  }{\it  We assume that $\kappa\ge 3$. There exists a subsequence (still indexed by $\delta$) of the sequence given in Lemma III.3.1 such that
$$\left\{\eqalign{
{1\over \delta^{\kappa-2}}U_{I,\delta} &\longrightarrow  {\cal U} \quad \hbox{strongly in}\quad \big(H^1(0,L)\big)^3 \cr 
{1\over \delta^{\kappa-1}}U_{E,\delta}& \rightharpoonup   U_E \quad \hbox{weakly in}\quad \big(H^1(0,L)\big)^3 \cr
{1\over \delta^{\kappa-2}}{\cal R}_{\delta} &\rightharpoonup  {\cal R} \quad \hbox{weakly in}\quad \big(H^1(0,L)\big)^3 \cr 
{1\over \delta^{\kappa-1}}\Big({d{\cal U}_\delta\over ds_3}&-{\cal R}_\delta\land \Gt\Big)\cdot \Gn_\alpha \rightharpoonup  {\cal Z}_\alpha \quad \hbox{weakly in}\quad L^2(0,L)\cr
{1\over \delta^{\kappa}}\Pi_\delta \overline{u}_{\delta} &\rightharpoonup \overline{u} \quad \hbox{weakly in}\quad \big(L^2(0,L;H^1(\omega))\big)^3 \cr
{1\over 2 \delta^{\kappa-1}}\Pi_\delta\bigl(\nabla u_\delta &+(\nabla u_{\delta})^T\big) \rightharpoonup \GE^{'}\quad \hbox{weakly in}\quad \big(L^2(\Omega)\big)^{3\times 3} \cr }\right.\leqno(III.4.5)$$  with $\GE^{'}=(\Gn_1\,|\, \Gn_2\,|\, \Gt)\widehat{\GE}^{'} (\Gn_1\,|\, \Gn_2\,|\, \Gt)^T$ where the symmetric matrix $\widehat\GE^{'}$ is  given by
$$\widehat{\GE}^{'}=\pmatrix{
\ds{\partial\overline{u}\over \partial S_1}\cdot\Gn_1 & \ds{1\over 2}\Big\{{\partial\overline{u}\over \partial S_1}\cdot\Gn_2+{\partial\overline{u}\over \partial S_2}\cdot\Gn_1\Big\} & \ds {1\over 2}\Big\{{\partial\overline{u}\over \partial S_1}\cdot\Gt-S_2{d{\cal R}\over ds_3}\cdot \Gt+{\cal Z}_1 \Big\}\cr
* & \ds{\partial\overline{u}\over \partial S_2}\cdot\Gn_2   &  \ds {1\over 2}\Big\{{\partial\overline{u}\over \partial S_2}\cdot\Gt+S_1{d{\cal R}\over ds_3} \cdot \Gt+{\cal Z}_2\Big\}\cr
 * & *  & \ds  -S_1{d{\cal R}\over ds_3} \cdot\Gn_2+S_2{d{\cal R}\over ds_3} \cdot\Gn_1+{dU_E\over ds_3}\cdot \Gt}.$$ Moreover the symmetric matrices $\GE$ and $\widehat{\GE}$ given in (III.3.10) satisfy
$$\widehat{\GE}=\left\{\eqalign{
&\widehat{\GE}^{'}+{1\over 2}\big(||{\cal R}||_2^2\GI_3-{\cal R}.^T{\cal R}\big)\qquad \hbox{if}\quad \kappa=3,\cr 
&\widehat{\GE}^{'}\qquad \hbox{if}\quad \kappa>3.\cr}\right.\leqno(III.4.6)$$}
\medskip
\noindent{\ggras  Proof. }  The convergences (III.4.5) and the expression of $\widehat{\GE}^{'}$ are proved in [14] taking into account (III.4.4) and the fact that $|u_\delta|_{\cal E}\le C\delta^\kappa$ with $\kappa\ge3$. Let us notice that the limit ${\cal U}$ in the first convergence (III.4.5) is the same that in (III.3.3). In order to prove (III.4.6) we first write the Green-St Venant  deformation tensor as (using $u_\delta=v_\delta-I_d$)
$${1\over 2\delta^{\kappa-1}}\Pi _\delta \big((\nabla_xv_\delta)^T\nabla_x v_\delta-\GI_3\big)={1\over 2\delta^{\kappa-1}}\Pi _\delta \big((\nabla_xu_\delta)^T+\nabla_x u_\delta\big)+{\delta^{\kappa-3}\over 2}\Pi _\delta \big({1\over \delta^{\kappa-2}}(\nabla_xu_\delta)^T{1\over \delta^{\kappa-2}}\nabla_x u_\delta\big). $$
Using convergences (III.3.5) and  (III.3.9), we deduce that 
$$\GE=\left\{\eqalign{
&\GE^{'}+{1\over 2}\GA^T\GA\quad \hbox{if}\quad \kappa=3\cr
&\GE^{'}\qquad \hbox{if}\quad \kappa>3,\cr}\right.$$
where the matrix $\GA$ is defined in (III.3.3) and (III.3.5). Recalling relation (III.3.6) between $\GA$ and ${\cal R}$ leads to (III.4.6).
  \fin
In the case $\kappa>3$, Lemma III.4.1 shows first that starting from nonlinear deformations leads exactly to the same deformation model that starting from linearized deformations. The comparison in the case where $\kappa=3$ is more intricate. This is due to the definitions of the two warping $\overline {v}$ and  $\overline {u}$ which do not satisfy the same geometrical conditions (see (II.2.7) for $\overline {v}$  and the beginning of this section  for $\overline {u}$). The second difference concerns the comparison between the stretching deformation ${\cal V}_S$ and the extentionnal displacement $U_E$. While it is easy to see that 
$\ds{dU_E\over ds_3}\cdot\Gt={d{\cal V}_S\over ds_3}\cdot\Gt$ for $\kappa>3$, in the case where $\kappa=3$ one obtains $\ds{dU_E\over ds_3}\cdot\Gt={d{\cal V}_S\over ds_3}\cdot\Gt-{1\over 2}\Big\|{d{\cal U}\over ds_3}\Big\|^2_2$. The correcting term $\ds{1\over 2}\Big\|{d{\cal U}\over ds_3}\Big\|^2_2$ actually comes from the limit contribution of the term $\ds{1\over 2\delta^{2}}(\GR_\delta -\GI_3)\Gt  \cdot (\GR_\delta -\GI_3)\Gt$.
\bigskip
\noindent {\GGgras IV. Asymptotic behavior of an elastic curved rod}
\medskip
This section is devoted to use the above geometrical results in order to analyze the asymptotic behavior of an elastic rod when its thickness tends to $0$. As usual, we consider a elastic energy density which is bounded below by $dist^2\big(F,SO(3)\big)$ (see e.g. [7], [11], [12], [19] and [20]). As mentioned in the introduction, our decomposition of the deformation permit us to scale the applied forces in order to obtain estimates on the deformation and on the total elastic energy (see (IV.1.8)).  Then, to simplify the argument,  we specify the energy density through choosing a St Venant-Kirchhoff's material (see (IV.1.9)). In the following, we derive the limit elastic energy in the two cases $\kappa=2$ and $\kappa>2$ using $\Gamma$-convergence techniques. The limit energy is expressed as a functional of  the fields ${\cal V}$, $\GR$, ${\cal V}_S$ and $\overline {v}$. Such limit energies depending on more variables than the 3D ones have also be derived in [19] by different techniques. In the present paper we also eliminate the fields ${\cal V}_S$ and $\overline {v}$ to be in a position to obtain a minimization problem for the rotation field $\GR$ and the field ${\cal V}$.  Let us emphazise that in the $\Gamma$-limit procedure, the decomposition of the deformations is again helpful in two directions. Firstly it provides an explicit expression of the limit Green-St Venant deformation  tensor and secondly it simplifies the proof of the two  conditions involved in the identification of the limit energy by $\Gamma$-convergence.
\medskip
\noindent {\Ggras  IV.1  Assumption  on the forces}
\medskip
In this part we assume that the curved rod ${\cal P}_\delta$ is made of an elastic material. As in [5] and [11] we assume that the elastic energy $W$ satisfy (actually we will consider an explicit energy)
$$\forall \GF \in \GM_3, \qquad W(\GF)\ge C\,\hbox{dist}^2\big(\GF,SO(3)\big),\leqno(IV.1.1)$$ where $C$ is a strictly positive constant.

Let us denote by $f_\delta\in(L^2(\Omega_\delta))^3$ the applied forces and by $J(\phi)$ the total energy 
$$J(\phi)=\int_{{\cal P}_\delta}W(\nabla \phi)-\int_{{\cal P}_\delta}f_\delta\cdot \phi\leqno(IV.1.2)$$ This energy is considered over the set of  admissible deformations:
$$\GU_\delta=\Bigl\{\phi\in (H^1({\cal P}_\delta))^3\; |\; \phi=I_d\quad \hbox{on}\quad \Gamma_{0,\delta}\Big\}.\leqno(IV.1.3)$$ 
For different boundary conditions see Subsection II.2.4.
\noindent As far as the behavior of the forces $f_\delta$ is concerned we split the forces into two parts. The first one does not depends on the variables $(s_1,s_2)$ and the second part has  a resultant equal to $0$. Due to  estimates (II.3.4), (II.3.6)  and (II.3.7) the admissible order of theses forces  can be chosen different. 
\medskip
\noindent  Let   $f$ be in $(L^2(0,L))^3$ and let $g $ be in  $(L^2(\Omega))^3$ such that
$$ \int_\omega g (S_1,S_2,s_3)dS_1dS_2=0\qquad \hbox{for a.e.}\enskip s_3\in ]0,L[.\leqno(IV.1.4)$$ We assume that $f_\delta$ is defined by
$$f_\delta(s)=\delta^\kappa f(s_3)+\delta^{\kappa-1}g \Big({s_1\over \delta},{s_2\over \delta},s_3\Big)\qquad \hbox{for a.e.}\enskip s\in \Omega_\delta.\leqno(IV.1.5)$$ The fact remains that to find a minimizer or to find a deformation that approaches  the minimizer of $J(\phi)$ or of $J(\phi)-J(I_d)$ is the same.  Let $v$ be in $\GU_\delta$, thanks to (II.3.7), (II.3.9), (IV.1.4) and (IV.1.5), we obtain 
$$\Big|\int_{{\cal P}_\delta}f_\delta\cdot(v-I_d)\Big| \le C\delta^\kappa(||f_{(L^2(0,L))^3}+
||g ||_{(L^2(\Omega))^3})||\hbox{dist}(\nabla v, SO(3))||_{L^2({\cal P}_\delta)}. \leqno(IV.1.6)$$ 
Actually one can think to use estimate (II.3.4) instead of (II.3.7) and (II.3.9) in the above inequality. The reader will easily see that this gives a better estimate only in the case where $\kappa<2$ which is not considered in the following (see Remark below).

 It is well known that generally a minimizer of $J$ does not exist  on $\GU_\delta$. In the next sections we will investigate the behavior of the functional $\ds {1\over \delta^{2\kappa}}\big(J(\phi)-J(I_d)\big)$ in the framework of the $\Gamma$-convergence. Hence, we assume that 
 $$ {1\over \delta^{2\kappa}}\big(J(v)-J(I_d)\big)\le C_1 \leqno(IV.1.7)$$ where $C_1 $ does not depend on $\delta$. Using (IV.1.1) and (IV.1.6) we obtain for such $v$
$$C||\hbox{dist}(\nabla v, SO(3))||^2_{L^2({\cal P}_\delta)}-C\delta^\kappa(||f||_{(L^2(\Omega))^3}+
||g ||_{(L^2(\Omega))^3})||\hbox{dist}(\nabla v, SO(3))||_{L^2({\cal P}_\delta)}\le C_1 \delta^{2\kappa}.$$  Hence,  we have
$$||\hbox{dist}(\nabla v, SO(3))||_{L^2({\cal P}_\delta)}\le C\delta^\kappa.\leqno(IV.1.8)$$ where the constant  $C$ depends on the sum $||f||_{(L^2(0,L))^3}+ ||g ||_{(L^2(\Omega))^3}$ and of $C_1$. 
\smallskip
\noindent {\ggras Remark.} 
If one uses estimates (II.3.4) to bound the contribution of the forces in the energy, one alternatively obtains through similar calculations
$$||\hbox{dist}(\nabla v, SO(3))||_{L^2({\cal P}_\delta)}\le C\delta^{1+\kappa/2}.$$
Comparing to (IV.1.8), one gets a better estimate only if $\kappa<2$.
\medskip
 Let us notice that once the assumption (IV.1.5) on the applied forces is adopted, the estimate (IV.1.8) and the results of Section II permit to obtain estimates  of ${\cal V}$, $\GR$, $\overline{v}$ and $\nabla_x v-\GR$ with respect to $\delta$.  To emphasize how these estimates can help pass to the limit as $\delta$ tends to $0$, we will restrict the following analysis to a classical and simple elastic energy.  
 We denote by $tr(\GA)$  the sum of the elements on the main diagonal  of the $3\times 3$ matrix $\GA$.

{\it In order to simplify the derivation of the limit model we choose $$W(\GF)=\left\{\eqalign{
&{\lambda\over 8}\big(tr(\GF^T\GF-\GI_3) \big)^2+{\mu\over 4}tr\big((\GF^T\GF-\GI_3)^2\big)\quad\hbox{if}\quad \det(\GF)>0\cr
&+\infty\qquad \hbox{if}\qquad \det(\GF)\le 0\cr}\right.\leqno(IV.1.9)$$ which corresponds to a St Venant-Kirchhoff's material (see [7], [8]).  }
\smallskip
For all $3\times 3$ matrices such that $\det(\GF)>0$ we have $$W(\GF)\ge {\mu\over 4}|||\GF^T\GF-\GI_3|||^2\ge{\mu\over 4}\hbox{dist}^2(\GF,SO(3))$$ hence assumption (IV.1.1) is satisfied. For every $\phi\in \GU_\delta$ satisfying (IV.1.7), we have using (IV.1.6)
$$\eqalign{
{\mu\over 4}||\nabla\phi^T\nabla\phi-\GI_3||^2_{(L^2({\cal P}_\delta))^{3\times 3}} & \le J(\phi)-J(I_d)+\int_{{\cal P}_\delta}f_\delta\cdot( \phi-I_d)\cr
&\le C_1 \delta^{2\kappa}+C\delta^\kappa( ||f||_{(L^2(0,L))^3}+||g||_{(L^2(\Omega))^3})||\hbox{dist}(\nabla\phi, SO(3))||_{L^2({\cal P}_\delta)}.\cr}$$ Due to estimate (IV.1.8) we obtain the following estimate of the Green-St Venant's tensor:
$$\big\|{1\over 2}\big\{\nabla\phi^T\nabla\phi-\GI_3\big\}\big\|_{(L^2({\cal P}_\delta))^{3\times 3}}\le  C\delta^\kappa.\leqno(IV.1.10)$$ It results that $\phi$ belongs to $(W^{1,4}({\cal P}_\delta))^3$ and moreover
$$||\nabla\phi||_{(L^4({\cal P}_\delta))^{3\times 3}}\le  C\delta^{1\over 2}.\leqno(IV.1.11)$$ Furthermore, there exists two strictly positive constants $c$ and $C$ which does not depend  on $\delta$  such that for any $\phi\in \GU_\delta$ satisfying (IV.1.7) we have
$$-c\delta^{2\kappa}\le J(\phi)-J(I_d)\le C \delta^{2\kappa}.\leqno(IV.1.12)$$  We set
$$m_\delta=\inf_{\phi\in \GU_\delta}\big(J(\phi)-J(I_d)\big).\leqno(IV.1.13)$$
As a consequence of the  inequality in (IV.1.12) we  have 
$$-c\le {m_\delta\over \delta^{2\kappa}}\le 0.\leqno(IV.1.14)$$ We denote
$$m_\kappa=\liminf_{\delta\to 0}{m_\delta\over \delta^{2\kappa}}.\leqno(IV.1.15) $$

\noindent {\Ggras  IV.2   Limit  model in the case $\kappa=2$}
\medskip
Let   $\big(v_\delta\big)_{0<\delta\le \delta_0 }$  be  a sequence of deformations belonging to $\GU_\delta$ and such that
$$\liminf_{\delta\to 0}{J(v_\delta)-J(I_d)\over \delta^4}<+\infty.\leqno(IV.2.1)$$ Upon extracting a subsequence (still indexed by $\delta$) we can assume that 
 the sequence $(v_\delta)$ satisfies the condition (IV.1.7). From the estimates of the previous section we obtain
$$\left\{\eqalign{
&||\hbox{dist}(\nabla v_\delta, SO(3))||_{L^2({\cal P}_\delta)}\le C\delta^2,\cr
& \big\|{1\over 2}\big\{\nabla v_\delta^T\nabla v_\delta-\GI_3\big\}\big\|_{(L^2({\cal P}_\delta))^{3\times 3}}\le  C\delta^2,\cr
 &||\nabla v_\delta||_{(L^4({\cal P}_\delta))^{3\times 3}}\le  C\delta^{1\over 2}\cr}\right.\leqno(IV.2.3)$$ For any fixed $\delta \in (0,\delta_0]$, the deformation $v_\delta$ is decomposed following (II.2.1) in such a way that Theorem II.2.2 is satisfied.  There exists a subsequence still indexed by $\delta$ such that  (see Section II.6)
$$\left\{\eqalign{
&\GR_\delta \rightharpoonup \GR  \quad \hbox{weakly in}\quad \big(H^1(0,L)\big)^{3\times 3}\cr
&{\cal V}_\delta\longrightarrow   {\cal V}  \quad \hbox{strongly in}\quad \big(H^1(0,L)\big)^3\cr
&{\cal V}_{B,\delta}\longrightarrow   {\cal V}  \quad \hbox{strongly in}\quad \big(H^1(0,L)\big)^3\cr
&{1\over \delta}{\cal V}_{S,\delta}\rightharpoonup {\cal V}_{S} \quad \hbox{weakly in}\quad \big(H^1(0,L)\big)^3\cr
&{1\over \delta^2}\Pi _\delta\overline {v}_\delta \rightharpoonup  \overline {v}  \quad \hbox{weakly in}\quad \big(L^2(0,L;H^1(\omega))\big)^3\cr}\right.\leqno(IV.2.4)$$ where $\GR (s_3)$ belongs to $SO(3)$ for any $s_3 \in [0,L]$, ${\cal V} \in \big(H^2(0,L)\big)^3 $ together with
$${\cal V} (0)=M(0), \qquad  \GR (0)=\GI_3,\qquad {\cal V}_{S} (0)=0,\qquad\hbox{and} \qquad{d{\cal V} \over ds_3}=\GR  \Gt .\leqno (IV.2.5)$$ 
\noindent Furthermore, we also have (see (III.2.6) and (III.2.7))
$$\left\{\eqalign{
&\Pi _\delta v_\delta \longrightarrow    {\cal V}  \quad \hbox{strongly in}\quad \big(H^1(\Omega)\big)^3,\cr
&\Pi _\delta (\nabla_x v_\delta)  \longrightarrow  \GR \quad \hbox{strongly in}\quad \big(L^2(\Omega)\big)^{3\times 3}.\cr
&{1\over 2\delta}\Pi _\delta \big((\nabla_xv_\delta)^T\nabla_x v_\delta-\GI_3\big) \rightharpoonup   \GE  \qquad\hbox{weakly in}\quad (L^2(\Omega))^{3\times 3},\cr}\right.\leqno(IV.2.6)$$ where 
$$\eqalign{
\GE =&{1\over 2}\Big\{(\Gn_1\,|\, \Gn_2\,|\, \Gt) \Bigl({\partial\overline{v}\over \partial S_1}\,|\, {\partial\overline{v} \over \partial S_2}\,|\,{d\GR \over ds_3}\big(S_1\Gn_1+S_2\Gn_2)+{d{\cal V}_{S} \over ds_3}\Big)^T\GR \cr
+& \GR^T  \Bigl({\partial\overline{v}\over \partial S_1}\,|\, {\partial\overline{v} \over \partial S_2}\,|\,{d\GR \over ds_3}\big(S_1\Gn_1+S_2\Gn_2)+{d{\cal V}_{S} \over ds_3}\Big) (\Gn_1\,|\, \Gn_2\,|\, \Gt)^T\Big\}.\cr}$$ Remark that, due to the decomposition (II.2.1), the convergences (IV.2.4) and (IV.2.6) imply that
$${\Pi_\delta(v_\delta-{\cal V}_\delta)\over \delta}\longrightarrow S_1(\GR-\GI_3)\Gn_1+S_2(\GR-\GI_3)\Gn_2\quad \hbox{strongly in}\quad \big(L^2(\Omega)\big)^3.\leqno (IV.2.7) $$ Now, recall that 
$$\left\{\eqalign{
&{J(v_\delta)-J(I_d)\over \delta^4}\cr
=&\int_{\Omega} \Big\{{\lambda\over 2}\Big[tr\Big({1\over 2\delta}\Pi _\delta \big((\nabla_xv_\delta)^T\nabla_x v_\delta-\GI_3\big)\Big) \Big]^2+\mu |||{1\over 2\delta}\Pi _\delta \big((\nabla_xv_\delta)^T\nabla_x v_\delta-\GI_3\big) |||^2\Big\}|\Pi_\delta\det(\nabla\Phi)|\cr
-&\int_\Omega f \cdot\Pi_\delta(v_\delta-I_d)|\Pi_\delta\det(\nabla\Phi)|-\int_\Omega g \cdot{\Pi_\delta(v_\delta-I_d)\over \delta}|\Pi_\delta\det(\nabla\Phi)|\cr}\right.\leqno(IV.2.8)$$  In order to obtain the limit of the terms involving the forces, we recall that $\ds \det(\nabla\Phi)=1+s_1\det\Big(\Gn_1\; |\; \Gn_2\;|\;{d\Gn_1\over ds_3}\Big)+s_2\det\Big(\Gn_1\; |\; \Gn_2\;|\;{d\Gn_2\over ds_3}\Big)$ so that indeed $\Pi_\delta\det(\nabla\Phi)$ strongly converges to $1$ in $L^\infty(\Omega)$ as $\delta$ tends to $0$. As a consequence and using the convergences (IV.2.6) and (IV.2.7), it follows that
$$\eqalign{
&\lim_{\delta \to 0}\Big(\int_\Omega f \cdot\Pi_\delta(v_\delta-I_d)|\Pi_\delta\det(\nabla\Phi)|+\int_\Omega g \cdot{\Pi_\delta(v_\delta-I_d)\over \delta}|\Pi_\delta\det(\nabla\Phi)|\Big)\cr
=&\int_0^L\Big(|\omega|f +\sum_{\alpha=1}^2\int_\omega g  S_\alpha \det\big(\Gn_1\; |\; \Gn_2\;|\;{d\Gn_\alpha\over ds_3}\big)dS_1dS_2\Big)\cdot({\cal V} -M) 
+ \sum_{\alpha=1}^2\int_0^L\Big(\int_\omega g S_\alpha dS_1dS_2\Big)\cdot(\GR-\GI_3)\Gn_\alpha.\cr}$$
In order to pass to the limit-inf in the left hand side of (IV.2.8), we recall that the map $\ds M\longmapsto {\lambda\over 2}(tr(M))^2+\mu |||M|||^2$ is continuous and convex from $\GM_3$ into $\R$, so that the map $A\longmapsto \ds \int_{\Omega} \big({\lambda\over 2}(tr(A))^2+\mu |||A|||^2\big)$ from $(L^2(\Omega))^{3\times 3} $ into $\R$ is weak lower semi-continuous. The above strong convergence of $\Pi_\delta\det(\nabla\Phi)$ together with convergences (IV.2.6) finally give
$$\left\{\eqalign{
&\int_\Omega \Big\{{\lambda\over 2}(tr(\GE ))^2+\mu |||\GE |||^2\Big\}\cr
-&\int_0^L\Big(|\omega| f+\sum_{\alpha=1}^2\int_\omega g  S_\alpha \det\big(\Gn_1\; |\; \Gn_2\;|\;{d\Gn_\alpha\over ds_3}\big)dS_1dS_2\Big)\cdot({\cal V} -M)\cr
-&\sum_{\alpha=1}^2\int_0^L\Big(\int_\omega g  S_\alpha dS_1dS_2\Big)\cdot(\GR-\GI_3)\Gn_\alpha\le \liminf_{\delta\to 0}{1\over \delta^4}\big(J(v_\delta)-J(I_d)\big).\cr}
\right.\leqno(IV.2.9)$$

Let $\GU_{nlin}$ be the set
$$\eqalign{
\GU_{nlin}=\Bigl\{ &\big({\cal V}^{'},\GR^{'},{\cal V}^{'}_S,\overline{v}^{'}\big)\in (H^2(0,L))^3\times (H^1(0,L))^{3\times 3}\times (H^1(0,L))^3\times (L^2(0,L;H^1(\omega)))^3\;\; |\;\; \cr
&{\cal V}^{'}(0)=M(0),\quad {\cal V}^{'}_S(0)=0, \quad \GR^{'}(0)=\GI_3,\quad\int_\omega \overline{v}^{'}(S_1,S_2,s_3)dS_1dS_2=0\quad \hbox{for a.e.}\enskip s_3\in (0,L)\cr
&\GR^{'}(s_3)\in SO(3)\; \hbox{for any}\enskip s_3\in [0,L],\quad {d{\cal V}^{'}\over ds_3}=\GR^{'}\Gt\Big\}.\cr}$$ The set $\GU_{nlin}$ is closed in the product space. For any $\big({\cal V}^{'},\GR^{'},{\cal V}^{'}_S,\overline{v}^{'}\big)\in \GU_{nlin}$, we denote by ${\cal J}_{NL}$ the following limit energy
$$\left\{\eqalign{
&{\cal J}_{NL}\big( {\cal V}^{'},\GR^{'},{\cal V}^{'}_S,\overline{v}^{'}\big)=\int_\Omega \Big\{{\lambda\over 2}(tr(\GE^{'}))^2+\mu |||\GE^{'}|||^2\Big\}\cr
&-\int^L_0\Big(|\omega |f +\sum_{\alpha=1}^2 \int_\omega g S_\alpha \det\big(\Gn_1\; |\; \Gn_2\;|\;{d\Gn_\alpha\over ds_3}\big)dS_1dS_2\Big)
\cdot({\cal V}^{'}-M)\cr
&- \sum_{\alpha=1}^2\int_0^L\Big(\int_\omega g S_\alpha dS_1dS_2\Big)\cdot(\GR^{'}-\GI_3)\Gn_\alpha\cr}\right.\leqno(IV.2.10)$$ where
$$\left\{\eqalign{
\GE^{'}=&{1\over 2}\Big\{(\Gn_1\,|\, \Gn_2\,|\, \Gt) \Bigl({\partial\overline{v}^{'}\over \partial S_1}\,|\, {\partial\overline{v}^{'}\over \partial S_2}\,|\,{d\GR^{'}\over ds_3}\big(S_1\Gn_1+S_2\Gn_2)+{d{\cal V}^{'}_S\over ds_3}\Big)^T\GR^{'}\cr
+& \GR^{'T} \Bigl({\partial\overline{v}^{'}\over \partial S_1}\,|\, {\partial\overline{v}\over \partial S_2}\,|\,{d\GR^{'}\over ds_3}\big(S_1\Gn_1+S_2\Gn_2)+{d{\cal V}^{'}_S\over ds_3}\Big) (\Gn_1\,|\, \Gn_2\,|\, \Gt)^T\Big\}.\cr}\right.\leqno(IV.2.11)$$ With this notation  (IV.2.9)  reads as
$${\cal J}_{NL}\big({\cal V} ,\GR ,{\cal V}_{S} ,\overline{v} \big)\le\liminf_{\delta\to 0}{J(v_\delta)-J(I_d)\over \delta^4}.\leqno(IV.2.12)$$
\medskip
 Now let $\big({\cal V},\GR,{\cal V}_S,\overline{v}\big)$ be in $\GU_{nlin}$ and let $\Big(\big({\cal V}_N,\GR_N,{\cal V}_{S,N},\overline{v}_N\big)\Big)_{N\in \N^*}$ be a sequence  of elements belonging to $\GU_{nlin}$ such that
$$\left\{\eqalign{
{\cal V}_N\in (W^{2,\infty}(0,L))^{3},\qquad {\cal V}_N&\longrightarrow {\cal V}\quad \hbox{strongly in }\quad (H^2(0,L))^{3}\cr
\GR_N\in (W^{1,\infty}(0,L))^{3\times 3},\qquad \GR_N&\longrightarrow \GR\quad \hbox{strongly in }\quad (H^1(0,L))^{3\times 3}\cr
{\cal V}_{S,N}\in (W^{1,\infty}(0,L))^{3},\qquad {\cal V}_{S,N}&\longrightarrow {\cal V}_S\quad \hbox{strongly in }\quad (H^1(0,L))^{3}\cr
\overline{v}_N\in (W^{1,\infty}(\Omega))^{3},\enskip  \overline{v}_N(S_1,S_2,0)&=0,\; \hbox{ for a.e. }\; (S_1,S_2)\in \omega,\cr
 \overline{v}_N&\longrightarrow \overline{v}\quad \hbox{strongly in }\quad (L^2(0,L;H^1(\omega))^{3}.\cr}\right.\leqno(IV.2.13 )$$ To prove the existence of the sequence $\big(\GR_N\big)_{N\in \N^*}$ see the appendix at the end of the paper.
\smallskip
We consider the deformations ($\delta\in (0,\delta_0]$)
$$v_{N,\delta}(s)={\cal V}_N(s_3)+\GR_N(s_3)(s_1\Gn_1+s_2\Gn_2)+\delta{\cal V}_{S,N}(s_3)+\delta^2\overline{v}_N\Big({s_1\over \delta},{s_2\over \delta},s_3\Big),\qquad s\in \Omega_\delta.\leqno(IV.2.14 )$$ Using  convergences (IV.2.13), the fact that $\big({\cal V}_N,\GR_N,{\cal V}_{S,N},\overline{v}_N\big)$ belongs to $\GU_{nlin}$ and proceeding as in Subsection III.2, we have
$$\left\{\eqalign{
&\Pi _\delta v_{N,\delta} \longrightarrow    {\cal V}_N \quad \hbox{strongly in}\quad \big(W^{1,\infty}(\Omega)\big)^3,\cr
&\Pi _\delta (\nabla_x v_{N,\delta})  \longrightarrow  \GR_N\quad \hbox{strongly in}\quad \big(L^\infty(\Omega)\big)^{3\times 3},\cr
&{1\over 2\delta}\Pi _\delta \big((\nabla_xv_{N,\delta})^T\nabla_x v_{N,\delta}-\GI_3\big) \longrightarrow \GE_N \qquad\hbox{strongly in}\quad (L^\infty(\Omega))^{3\times 3},\cr}\right.\leqno(IV.2.15 )$$ where 
$$\eqalign{
\GE_N=&{1\over 2}\Big\{(\Gn_1\,|\, \Gn_2\,|\, \Gt) \Bigl({\partial\overline{v}_N\over \partial S_1}\,|\, {\partial\overline{v}_N\over \partial S_2}\,|\,{d\GR_N\over ds_3}\big(S_1\Gn_1+S_2\Gn_2)+{d{\cal V}_{S,N}\over ds_3}\Big)^T\GR_N\cr
+& \GR^T_N \Bigl({\partial\overline{v}_N\over \partial S_1}\,|\, {\partial\overline{v}_N \over \partial S_2}\,|\,{d\GR_N\over ds_3}\big(S_1\Gn_1+S_2\Gn_2)+{d{\cal V}_{S,N}\over ds_3}\Big) (\Gn_1\,|\, \Gn_2\,|\, \Gt)^T\Big\}\cr}$$ If $\delta$ is sufficiently  small we have $\det\big(\nabla_x v_{N,\delta}(x)\big)>0$ for a.e. $x\in  {\cal P}_\delta$ because of the second convergence in (IV.2.15). It follows that  $J(v_{N,\delta})<+\infty$.  We divide $J(v_{N,\delta})-J(I_d)$ by $\delta^4$ and we pass to the limit using the strong convergences  (IV.2.15). We obtain
$$\lim_{\delta\to 0}{1\over \delta^4}\big(J(v_{N,\delta})-J(I_d)\big)={\cal J}_{NL}\big({\cal  V}_N,\GR_N,{\cal V}_{N,S},\overline{v}_N\big).\leqno(IV.2.16 )$$ Letting $N$ tend to $+\infty$ and using (IV.2.13), it follows that for any $\big({\cal V},\GR,{\cal V}_S,\overline{v}\big)\in \GU_{nlin}$
$${\cal J}_{NL}\big({\cal V},\GR,{\cal V}_S,\overline{v}\big)=\lim_{N\to+\infty}{\cal J}_{NL}\big({\cal  V}_N,\GR_N,{\cal V}_{N,S},\overline{v}_N\big).\leqno(IV.2.17 )$$ Hence, through a standard diagonal process for any $\big({\cal V},\GR,{\cal V}_{S},\overline{v}\big)\in \GU_{nlin}$ there exists a sequence of admissible deformations $v_\delta\in(H^1({\cal P}_\delta))^3$ such that
$${\cal J}_{NL}\big({\cal V},\GR,{\cal V}_S,\overline{v}\big)=\lim_{\delta\to 0}{J(v_\delta)-J(I_d)\over \delta^4}.\leqno(IV.2.18 )$$ The following theorem summarizes the above results.
\smallskip
\noindent{\ggras  Theorem IV.2.1. }{\it  The functional ${\cal J}_{NL}$ is the $\Gamma$-limit of $\ds {J(.)-J(I_d)\over \delta^4}$ in the following sense:

$\bullet$ consider any sequence of deformations  $\big(v_\delta\big)_{0<\delta\le \delta_0 }$   belonging to $\GU_\delta$ and satisfying 
$$\ds \liminf_{\delta\to 0}{J(v_\delta)-J(I_d)\over \delta^4}<+\infty$$ and let $\big({\cal V}_\delta,\GR_\delta,{\cal V}_{S,\delta},\overline{v}_\delta\big)$ be the terms of the decomposition of $v_\delta$ given by Theorem II.2.2.   Then there exists $\big({\cal V},\GR,{\cal V}_S,\overline{v}\big)\in \GU_{nlin}$ such that (up to a subsequence )
$$\eqalign{
&\GR_\delta \rightharpoonup \GR  \quad \hbox{weakly in}\quad \big(H^1(0,L)\big)^{3\times 3}\cr
&{\cal V}_\delta\longrightarrow   {\cal V}  \quad \hbox{strongly in}\quad \big(H^1(0,L)\big)^3\cr
&{\cal V}_{B,\delta}\longrightarrow   {\cal V}  \quad \hbox{strongly in}\quad \big(H^1(0,L)\big)^3\cr
&{1\over \delta}{\cal V}_{S,\delta}\rightharpoonup {\cal V}_{S} \quad \hbox{weakly in}\quad \big(H^1(0,L)\big)^3\cr
&{1\over \delta^2}\Pi _\delta\overline {v}_\delta \rightharpoonup  \overline {v}  \quad \hbox{weakly in}\quad \big(L^2(0,L;H^1(\omega))\big)^3\cr}$$ and we have
$${\cal J}_{NL}\big({\cal V},\GR,{\cal V}_S,\overline{v}\big)\le \liminf_{\delta\to 0}{J(v_\delta)-J(I_d)\over \delta^4}$$

$\bullet$ for any $\big({\cal V},\GR,{\cal V}_S,\overline{v}\big)\in \GU_{nlin}$ there exists a sequence $\big(v_\delta\big)_{0<\delta\le \delta_0 }$   belonging to $\GU_\delta$such that
$${\cal J}_{NL}\big({\cal V},\GR,{\cal V}_S,\overline{v}\big)= \lim_{\delta\to 0}{J(v_\delta)-J(I_d)\over \delta^4}.$$ Moreover, there exists $\big({\cal V}_0 ,\GR_0 ,{\cal V}_{S,0} ,\overline{v}_0 
\big)\in \GU_{nlin}$ such that
$$m_2=\lim_{\delta\to 0}{m_\delta\over \delta^{4}}={\cal J}_{NL}\big({\cal V}_0,\GR_0,{\cal V}_{S,0} ,\overline{v}_0\big)=\min_{\big({\cal V},\GR,{\cal V}_{S},\overline{v}\big)\in\GU_{nlin}}{\cal J}_{NL}\big({\cal V},\GR,{\cal V}_{S},\overline{v} \big).\leqno(IV.2.19)$$ } 
\medskip
The next theorem shows that the variables ${\cal V}_S$ and $\overline{v}$  can be eliminated in the minimization problem (IV.2.19). To this end let us first introduce a few notations. We denote by $E$ the Young's modulus of the material and by $\chi$ the solution of the following torsion problem:
 $$\left\{\eqalign{
&\chi\in H^1(\omega),\qquad \int_\omega\chi=0\cr
&\int_\omega\nabla\chi\nabla\psi=-\int_\omega\Bigl\{-S_2{\partial\psi\over \partial S_1}+S_1{\partial\psi\over \partial S_2}\Big\}\cr
&\forall \psi\in H^1(\omega).\cr}\right.\leqno(IV.2.22)$$
At least we set 
$$K=\int_\omega\Big[\Big({\partial\chi\over \partial S_1}-S_2\Big)^2+\Big({\partial\chi\over \partial S_2}+S_1\Big)^2\Big],\qquad I_1=\int_\omega S_1^2,\qquad I_2=\int_\omega S_2^2.$$
\smallskip
\noindent{\ggras  Theorem  IV.2.2 }{\it Let $({\cal V}_0,\GR_0)$ be given by Theorem IV.2.1. The minimum $m_2$ of the functional ${\cal J}_{NL}$ over $\GU_{nlin}$ satisfies the following minimization problem:
$$m_2={\cal F}_{NL}\big({\cal V}_0,\GR_0\big)=\min_{\big({\cal V},\GR\big)\in\GV_{nlin}}{\cal F}_{NL}
\big({\cal V},\GR\big), \leqno(IV.2.20) $$ where 
$$\eqalign{
\GV_{nlin}=\Bigl\{ &\big({\cal V},\GR\big)\in (H^2(0,L))^3\times (H^1(0,L))^{3\times 3}\;\; |\;\; 
{\cal V}(0)=M(0),\quad \GR(0)=\GI_3,\cr
&\GR (s_3)\in SO(3)\; \hbox{for any}\; s_3\in [0,L],\quad {d{\cal V} \over ds_3}=\GR \Gt\Big\},\cr}$$ and  
$$\left\{\eqalign{
& {\cal F}_{NL}\big({\cal V},\GR\big)={EI_1\over 2}\int_0^L\Big({d\GR\over ds_3}\Gt\cdot \GR \Gn_1\Big)^2+{EI_2\over 2}\int_0^L\Big({d\GR\over ds_3}\Gt\cdot \GR \Gn_2\Big)^2+{\mu K\over 4}  \int_0^L\Big({d\GR\over ds_3}\Gn_1\cdot \GR \Gn_2\Big)^2\cr
&-\int_0^L\Big(|\omega |f +\sum_{\alpha=1}^2  \int_\omega g  S_\alpha \det\big(\Gn_1\; |\; \Gn_2\;|\;{d\Gn_\alpha\over ds_3}\big)dS_1dS_2\Big)\cdot({\cal V} -M)\cr
&-\sum_{\alpha=1}^2 \int_0^L\Big(\int_\omega g S_\alpha dS_1dS_2\Big)\cdot(\GR-\GI_3)\Gn_\alpha.\cr}\right.\leqno(IV.2.21)$$ }
\smallskip
\noindent{\ggras Proof of Theorem IV.2.2. } 
Let us first notice that in the expression (IV.2.10) of ${\cal J}_{NL}({\cal V},\GR,{\cal V}_S, \overline v)$, one can replace $\GE$ by $\widehat \GE$ where $\GE$ and $\widehat \GE$ are given by (III.2.7), (III.2.8) and (III.2.9).

In order to eliminate $({\cal V}_S, \overline v)$, we fix  $\big({\cal V} ,\GR \big)\in \GV_{nlin}$  and we minimize the functional ${\cal J}_{NL}\big({\cal V} ,\GR ,\cdot,\cdot\big)$ over the space
$$\eqalign{
\GW=\Bigl\{ &\big({\cal V} _S,\overline{v} \big)\in (H^1(0,L))^3\times (L^2(0,L;H^1(\omega)))^3\;\; |\;\;  {\cal V} _S(0)=0, \cr
&\int_\omega \overline{v} (S_1,S_2,s_3)dS_1dS_2=0\quad \hbox{for a.e.} \enskip s_3\in (0,L)\Big\}.\cr}$$  Through  solving simple variational problems    (see [14]), we find that the minimum of the functional ${\cal J}_{NL}\big({\cal V} ,\GR ,\cdot,\cdot\big)$ over the space $\GW$ is obtained with  $\ds{d{\cal V} _S\over ds_3}\cdot\GR \Gt=0$ and
$$\left\{\eqalign{
&\overline{v} (S_1,S_2,.)\cdot\GR \Gn_1 = \nu\Big\{{S_1^2-S_2^2\over 2}{d\GR \over ds_3}\Gt\cdot \GR  \Gn_1+S_1S_2{d\GR \over ds_3}\Gt\cdot \GR  \Gn_2\Big\}\cr
&\overline{v} (S_1,S_2,.)\cdot\GR \Gn_2= \nu\Big\{S_1S_2{d\GR \over ds_3}\Gt\cdot\GR \Gn_1+{S_2^2-S_1^2\over 2}{d\GR \over ds_3}\Gt\cdot \GR  \Gn_2\Big\}\cr
&\overline{v} (S_1,S_2,.)\cdot\GR \Gt+S_1{d{\cal V}_S\over ds_3}\cdot\GR\Gn_1 +S_2{d{\cal V}_S\over ds_3}\cdot\GR\Gn_2 =\big\{\chi(S_1,S_2){d\GR \over ds_3}\Gn_1\cdot \GR  \Gn_2\big\}\GR \Gt \cr}\right.\leqno(IV.2.23)$$ where $\ds \nu={\lambda\over 2(\lambda+\mu)}$ is the Poisson's coefficient of the material. Then the symmetric tensor $\widehat{\GE}$ (see again (III.2.9)) at the minimum is given by
$$\widehat{\GE}(\GR)=\pmatrix{
\ds-\nu\widehat{\GE}_{33}(\GR) & 0 & \ds {1\over 2}\Big({\partial\chi\over \partial S_1}-S_2\Big){d\GR\over ds_3}\Gn_1\cdot \GR\Gn_2 \cr
* & -\nu\widehat{\GE}_{33}(\GR)  &  \ds {1\over 2}\Big({\partial\chi\over \partial S_1}+S_1\Big){d\GR\over ds_3}\Gn_1\cdot \GR\Gn_2\cr
 * & *  & \ds  \widehat{\GE}_{33}(\GR)\cr},\leqno(IV.2.24)$$ where $\ds \widehat{\GE}_{33}(\GR)=-S_1{d\GR\over ds_3}\Gt\cdot \GR\Gn_1-S_2{d\GR\over ds_3}\Gt\cdot \GR\Gn_2$. Upon replacing $\widehat \GE$ by $\widehat{\GE}(\GR)$ in the expression of ${\cal J}_{NL}$ we obtain 
$$\min_{\big({\cal V} _S,\overline{v} \big)\in\GW}{\cal J}_{NL}\big({\cal V} ,\GR ,{\cal V} _S,\overline{v} \big)={\cal F}_{NL}\big({\cal V} ,\GR \big),$$ where the functional ${\cal F}_{NL}$ is given by (IV.2.21).\fin
\noindent {\ggras Remark.}  The above analysis shows that if $\big(v_\delta\big)_{0<\delta\le \delta_0}$ is a sequence such that 
$$m_2=\lim_{\delta\to 0}{J(v_\delta)-J(I_d)\over \delta^{4}},$$
then there exists a subsequence and $({\cal V}_0,\GR_0)\in \GV_{nlin}$, which is a solution of Problem (IV.2.20), such that  the sequence of the Green-St Venant's deformation tensors satisfies
 $${1\over 2\delta}\Pi _\delta \big((\nabla_xv_\delta)^T\nabla_x v_\delta-\GI_3\big) \longrightarrow   (\Gn_1\,|\, \Gn_2\,|\, \Gt)\widehat{\GE}(\GR_0) (\Gn_1\,|\, \Gn_2\,|\, \Gt)^T  \qquad\hbox{strongly in}\quad (L^2(\Omega))^{3\times 3},$$
 where $\widehat{\GE}(\GR_0)$ is defined in (IV.2.24). \medskip
\noindent {\Ggras  IV.3   Limit  model in the case $\kappa>2$}
\medskip
Let   $\big(v_\delta\big)_{0<\delta\le \delta_0}$  be  a sequence of deformations belonging to $\GU_\delta$ and such that
$$\liminf_{\delta\to 0}{J(v_\delta)-J(I_d)\over \delta^{2\kappa}}<+\infty.\leqno(IV.3.1)$$ Upon extracting a subsequence (still indexed by $\delta$) we can assume that the sequence $(v_\delta)$ satisfies the condition (IV.1.7). From the estimates of the  section IV.1 we obtain
$$\left\{\eqalign{
&||\hbox{dist}(\nabla v_\delta, SO(3))||_{L^2({\cal P}_\delta)}\le C\delta^\kappa,\cr
& \big\|{1\over 2}\big\{\nabla v_\delta^T\nabla v_\delta-\GI_3\big\}\big\|_{(L^2({\cal P}_\delta))^{3\times 3}}\le  C\delta^\kappa,\cr
 &||\nabla v_\delta||_{(L^4({\cal P}_\delta))^{3\times 3}}\le  C\delta^{1\over 2}.\cr}\right.\leqno(IV.3.3)$$  For any fixed $\delta \in (0,\delta_0]$, the displacement  $u_\delta=v_\delta-I_d$ is decomposed following (II.2.1) and (III.3.2) in such a way that Theorem II.2.2 is satisfied. There exists a subsequence still indexed by $\delta$ such that  (see Section III.3)
$$\left\{\eqalign{
&{1\over \delta^{\kappa-2}}(\GR_\delta-\GI_3) \rightharpoonup \GA \quad \hbox{weakly in}\quad \big(H^1(0,L)\big)^{3\times 3}\cr
&{1\over \delta^{\kappa-2}}{\cal U}_\delta\longrightarrow   {\cal U}  \quad \hbox{strongly in}\quad \big(H^1(0,L)\big)^3\cr
&{1\over \delta^{\kappa-1}}{\cal V}_{S,\delta}\rightharpoonup {\cal V}_{S} \quad \hbox{weakly in}\quad \big(H^1(0,L)\big)^3\cr
&{1\over \delta^\kappa}\Pi _\delta\overline{v}_\delta \rightharpoonup  \overline{v}  \quad \hbox{weakly in}\quad \big(L^2(0,L;H^1(\omega))\big)^3\cr}\right.\leqno(IV.3.4)$$ where  $\GA$ is an antisymetric matrix and ${\cal U} (0)={\cal V}_S(0)=0$. Moreover, ${\cal U}$ belongs to $(H^2(0,L))^3$ and there  exists ${\cal R}\in (H^1(0,L))^3$ with ${\cal R} (0)=0$ such that  
$${d{\cal U} \over ds_3}=\GA\Gt={\cal R}\land \Gt .\leqno (IV.3.5)$$ 
\noindent Furthermore, we also have
$$\left\{\eqalign{
&{1\over \delta^{\kappa-2}}\Pi _\delta u_\delta \longrightarrow    {\cal U} \quad \hbox{strongly in}\quad \big(H^1(\Omega)\big)^3,\cr
&\Pi _\delta (\nabla_x v_\delta)  \longrightarrow  \GI_3 \quad \hbox{strongly in}\quad \big(L^2(\Omega)\big)^{3\times 3},\cr
&{\Pi_\delta(u_\delta-{\cal U}_\delta)\over \delta^{\kappa-1}}\longrightarrow S_1{\cal R}\land\Gn_1+S_2{\cal R}\land\Gn_2\quad \hbox{strongly in}\quad \big(L^2(\Omega)\big)^3,\cr
&{1\over 2\delta^{\kappa-1}}\Pi _\delta \big((\nabla_xv_\delta)^T\nabla_x v_\delta-\GI_3\big) \rightharpoonup   \GE  \qquad\hbox{weakly in}\quad (L^2(\Omega))^{3\times 3},\cr}\right.\leqno(IV.3.6)$$ where (see Section III.3)
$$\left\{\eqalign{
\GE&=(\Gn_1\;|\; \Gn_2\; |\Gt)\widehat{\GE}(\Gn_1\;|\; \Gn_2\; |\Gt)^T\cr\cr
\widehat{\GE}&=\pmatrix{
\ds{\partial\overline{v}\over \partial S_1}\cdot\Gn_1 & \ds{1\over 2}\Big\{{\partial\overline{v}\over \partial S_1}\cdot\Gn_2+{\partial\overline{v}\over \partial S_2}\cdot\Gn_1\Big\} & \ds {1\over 2}\Big\{{\partial\overline{v}\over \partial S_1}\cdot\Gt-S_2{d{\cal R}\over ds_3}\cdot \Gt+{d{\cal V}_S\over ds_3}\cdot \Gn_1 \Big\}\cr
* & \ds{\partial\overline{v}\over \partial S_2}\cdot\Gn_2   &  \ds {1\over 2}\Big\{{\partial\overline{v}\over \partial S_2}\cdot\Gt+S_1{d{\cal R}\over ds_3} \cdot \Gt+{d{\cal V}_S\over ds_3}\cdot \Gn_2\Big\}\cr
 * & *  & \ds  -S_1{d{\cal R}\over ds_3} \cdot\Gn_2+S_2{d{\cal R}\over ds_3} \cdot\Gn_1+{d{\cal V}_S\over ds_3}\cdot \Gt}.\cr}\right.\leqno(IV.3.7)$$ Proceeding as in the previous section, we pass to the limit-inf in $\ds {J(v_\delta)-J(I_d)\over \delta^{2\kappa}}$ and  we obtain
$$\left\{\eqalign{
&\int_\Omega \Big\{{\lambda\over 2}(tr(\GE ))^2+\mu |||\GE |||^2\Big\}\cr
&-\int_0^L \Big(|\omega| f+\sum_{\alpha=1}^2  \int_\omega g (S_1,S_2,.)S_\alpha \det\big(\Gn_1\; |\; \Gn_2\;|\;{d\Gn_\alpha\over ds_3}\big)dS_1dS_2\Big)\cdot {\cal U}\cr
&-\sum_{\alpha=1}^2\int_0^L\Big(\int_\omega g (S_1,S_2,.)S_\alpha dS_1dS_2\Big)\cdot \big({\cal R}\land\Gn_\alpha\big)  \le \liminf_{\delta\to 0}{1\over \delta^{2\kappa}}\big(J(v_\delta)-J(I_d)\big).\cr}\right.\leqno(IV.3.8)$$

Let $\GU_{lin}$ be the space
$$\eqalign{
\GU_{lin}=\Bigl\{ &\big({\cal U}^{'},{\cal R}^{'},{\cal V}^{'}_S,\overline{v}^{'}\big)\in (H^2(0,L))^3\times (H^1(0,L))^{3}\times (H^1(0,L))^3\times (L^2(0,L;H^1(\omega)))^3\;\; |\;\; \cr
&{\cal U}^{'}(0)={\cal R}^{'}(0)={\cal V}^{'}_S(0)=0,\quad{d{\cal U}^{'}\over ds_3}={\cal R}^{'}\land\Gt,\qquad\int_\omega \overline{v}^{'}(S_1,S_2,s_3)dS_1dS_2=0\qquad \hbox{for a.e.}\;\; s_3\in (0,L)\Big\}.\cr}$$ For any $\big({\cal U}^{'},{\cal R}^{'},{\cal V}^{'}_S,\overline{v}^{'}\big)\in \GU_{lin}$, we set
$$\left\{\eqalign{
&{\cal J}_{L}\big({\cal U}^{'},{\cal R}^{'},{\cal V}^{'}_S,\overline{v}^{'}\big)=\int_\Omega \Big\{{\lambda\over 2}(tr(\GE^{'}))^2+\mu |||\GE^{'}|||^2\Big\}\cr
&-\int^L_0\Big(|\omega| f +\sum_{\alpha=1}^2  \int_\omega g (S_1,S_2,.)S_\alpha \det\big(\Gn_1\; |\; \Gn_2\;|\;{d\Gn_\alpha\over ds_3}\big)dS_1dS_2\Big)\cdot {\cal U}^{'}\cr
&-\sum_{\alpha=1}^2\int_0^L \Big(\int_\omega g (S_1,S_2,.)S_\alpha dS_1dS_2\Big)\cdot \big({\cal R}^{'}\land\Gn_\alpha\big),\cr}
\right.\leqno(IV.3.9)$$ with $\GE^{'}$ is given by (IV.3.7) where we have replaced $\big({\cal U},{\cal R},{\cal V}_S,\overline{v}\big)$ by $\big({\cal U}^{'},{\cal R}^{'},{\cal V}^{'}_S,\overline{v}^{'}\big)$.  From (IV.3.8)  it results that
$${\cal J}_{L}\big({\cal U} ,{\cal R} ,{\cal V}_S,\overline{v}\big)\le\liminf_{\delta\to 0}{J(v_\delta)-J(I_d)\over \delta^{2\kappa}}.\leqno(IV.3.10)$$
\medskip
 Now, let $\big({\cal U} ,{\cal R} , {\cal V}_S,\overline{v}\big)$ be in $\GU_{lin}$ and let $\Big({\cal R}_N , {\cal V}_{S,N} ,\overline{v}_N \Big)_{N\in \N^*}$ sequences of elements  such that
$$\left\{\eqalign{
{\cal R}_N(0)=0,\qquad {\cal V}_{S,N}(0)=0,&\qquad \overline{v}_N(S_1,S_2,0)=0,\; \hbox{ for a.e. }\; (S_1,S_2)\in \omega,\cr
{\cal R}_N\in (W^{1,\infty}(0,L))^{3},\qquad {\cal R}_N&\longrightarrow {\cal R}\quad \hbox{strongly in }\quad (H^1(0,L))^{3}\cr
{\cal V}_{S,N}\in (W^{1,\infty}(0,L))^{3},\qquad {\cal V}_{S,N}&\longrightarrow {\cal V}_S\quad \hbox{strongly in }\quad (H^1(0,L))^{3}\cr
\overline{v}_N\in (W^{1,\infty}(\Omega))^{3},\enskip \overline{v}_N&\longrightarrow \overline{v}\quad \hbox{strongly in }\quad (L^2(0,L;H^1(\omega))^{3}\cr}\right.\leqno(IV.3.11 )$$  Moreover we set
$${d{\cal U}_N\over ds_3}={\cal R}_N\land \Gt,\qquad {\cal U}_N(0)=0).$$

We consider the deformations ($\delta\in (0,\delta_0]$ and $ s\in \Omega_\delta$)
$$v_{N,\delta}(s)= {\cal V}_{N,\delta}(s_3)+\GR_{N,\delta}(s_3)(s_1\Gn_1+s_2\Gn_2)+\delta^{\kappa-1}{\cal V}_{S,N}(s_3)+\delta^\kappa\overline{v}_N\Big({s_1\over \delta},{s_2\over \delta},s_3\Big)\leqno(IV.3.12 )$$ where $\GR_{N,\delta}$ and  ${\cal V}_{N,\delta}$ are defined below
$$\left\{\eqalign{
&{d \GR_{N,\delta}\over ds_3}=\delta^{\kappa-2}\GR_{N,\delta}\GB_N\cr
&\GR_{N,\delta}(0)=\GI_3\cr}\right.,\qquad {\cal V}_{N,\delta}(s_3)=M(0)+\int_0^{s_3}\GR_{N,\delta}(z)\Gt(z)dz.$$ Here $\GB_N$ is the $3\times 3$ antisymmetric matrix such that 
$$\forall x\in \R^3,\qquad \GB_Nx={d{\cal R}_N\over ds_3}\land x.$$ Using the above convergences and the fact that $\big({\cal U}_N,{\cal R}_N,{\cal V}_{S,N},\overline{v}_N\big)$ belongs to $\GU_{lin}$, we have
$$\left\{\eqalign{
&{1\over \delta^{\kappa-2}}\Pi _\delta u_{N,\delta} \longrightarrow    {\cal U}_N \quad \hbox{strongly in}\quad \big(W^{1,\infty}(\Omega)\big)^3,\cr
&\Pi _\delta (\nabla_x v_{N,\delta})  \longrightarrow  \GI_3\quad \hbox{strongly in}\quad \big(L^\infty(\Omega)\big)^{3\times 3}.\cr
&{1\over 2\delta^{\kappa-1}}\Pi _\delta \big((\nabla_xv_{N,\delta})^T\nabla_x v_{N,\delta}-\GI_3\big) \longrightarrow \GE_N \qquad\hbox{strongly in}\quad (L^\infty(\Omega))^{3\times 3},\cr}\right.\leqno(IV.3.13 )$$ where  $\GE_N$ is given by (IV.3.7) where we have replaced $\big({\cal U},{\cal R},{\cal V}_S,\overline{v}\big)$ by $\big({\cal U}_N,{\cal R}_N,{\cal V}_{S,N},\overline{v}_N\big)$.  
 If $\delta$ is sufficiently  small  we have $\det\big(\nabla_x v_{N,\delta}(x)\big)>0$ for a.e. $x\in  {\cal P}_\delta$.  We divide $J(v_{N,\delta})-J(I_d)$ by $\delta^{2\kappa}$ and we pass to the limit. We obtain
$$\lim_{\delta\to 0}{1\over \delta^{2\kappa}}\big(J(v_{N,\delta})-J(I_d)\big)={\cal J}_{L}\big({\cal  U}_N,{\cal R}_N, {\cal V}_{S,N},\overline{v}_N\big).\leqno(IV.3.14 )$$ Now  letting $N$ tend to $+\infty$ gives that for any $\big({\cal U} ,{\cal R} , {\cal V}_{S} ,\overline{v} \big)\in \GU_{lin}$
$${\cal J}_{L}\big({\cal U} ,{\cal R} , {\cal V}_{S} ,\overline{v} \big)=\lim_{N\to+\infty}{\cal J}_{L}\big({\cal  U}_N,{\cal R}_N, {\cal V}_{S,N},\overline{v}_N\big).\leqno(IV.3.15 )$$ Hence, through a standard diagonal process for any $\big({\cal U} ,{\cal R} , {\cal V}_{S} ,\overline{v} \big)\in \GU_{lin}$ there exists a sequence of admissible deformations $v_\delta\in(H^1({\cal P}_\delta))^3$ such that
$${\cal J}_{L}\big({\cal U} ,{\cal R} , {\cal V}_{S} ,\overline{v} \big)=\lim_{\delta\to 0}{J(v_\delta)-J(I_d)\over \delta^{2\kappa}}.\leqno(IV.3.16)$$ The following theorem summarizes the results of the case $\kappa>2$.
\smallskip
\noindent{\ggras  Theorem IV.3.1. }{\it  The functional ${\cal J}_{L}$ is the $\Gamma$-limit of $\ds {J(.)-J(I_d)\over \delta^{2\kappa}}$ in the following sense: 

\noindent $\bullet$ for any sequence of deformations  $\big(v_\delta\big)_{0<\delta\le \delta_0 }$   belonging to $\GU_\delta$ and satisfying $$\ds \liminf_{\delta\to 0}{J(v_\delta)-J(I_d)\over \delta^{2\kappa}}<+\infty$$ and let $\big({\cal U}_\delta,\GR_\delta,{\cal V}_{S,\delta},\overline{v}_\delta\big)$ be the  terms of the decomposition of the displacement  $u_\delta=v_\delta-I_d$ given by (III.3.2).   Up to a subsequence there exists $\big({\cal U} , {\cal R} , {\cal V}_{S} ,\overline{v} \big)\in \GU_{lin}$ such that
$$\eqalign{
&{1\over \delta^{\kappa-2}}(\GR_\delta-\GI_3) \rightharpoonup \GA \quad \hbox{weakly in}\quad \big(H^1(0,L)\big)^{3\times 3}\cr
&{1\over \delta^{\kappa-2}}{\cal U}_\delta\longrightarrow   {\cal U}  \quad \hbox{strongly in}\quad \big(H^1(0,L)\big)^3\cr
&{1\over \delta^{\kappa-1}}{\cal V}_{S,\delta}\rightharpoonup {\cal V}_{S} \quad \hbox{weakly in}\quad \big(H^1(0,L)\big)^3\cr
&{1\over \delta^\kappa}\Pi _\delta\overline{v}_\delta \rightharpoonup  \overline{v}  \quad \hbox{weakly in}\quad \big(L^2(0,L;H^1(\omega))\big)^3\cr}$$ where  for any $x\in \R^3$, $\GA x={\cal R}\land x$ and we have
$${\cal J}_{L}\big({\cal U} ,{\cal R} , {\cal V}_{S} ,\overline{v} \big)\le \liminf_{\delta\to 0}{J(v_\delta)-J(I_d)\over \delta^{2\kappa}}$$
$\bullet$ for any $\big({\cal U} ,{\cal R} , {\cal V}_{S} ,\overline{v} \big)\in \GU_{lin}$ there exists a sequence $\big(v_\delta\big)_{0<\delta\le \delta_0 }$   belonging to $\GU_\delta$ such that
$${\cal J}_{L}\big({\cal U} ,{\cal R} , {\cal V}_{S} ,\overline{v} \big)= \lim_{\delta\to 0}{J(v_\delta)-J(I_d)\over \delta^{2\kappa}}.$$ Moreover, there exists  $\big({\cal U}_0 ,{\cal R}_0 , {\cal V}_{S,0} ,\overline{v}_0 
\big)\in \GU_{lin}$ such that
$$m_\kappa=\lim_{\delta\to 0}{m_\delta\over \delta^{2\kappa}}={\cal J}_{L}\big({\cal U}_0,{\cal R}_0, {\cal V}_{S,0} ,\overline{v}_0\big)=\min_{\big({\cal U} ,{\cal R} , {\cal V}_{S} ,\overline{v} \big)\in\GU_{lin}}{\cal J}_{L}\big({\cal U} ,{\cal R} , {\cal V}_{S} ,\overline{v} \big).$$ }
The next theorem is the analog of Theorem IV.2.2.
\smallskip
\noindent{\ggras  Theorem  IV.3.2 }{\it Let $({\cal U}_0 ,{\cal R}_0)$ be given by Theorem IV.3.1. The minimum $m_\kappa$ of the functional ${\cal J}_{L}$ over $\GU_{lin}$ satisfies the following minimization problem which admits a unique solution:
$$m_\kappa={\cal F}_{L}\big({\cal U}_0,{\cal R}_0\big)=\min_{\big({\cal U},{\cal R}\big)\in\GV_{lin}}{\cal F}_{L}\big({\cal U},{\cal R}\big), \leqno(IV.3.17) $$ where 
$$\eqalign{
\GV_{lin}=\Bigl\{ &\big({\cal U},{\cal R}\big)\in (H^2(0,L))^3\times (H^1(0,L))^{3}\;\; |\;\; 
{\cal U}(0)={\cal R}(0)=0,\qquad {d{\cal U} \over ds_3}={\cal R}\land\Gt\Big\},\cr}$$ and 
$$\left\{\eqalign{
& {\cal F}_{L}\big({\cal U},{\cal R}\big)={EI_1\over 2}\int_0^L\Big({d{\cal R}\over ds_3} \cdot  \Gn_2\Big)^2+{EI_2\over 2}\int_0^L\Big({d{\cal R}\over ds_3} \cdot  \Gn_1\Big)^2+{\mu K\over 4}\int_0^L\Big({d{\cal R}\over ds_3} \cdot \Gt\Big)^2\cr
&-\int_0^L\Big(|\omega| f +\sum_{\alpha=1}^2  \int_\omega g (S_1,S_2,.)S_\alpha \det\big(\Gn_1\; |\; \Gn_2\;|\;{d\Gn_\alpha\over ds_3}\big)dS_1dS_2\Big)\cdot {\cal U}\cr
&-\sum_{\alpha=1}^2 \int_0^L\Big(\int_\omega g (S_1,S_2,.)S_\alpha dS_1dS_2\Big)\cdot({\cal R}\land\Gn_\alpha)\cr}\right.\leqno(IV.3.18)$$ $E$ is the Young's modulus  and $K$ is given in Theorem IV.1.3. }
\medskip
\noindent{\ggras Proof of Theorem IV.3.2. } We proceed as in Theorem IV.2.2. We fix  $\big({\cal V} , {\cal R} \big)\in \GV_{lin}$  and we minimize the functional ${\cal J}_{L}\big({\cal U} ,{\cal R} ,\cdot,\cdot\big)$ over the space $\GW$.  Through  solving simple variational problems    (see [14] again), we find that the minimum of the functional ${\cal J}_{L}\big({\cal U} ,{\cal R} ,\cdot,\cdot\big)$ over the space $\GW$ is obtained with  $\ds{d{\cal V} _S\over ds_3}\cdot \Gt=0$ and
$$\left\{\eqalign{
&\overline{v} (S_1,S_2,.)\cdot \Gn_1 = -\nu\Big\{{S_2^2-S_1^2\over 2}{d{\cal R} \over ds_3} \cdot  \Gn_2+S_1S_2{d{\cal R} \over ds_3}\cdot  \Gn_1\Big\}\cr
&\overline{v} (S_1,S_2,.)\cdot \Gn_2= -\nu\Big\{-S_1S_2{d{\cal R}  \over ds_3}\cdot \Gn_2+{S_2^2-S_1^2\over 2}{d{\cal R}  \over ds_3}\cdot   \Gn_1\Big\}\cr
&\overline{v} (S_1,S_2,.)\cdot \Gt+S_1{d{\cal V}_S\over ds_3}\cdot \Gn_1 +S_2{d{\cal V}_S\over ds_3}\cdot \Gn_2 = \chi(S_1,S_2){d{\cal R} \over ds_3} \cdot \Gt \cr}\right.\leqno(IV.3.19)$$  Then the symmetric tensor $\widehat{\GE}$ (see (IV.3.7)) at the minimum is given by
$$\widehat{\GE}({\cal R})=\pmatrix{
\ds-\nu\widehat{\GE}_{33}({\cal R}) & 0 & \ds {1\over 2}\Big({\partial\chi\over \partial S_1}-S_2\Big){d{\cal R} \over ds_3} \cdot \Gt  \cr
* & -\nu\widehat{\GE}_{33}({\cal R}) &  \ds {1\over 2}\Big({\partial\chi\over \partial S_1}+S_1\Big){d{\cal R} \over ds_3} \cdot \Gt \cr
 * & *  & \ds  \widehat{\GE}_{33}({\cal R})\cr},\leqno(IV.3.20)$$ where $\ds \widehat{\GE}_{33}({\cal R})=-S_1{d{\cal R} \over ds_3}\cdot  \Gn_2+S_2{d{\cal R} \over ds_3}\cdot  \Gn_1$. Upon replacing $\GE$ by $\widehat{\GE}({\cal R})$  in the expression of ${\cal J}_{L}$ we obtain 
$$\min_{\big({\cal V} _S,\overline{v} \big)\in\GW}{\cal J}_{L}\big({\cal U} ,{\cal R}  ,{\cal V} _S,\overline{v} \big)={\cal F}_{L}\big({\cal U} , {\cal R}  \big)$$  where the functional ${\cal F}_{L}$ is given by (IV.3.18).\fin
\noindent {\ggras Remark.}  The above analysis shows that if $\big(v_\delta\big)_{0<\delta\le \delta_0}$ is a sequence such that 
$$m_\kappa=\lim_{\delta\to 0}{J(v_\delta)-J(I_d)\over \delta^{2\kappa}},$$
then the sequence of the Green-St Venant's deformation tensors satisfies
 $${1\over 2\delta^{\kappa-1}}\Pi _\delta \big((\nabla_xv_\delta)^T\nabla_x v_\delta-\GI_3\big) \longrightarrow   (\Gn_1\,|\, \Gn_2\,|\, \Gt)\widehat{\GE}({\cal R}_0) (\Gn_1\,|\, \Gn_2\,|\, \Gt)^T  \qquad\hbox{strongly in}\quad (L^2(\Omega))^{3\times 3},$$
 where $\widehat{\GE}({\cal R}_0)$ is defined in (IV.3.20) and ${\cal R}_0$ is the solution of (IV.3.17).
 \medskip

\noindent {\Ggras  IV.4   Extentional models for special forces.}
\medskip
\noindent     In this subsection we investigate the  case where $f_\delta$ is given by
$$f_\delta(s)=\delta^{\kappa-1}f(s_3)\qquad \hbox{for a.e.}\enskip s\in \Omega_\delta,\leqno(IV.4.1)$$ where $f$ belongs to  $(L^2(0,L))^3$. Without any additional assumption on $f$, Subsection IV.1 shows that this leads to 
$$||\hbox{dist}(\nabla v_\delta, SO(3))||_{L^2({\cal P}_\delta)}\le C\delta^{\kappa-1},$$ if $\ds{1\over \delta^{2\kappa-2}}\big(J(v_\delta)-J(I_d)\big)\le C_1$. As a consequence, the results of Subsections IV.2 et IV.3 can be applied if $\kappa \ge 3)$. Let us for example consider the case $\kappa > 3$ and   remark that due to the choice of $f$ the contribution of the forces in the limit energy ${\cal F}_L({\cal U}_0, {\cal R}_0)$ is equal to 
$$-|\omega|\int_0^Lf(s_3)\cdot {\cal U}_0(s_3)ds_3=-|\omega|\int_0^L\Big(\int_{s_3}^L f(s)ds\Big)\cdot {\cal R}_0(s_3)\land\Gt(s_3)ds_3.$$  
Then if the quantity $\ds \int_{s_3}^L f(s)ds$ is proportionnal to $\Gt(s_3)$, this contribution vanishes and then $ {\cal R}_0={\cal U}_0=0$ and the minimum is null. This example shows that for this kind of special forces, the energy have a smaller order than $2\kappa-2$ or equivalently that the estimates on $v_\delta$ can be improved in this case.

\medskip
\noindent We assume  that there exists $\widetilde{f} \in H^1(0,L)$ such that
$$ \int_{s_3}^Lf(l)dl=\widetilde{f}(s_3)\Gt(s_3) \qquad \hbox{for any}\enskip s_3\in [0,L].\leqno(IV.4.2)$$ 

Let $v$ an admissible deformation of the rod ${\cal P}_\delta$. Now, using (IV.4.2) we derive a new estimate of $\ds \int_{{\cal P}_\delta}f_\delta\cdot(v-I_d)$.
\noindent Notice first that $\ds \det(\nabla\Phi)=1+s_1\det\Big(\Gn_1\; |\; \Gn_2\;|\;{d\Gn_1\over ds_3}\Big)+s_2\det\Big(\Gn_1\; |\; \Gn_2\;|\;{d\Gn_2\over ds_3}\Big)$, then using the decomposition (II.2.1) for the admissible deformation $v$, estimates of Theorem II.2.2 and (II.3.5) together with $\ds \int_{\omega_\delta} s_\alpha ds_1ds_2=0$ we deduce that
$$\left\{\eqalign{
&\Big|\int_{{\cal P}_\delta}f_\delta\cdot(v-I_d)-|\omega|\delta^{\kappa+1}\int_0^Lf(s_3)\cdot\big({\cal V}(s_3)-M(s_3)\big)ds_3\Big|\cr
\le& C\delta^{\kappa+1}||f||_{(L^2(0,L))^3}||\hbox{dist}(\nabla v, SO(3))||_{L^2({\cal P}_\delta)}.\cr}\right.\leqno(IV.4.3)$$  We obtain by integrating by parts and using the decomposition (II.2.16) of ${\cal V}$ (see also (III.4.7))
$$\left\{\eqalign{
&\int_0^L f(s_3)\cdot\big({\cal V}(s_3)-M(s_3)\big)ds_3
=\int_0^L\widetilde{f}(s_3)\Gt(s_3)
\cdot\Bigl({d{\cal V}\over ds_3}(s_3)-\Gt(s_3)\Big)ds_3\cr
=&\int_0^L\widetilde{f}(s_3)\Gt(s_3) \cdot (\GR(s_3)-\GI_3)\Gt(s_3)ds_3+\int_0^L\widetilde{f}(s_3)\Gt(s_3)
\cdot {d{\cal V}_{S}\over ds_3}(s_3)ds_3\cr
=&-{1\over 2}\int_0^L\widetilde{f}(s_3)(\GR(s_3)-\GI_3)\Gt(s_3) \cdot (\GR(s_3)-\GI_3)\Gt(s_3)ds_3+\int_0^L\widetilde{f}(s_3){d{\cal V}_{S}\over ds_3}(s_3)\cdot\Gt(s_3)ds_3.\cr}\right.\leqno(IV.4.4)$$ Finally, using the above estimate and (II.2.17) we get
$$\left\{\eqalign{
&\Big|\int_{{\cal P}_\delta}f_\delta\cdot(v-I_d)+{|\omega|\delta^{\kappa+1}\over 2}\int_0^L\widetilde{f}(s_3)(\GR(s_3)-\GI_3)\Gt(s_3) \cdot (\GR(s_3)-\GI_3)\Gt(s_3)ds_3\Big|\cr
\le & C\delta^{\kappa}||f||_{(L^2(0,L))^3}||\hbox{dist}(\nabla v, SO(3))||_{L^2({\cal P}_\delta)}.\cr}\right.\leqno(IV.4.5)$$ We assume that
$$ \int_{{\cal P}_\delta}W(\nabla v)-\int_{{\cal P}_\delta}f_\delta\cdot(v-I_d)= J(v)-J(I_d)<+\infty$$ which implies using (IV.1.9)
$${\mu\over 4}||\hbox{dist}(\nabla v,SO(3))||^2_{L^2({\cal P}_\delta)}-\int_{{\cal P}_\delta}f_\delta\cdot(v-I_d)\le  J(v)-J(I_d)<+\infty.\leqno (IV.4.6)$$ Hence
$$\left\{\eqalign{
&{\mu\over 4}||\hbox{dist}(\nabla v,SO(3))||^2_{L^2({\cal P}_\delta)}+{|\omega|\delta^{\kappa+1}\over 2}\int_0^L\widetilde{f}(\GR-\GI_3)\Gt  \cdot (\GR -\GI_3)\Gt\cr
&\le C\delta^{\kappa}||f||_{(L^2(0,L))^3}||\hbox{dist}(\nabla v, SO(3))||_{L^2({\cal P}_\delta)}+ J(v)-J(I_d)<+\infty \cr}\right.\leqno(IV.4.7)$$ Now, in view of the above inequality, let us consider a sequence  $\big(v_\delta\big)_{0<\delta\le\delta_0}$  satisfiyng $J(v_\delta)-J(I_d)\le C_1\delta^{2\kappa}$.

\noindent  From estimate (II.3.8) we deduce that $$\Big|\int_0^L\widetilde{f}(\GR-\GI_3)\Gt  \cdot (\GR-\GI_3)\Gt\Big|\le C^*\delta^{-4}||f||_{(L^2(0,L))^3}||\hbox{dist}(\nabla v, SO(3))||^2_{L^2({\cal P}_\delta)},\leqno(IV.4.8)$$
where the constant $C^*$ only  depends upon the geometry of the middle line of the rod.  According to  inequalities (IV.4.8) and (IV.4.7) the cases $\kappa=3$ and $\kappa>3$ lead two different energy estimates.

\noindent  If $\kappa>3$ we obtain
$$\left\{\eqalign{
&||\hbox{dist}(\nabla v_\delta, SO(3))||_{L^2({\cal P}_\delta)}\le C\delta^{\kappa},\cr
& \big\|{1\over 2}\big\{\nabla v_\delta^T\nabla v_\delta-\GI_3\big\}\big\|_{(L^2({\cal P}_\delta))^{3\times 3}}\le  C\delta^{\kappa},\cr
 &||\nabla v_\delta||_{(L^4({\cal P}_\delta))^{3\times 3}}\le  C\delta^{1\over 2}\cr}\right.\leqno(IV.4.9)$$ The constant  $C$ does not  depend   on $\delta$.
\medskip
\noindent  If $\kappa=3$, the energy estimate depends on $||f||_{(L^2(0,L))^3}$. Indeed , if 
$$\ds ||f||_{(L^2(0,L))^3}< {\mu \over 2C^*|\omega|}\leqno(IV.4.10)$$ estimate (IV.4.8) and (IV.4.7) give
$$\left\{\eqalign{
&\Big({\mu\over 4}-{C^*\over 2}|\omega|\,||f||_{(L^2(0,L))^3}\Big)||\hbox{dist}(\nabla v_\delta,SO(3))||^2_{L^2({\cal P}_\delta)}\cr
&\le C\delta^3||f||_{(L^2(0,L))^3}||\hbox{dist}(\nabla v, SO(3))||_{L^2({\cal P}_\delta)}+ C_1\delta^6 \cr}\right.\leqno(IV.4.11)$$ and then
$$\left\{\eqalign{
&||\hbox{dist}(\nabla v_\delta, SO(3))||_{L^2({\cal P}_\delta)}\le C\delta^{3},\cr
& \big\|{1\over 2}\big\{\nabla v_\delta^T\nabla v_\delta-\GI_3\big\}\big\|_{(L^2({\cal P}_\delta))^{3\times 3}}\le  C\delta^{3},\cr
 &||\nabla v_\delta||_{(L^4({\cal P}_\delta))^{3\times 3}}\le  C\delta^{1\over 2}\cr}\right.\leqno(IV.4.12)$$ The constant  $C$ does not depend on $\delta$.
 \medskip

\noindent In view of (IV.4.7), an alternative assumption to (IV.4.10) which also leads to estimates (IV.4.12) is to suppose that   
$$\widetilde{f}(s_3)\ge 0\quad\hbox{ for almost any }\; \; s_3\in (0,L).\leqno(IV.4.13)$$
\smallskip
 {\it In the sequel of this subsection, we will assume that the forces $f$ satisfy (IV.4.2) together with (IV.4.10) or (IV.4.13) if $\kappa=3$.}
\smallskip

Now we have to pass to the limit-inf in $\ds {J(v_\delta)-J(I_d)\over \delta^{2\kappa}}$. According to estimates (IV.4.9) and (IV.4.12), performing this process in the elastic energy term is identical to the one detailed in the previous Subsection. We just focus on the behavior of the terms involving the forces. In view of (IV.4.3), (IV.4.9) and (IV.4.12) we get
$$\lim_{\delta\to 0}{1\over \delta^{2\kappa} }\int_{{\cal P}_\delta}f_\delta\cdot(v_\delta-I_d)= \lim_{\delta\to 0}{|\omega|\over \delta^{\kappa-1}}\int_0^Lf(s_3)\cdot\big({\cal V}_\delta(s_3)-M(s_3)\big)ds_3.$$ Now we use the notations and results of Subsection III.4, we have
$$\int_0^Lf \cdot\big({\cal V}_\delta -M\big)=\int_0^Lf\,\cdot\,{\cal U}_\delta =\int_0^L\widetilde{f}\, \Gt\,\cdot{d{\cal U}_\delta\over ds_3} =\int_0^L\widetilde{f}\,  \Gt\cdot\,{dU_{E,\delta}\over ds_3}.$$ Thanks to the convergences of Lemma III.4.1 we deduce that
$$\lim_{\delta\to 0}{1\over \delta^{2\kappa} }\int_{{\cal P}_\delta}f_\delta\cdot(v_\delta-I_d)=|\omega|\int_0^L\widetilde{f}\,  \Gt\cdot\,{dU_{E}\over ds_3}.\leqno (IV.4.14)$$
\noindent Let us define the limit operator ${\cal J}_{LS}$ by
$$\forall({\cal U},{\cal R}, {\cal V}_S,\overline{v})\in\GU_{lin},\qquad
{\cal J}_{LS}({\cal U},{\cal R}, {\cal V}_S,\overline{v})=\int_\Omega \Big\{{\lambda\over 2}(tr(\widehat{\GE }))^2+\mu |||\widehat{\GE} |||^2\Big\}-|\omega|\int_0^L \widetilde{f}\,{dU_E\over ds_3}\cdot\Gt.\leqno(IV.4.15)$$
The matrix $\widehat{\GE}$ is given by (IV.3.7) and the displacement $U_E$ is such that (see Lemma III.4.1)
$${dU_E\over ds_3}\cdot\Gt=
\left\{\eqalign{
&{d{\cal V}_S\over ds_3}\cdot\Gt-{1\over 2}\Big\|{d{\cal U}\over ds_3}\Big\|^2_2\qquad\hbox{if}\; \kappa=3,\cr
&{d{\cal V}_S\over ds_3}\cdot\Gt\qquad\hbox{if}\; \kappa>3.\cr}\right.\leqno(IV.4.16)$$
The expression of ${\cal J}_{LS}$ shows that this functional has a unique minimizer.
We have  obtained the following result.

\noindent{\ggras  Theorem IV.4.1. }{\it  The functional ${\cal J}_{LS}$ is the $\Gamma$-limit of $\ds {J(.)-J(I_d)\over \delta^{2\kappa}}$ in the following sense: 

\noindent $\bullet$ for any sequence of deformations  $\big(v_\delta\big)_{0<\delta\le \delta_0 }$   belonging to $\GU_\delta$ and satisfying $$\ds \liminf_{\delta\to 0}{J(v_\delta)-J(I_d)\over \delta^{2\kappa}}<+\infty$$ and let $\big({\cal U}_\delta,\GR_\delta,{\cal V}_{S,\delta},\overline{v}_\delta\big)$ be the terms of the decomposition of the displacement  $u_\delta=v_\delta-I_d$ given by (III.3.2) and (III.4.3).   Up to a subsequence there exists $\big({\cal U} , {\cal R} , {\cal V}_{S} ,\overline{v} \big)\in \GU_{lin}$ such that
$$\eqalign{
&{1\over \delta^{\kappa-2}}(\GR_\delta-\GI_3) \rightharpoonup \GA \quad \hbox{weakly in}\quad \big(H^1(0,L)\big)^{3\times 3}\cr
&{1\over \delta^{\kappa-2}}{\cal U}_\delta\longrightarrow   {\cal U}  \quad \hbox{strongly in}\quad \big(H^1(0,L)\big)^3\cr
&{1\over \delta^{\kappa-1}}{\cal V}_{S,\delta}\rightharpoonup {\cal V}_{S} \quad \hbox{weakly in}\quad \big(H^1(0,L)\big)^3\cr
&{1\over \delta^{\kappa-1}}U_{E,\delta}\rightharpoonup U_E \quad \hbox{weakly in}\quad \big(H^1(0,L)\big)^3\cr
&{1\over \delta^{\kappa}}\Pi _\delta\overline{v}_\delta \rightharpoonup  \overline{v}  \quad \hbox{weakly in}\quad \big(L^2(0,L;H^1(\omega))\big)^3\cr}$$ where  for any $x\in \R^3$, $\GA x={\cal R}\land x$ and  where the relation between $U_E$, ${\cal V}_S$ and ${\cal U}$ is given by (IV.4.16). We have
$${\cal J}_{LS}\big({\cal U} ,{\cal R} , {\cal V}_{S} ,\overline{v} \big)\le \liminf_{\delta\to 0}{J(v_\delta)-J(I_d)\over \delta^{2\kappa}}$$
$\bullet$ for any $\big({\cal U} ,{\cal R} , {\cal V}_{S} ,\overline{v} \big)\in \GU_{lin}$ there exists a sequence $\big(v_\delta\big)_{0<\delta\le \delta_0 }$   belonging to $\GU_\delta$such that
$${\cal J}_{LS}\big({\cal U} ,{\cal R} , {\cal V}_{S} ,\overline{v} \big)= \lim_{\delta\to 0}{J(v_\delta)-J(I_d)\over \delta^{2\kappa}}.$$ Moreover, there exists a unique $\big({\cal U}_0 ,{\cal R}_0 , {\cal V}_{S,0} ,\overline{v}_0 
\big)\in \GU_{lin}$ such that
$$m_\kappa=\lim_{\delta\to 0}{m_\delta\over \delta^{2\kappa}}={\cal J}_{LS}\big({\cal U}_0,{\cal R}_0, {\cal V}_{S,0} ,\overline{v}_0\big)=\inf_{\big({\cal U} ,{\cal R} , {\cal V}_{S} ,\overline{v} \big)\in\GU_{lin}}{\cal J}_{LS}\big({\cal U} ,{\cal R} , {\cal V}_{S} ,\overline{v} \big).$$ }
The next theorem is the analog of Theorems IV.2.2 and  IV.3.2.
\smallskip
\noindent{\ggras  Theorem  IV.4.2 }{\it Let $({\cal U}_0,{\cal R}_0,{\cal V}_{S,0})$ be given by Theorem IV.4.1 and $U_{E,0}\in D_{Ex}$ defined by (IV.4.16). The minimum $m_\kappa$ of the functional ${\cal J}_{LS}$ over $\GU_{lin}$ is obtained with ${\cal U}_0={\cal R}_0=0$ and  it is given by the following minimization problem which admits a unique solution:
$$m_\kappa={\cal F}_{LS}\big(U_{E,0}\big)=\min_{U_E\in D_{Ex}}{\cal F}_{LS}\big(U_E\big), \leqno(IV.4.17) $$ where 
$${\cal F}_{LS}\big(U_E\big)=|\omega|\Big\{{E\over 2}\int_0^L\Big({dU_E\over ds_3} \cdot  \Gt\Big)^2-\int_0^L \widetilde{f}\,{dU_E\over ds_3}\cdot\Gt\Big\}.\leqno(IV.4.18)$$ $E$ is the Young's modulus.}
\medskip
\noindent{\ggras Proof of Theorem IV.4.2. } We proceed as in Theorem IV.3.2. We fix  $\big({\cal V} ,{\cal R},{\cal V}_S \big)$  and we minimize the functional $\overline{v}\longmapsto{\cal J}_{LS}\big({\cal U} ,{\cal R} ,{\cal V}_S,\overline{v}\big)$ over the space
$$\widetilde{\GW}=\Bigl\{ \overline{v}^{'}\in (L^2(0,L;H^1(\omega)))^3\;\; |\;\;  \int_\omega \overline{v}^{'}(S_1,S_2,s_3)dS_1dS_2=0\qquad \hbox{for a.e.}\;\; s_3\in (0,L)\Big\}.$$   Through  solving simple variational problems    (see [14] again), we find that the minimum of this functional  is obtained for
$$\left\{\eqalign{
&\overline{v} (S_1,S_2,.)\cdot \Gn_1 = -\nu\Big\{S_1{d{\cal V} _S\over ds_3}\cdot \Gt+{S_2^2-S_1^2\over 2}{d{\cal R} \over ds_3} \cdot  \Gn_2+S_1S_2{d{\cal R} \over ds_3}\cdot  \Gn_1\Big\}\cr
&\overline{v} (S_1,S_2,.)\cdot \Gn_2= -\nu\Big\{S_2{d{\cal V} _S\over ds_3}\cdot \Gt-S_1S_2{d{\cal R}  \over ds_3}\cdot \Gn_2+{S_2^2-S_1^2\over 2}{d{\cal R}  \over ds_3}\cdot   \Gn_1\Big\}\cr
&\overline{v} (S_1,S_2,.)\cdot \Gt+S_1{d{\cal V}_S\over ds_3}\cdot \Gn_1 +S_2{d{\cal V}_S\over ds_3}\cdot \Gn_2 = \chi(S_1,S_2){d{\cal R} \over ds_3} \cdot \Gt \cr}\right.\leqno(IV.4.19)$$  Then the symmetric tensor $\widehat{\GE}$ (see (IV.3.7)) at the minimum is given by
$$\widehat{\GE}({\cal R},{\cal V} _S)=\pmatrix{
\ds-\nu\widehat{\GE}_{33}({\cal R},{\cal V} _S) & 0 & \ds {1\over 2}\Big({\partial\chi\over \partial S_1}-S_2\Big){d{\cal R} \over ds_3} \cdot \Gt  \cr
* & -\nu\widehat{\GE}_{33}({\cal R},{\cal V} _S) &  \ds {1\over 2}\Big({\partial\chi\over \partial S_1}+S_1\Big){d{\cal R} \over ds_3} \cdot \Gt \cr
 * & *  & \ds  \widehat{\GE}_{33}({\cal R},{\cal V} _S)\cr},\leqno(IV.4.20)$$ where $\ds \widehat{\GE}_{33}({\cal R},{\cal V} _S)={d{\cal V} _S\over ds_3}\cdot \Gt-S_1{d{\cal R} \over ds_3}\cdot  \Gn_2+S_2{d{\cal R} \over ds_3}\cdot  \Gn_1$. Upon replacing $\GE$ by $\widehat{\GE}({\cal R},{\cal V} _S)$  in the expression of ${\cal J}_{LS}$ and using (IV.4.16), we obtain that the minimum of the functional  ${\cal J}_{LS}\big({\cal U} ,{\cal R} ,{\cal V}_S,\cdot\big)$ over the space $\widetilde{\GW}$ is equal to:
\smallskip 
\noindent  $\bullet$ if $\kappa>3$ 
 $$\eqalign{
& {EI_1\over 2}\int_0^L\Big({d{\cal R}\over ds_3} \cdot  \Gn_2\Big)^2+{EI_2\over 2}\int_0^L\Big({d{\cal R}\over ds_3} \cdot  \Gn_1\Big)^2+{\mu K\over 4}\int_0^L\Big({d{\cal R}\over ds_3} \cdot \Gt\Big)^2
+ {E|\omega|\over 2}\int_0^L\Big({dU_E\over ds_3} \cdot  \Gt\Big)^2-|\omega|\int_0^L \widetilde{f}\,{dU_E\over ds_3}\cdot\Gt\cr}$$ then, we immediately deduce that the minimum $m_\kappa$ of the above quantity is obtained with ${\cal U}_0={\cal R}_0=0$ and it is given by the minimum of the functional ${\cal F}_{LS}$ defined by (IV.4.18).
\smallskip
\noindent  $\bullet$ if $\kappa=3$
$$\eqalign{
& {EI_1\over 2}\int_0^L\Big({d{\cal R}\over ds_3} \cdot  \Gn_2\Big)^2+{EI_2\over 2}\int_0^L\Big({d{\cal R}\over ds_3} \cdot  \Gn_1\Big)^2+{\mu K\over 4}\int_0^L\Big({d{\cal R}\over ds_3} \cdot \Gt\Big)^2+
{|\omega|\over 2}\int_0^L \widetilde{f}||{\cal R}\land\Gt||_2^2\cr
+& {E|\omega|\over 2}\int_0^L\Big({d{\cal V}_S\over ds_3} \cdot  \Gt\Big)^2-|\omega|\int_0^L \widetilde{f}\,{d{\cal V}_S\over ds_3}\cdot\Gt,\cr}$$ where the relation between ${\cal V}_S$, ${\cal U}$ and $U_E$ is given by (IV.4.16). Now, if $\widetilde{f}(s_3)\ge 0$ for a.e. $s_3\in (0,L)$, then the minimum $m_3$ of the above quantity is obtained with ${\cal U}_0={\cal R}_0=0$. 

We now prove that under the condition (IV.4.10) we still have  ${\cal U}_0={\cal R}_0=0$. To this let $\big({\cal U} ,{\cal R} ,\,0\,,\overline{v}\big)$ (we have chosen ${\cal V}_S=0$) be in $\GU_{lin}$ and $v_\delta$ be a sequence of admissible deformations given by Theorem (IV.3.1) such that

$${\cal J}_{LS}\big({\cal U} ,{\cal R} ,\, 0\,,\overline{v}\big)=\lim_{\delta\to 0}{J(v_\delta)-J(I_d)\over \delta^{2\kappa}},$$
and
$$||\GE||_{(L^2(\Omega))^{3\times 3}}=||\widehat{\GE}||_{(L^2(\Omega))^{3\times 3}}=
\lim_{\delta\to 0}{1\over \delta^{\kappa}} ||\hbox{dist}(\nabla_x v_\delta,SO(3))||_{L^2({\cal P}_\delta)} . $$
 In view of (IV.4.3), (IV.4.4), (IV.4.6) and  (IV.4.8)   we obtain
$$\Big({\mu\over 4}-{C^*\over 2}|\omega|||f||_{(L^2(0,L))^3}\Big) ||\GE||^2_{(L^2(\Omega))^{3\times 3}} \le {\cal J}_{LS}\big({\cal U} ,{\cal R} ,\,0\,,\overline{v}\big).$$ Now we choose  ${\overline v}$ as the minimizer of ${\cal J}_{LS}\big({\cal U} ,{\cal R} ,\,0\,,\cdot \big)$ over $\widetilde {\GW}$ in the above inequality, it gives
$$\eqalign{
&\Big({\mu\over 4}-{C^*\over 2}|\omega|||f||_{(L^2(0,L))^3}\Big)\Big[\int_0^L(1+2\nu^2)\Big\{I_1\Big({d{\cal R}\over ds_3} \cdot  \Gn_2\Big)^2+ I_2\Big({d{\cal R}\over ds_3} \cdot  \Gn_1\Big)^2\Big\}+{K\over 2}\int_0^L\Big({d{\cal R}\over ds_3} \cdot \Gt\Big)^2\Big]\cr
 \le & {EI_1\over 2}\int_0^L\Big({d{\cal R}\over ds_3} \cdot  \Gn_2\Big)^2+{EI_2\over 2}\int_0^L\Big({d{\cal R}\over ds_3} \cdot  \Gn_1\Big)^2+{\mu K\over 4}\int_0^L\Big({d{\cal R}\over ds_3} \cdot \Gt\Big)^2+
{|\omega|\over 2}\int_0^L \widetilde{f}||{\cal R}\land\Gt||_2^2.\cr}$$
It follows from the above analysis that the minimum  $m_3$ is obtained for ${\cal U}_0={\cal R}_0=0$. In both cases the minimum $m_3$ is given by the minimum of the functional ${\cal F}_{LS}$ defined by (IV.4.18).)\fin
\smallskip
\noindent {\bf Remark}. Still in the case where $\kappa=3$, if $f$ satisfies (IV.4.2) without  any other assumption, the limit deformation is given by Theorem IV.2.2. 

\medskip
We assume $\kappa=3$. Let us end this section with an interesting model obtained by superposing  general forces given by (IV.1.5) denoted here by $(f^0_\delta,g^0_\delta)$  and  the special forces $f^1_\delta$ described in (IV.4.2) satisfying the assumptions of this section ((IV.4.10) or (IV.4.13)). The corresponding total limit energy in this case is then
$$\left\{\eqalign{
&\forall({\cal U},{\cal R}, {\cal V}_S,\overline{v})\in\GU_{lin},\qquad
{\cal J}_{LG}({\cal U},{\cal R}, {\cal V}_S,\overline{v})=\int_\Omega \Big\{{\lambda\over 2}(tr(\widehat{\GE }))^2+\mu |||\widehat{\GE} |||^2\Big\}\cr
&-\int_0^L\Big(|\omega| f^0 +\sum_{\alpha=1}^2  \int_\omega g^0(S_1,S_2,.)S_\alpha \det\big(\Gn_1\; |\; \Gn_2\;|\;{d\Gn_\alpha\over ds_3}\big)dS_1dS_2\Big)\cdot {\cal U}\cr
&-\sum_{\alpha=1}^2 \int_0^L\Big(\int_\omega g^0 (S_1,S_2,.)S_\alpha dS_1dS_2\Big)\cdot({\cal R}\land\Gn_\alpha)-|\omega|\int_0^L \widetilde{f}^1\,{dU_E\over ds_3}\cdot\Gt.\cr}\right.$$
Notice that the above energy leads to a linear limit  model which couples the inextentional displacement, the torsion and the extentional displacement (see also (IV.4.24)).

The analysis of the present subsection and of Subsection IV.3 permits to state the following theorem.
\smallskip
\noindent{\ggras  Theorem  IV.4.3 }{\it The functional ${\cal J}_{LG}$ is the $\Gamma$-limit of  $\ds {J(.)-J(I_d)\over \delta^{6}}$ in the sense of Theorem IV.4.1. 

Let $({\cal U}_0,{\cal R}_0,{\cal V}_{S,0}, \overline{v}_0)$ be a minimizer of the functional ${\cal J}_{LG}$ over $\GU_{lin}$ and define $U_{E,0}\in D_{Ex}$ by (IV.4.16). Then  $({\cal U}_0,{\cal R}_0, U_{E,0})$ is the unique solution of the following minimization problem:
$$m_3={\cal F}_{LG}\big({\cal U}_0,{\cal R}_0,U_{E,0}\big)=\min_{\big({\cal U},{\cal R} , U_E\big)\in\GV\GG_{lin}}{\cal F}_{LG}\big({\cal U},{\cal R},U_E\big), \leqno(IV.4.21) $$ where 
$$\eqalign{
\GV\GG_{lin}=\Bigl\{ &\big({\cal U},{\cal R},U_E\big)\in (H^2(0,L))^3\times (H^1(0,L))^{3}\times D_{Ex}\;\; |\;\; 
{\cal U}(0)={\cal R}(0)=0,\quad U_E(0)=0,\quad {d{\cal U} \over ds_3}={\cal R}\land\Gt\Big\},\cr}$$ and 
$$\left\{\eqalign{
{\cal F}_{LG}\big({\cal U},{\cal R},U_E\big)=&{EI_1\over 2}\int_0^L\Big[{d{\cal R}\over ds_3} \cdot  \Gn_2\Big]^2+{EI_2\over 2}\int_0^L\Big[{d{\cal R}\over ds_3} \cdot  \Gn_1\Big]^2+{\mu K\over 4}\int_0^L\Big[{d{\cal R}\over ds_3} \cdot \Gt\Big]^2\cr
+& {E|\omega|\over 2}\int_0^L\Big[{dU_E\over ds_3} \cdot  \Gt+{1\over 2}\Big\|{d{\cal U}\over ds_3}\Big\|^2_2\Big]^2\cr
-&\int_0^L\Big(|\omega| f^0 +\sum_{\alpha=1}^2  \int_\omega g^0 (S_1,S_2,.)S_\alpha \det\big(\Gn_1\; |\; \Gn_2\;|\;{d\Gn_\alpha\over ds_3}\big)dS_1dS_2\Big)\cdot {\cal U}\cr
-&\sum_{\alpha=1}^2 \int_0^L\Big(\int_\omega g^0(S_1,S_2,.)S_\alpha dS_1dS_2\Big)\cdot({\cal R}\land\Gn_\alpha)-|\omega|\int_0^L \widetilde{f}^1\,{dU_E\over ds_3}\cdot\Gt\cr}\right.\leqno(IV.4.22)$$ $E$ is the Young's modulus  and $K$ is given in Theorem IV.1.3. Moreover we have
$$ {dU_{E,0}\over ds_3} \cdot  \Gt=-{1\over 2}\Big\|{d{\cal U}_0\over ds_3}\Big\|^2_2+{2\over E} \widetilde {f}^1\leqno(IV.4.23)$$ and the couple $({\cal U}_0,{\cal R}_0)\in \GV_{lin}$ is the unique solution of the following variational problem:
$$\left\{\eqalign{
&E\int_0^L\sum_{\alpha=1}^2 I_{3-\alpha}\Big[{d{\cal R}_0\over ds_3} \cdot  \Gn_\alpha\Big]\Big[{d{\cal R}\over ds_3} \cdot  \Gn_\alpha\Big]+{\mu K\over 2}\int_0^L\Big[{d{\cal R}_0\over ds_3} \cdot \Gt\Big]\Big[{d{\cal R}\over ds_3} \cdot  \Gt\Big]+ {|\omega|\over 2}\int_0^L\widetilde{f}^1\, {d{\cal U}_0\over ds_3}\cdot {d{\cal U}\over ds_3}\cr
&=\int_0^L\Big(|\omega| f^0 +\sum_{\alpha=1}^2  \int_\omega g^0 (S_1,S_2,.)S_\alpha \det\big(\Gn_1\; |\; \Gn_2\;|\;{d\Gn_\alpha\over ds_3}\big)dS_1dS_2\Big)\cdot {\cal U}\cr
&+\sum_{\alpha=1}^2 \int_0^L\Big(\int_\omega g^0(S_1,S_2,.)S_\alpha dS_1dS_2\Big)\cdot({\cal R}\land\Gn_\alpha),\qquad \forall \big({\cal U},{\cal R}\big)\in \GV_{lin}\cr}\right.\leqno(IV.4.24)$$}
\smallskip
\noindent Indeed the remarks at the end of Sections IV.2 and IV.3 are still valid for the above  chosen forces.
\medskip
\noindent {\Ggras  V.  Solutions of the non-linear minimization problem (IV.2.21)}
\medskip
The results of this subsection are limited to the case where the curved rod is   fixed only on $\Gamma_{0,\delta}$ (see Subsection II.2.4). As a consequence, the other extremity (for $s_3=L$) is  free (or with little change submitted to a given load). For these boundary conditions, we replace the minimization problem (IV.2.21) by an integro-differential equation satisfies by $\GR$.
To do that, we write the minimization problem (IV.2.21) in terms of the unknown $\GR$. We denote by $\GG$ the matrix of  $(L^2(0,L))^{3\times 3}$ such that
$$\eqalign{
 \int_0^L<\GG , \GR-\GI_3>&= 
\int_0^L\Big( |\omega| f +\sum_{\alpha=1}^2 \int_\omega g (S_1,S_2,.)S_\alpha \det\big(\Gn_1\; |\; \Gn_2\;|\;{d\Gn_\alpha\over ds_3}\big)dS_1dS_2\Big)\cdot({\cal V} -M)\cr
&+\sum_{\alpha=1}^2 \int_0^L\Big(\int_\omega g (S_1,S_2,.)S_\alpha dS_1dS_2\Big)\cdot(\GR-\GI_3)\Gn_\alpha,\cr}$$ for any $(\GR, {\cal V})\in \GV_{nlin}$, where $<\cdot , \cdot >$ is the inner product associated to the Frobenius norm over the  space $\GM_3$.

We set
$$\eqalign{
{\cal A}_3=&\Bigl\{\GA\in (L^2(0,L))^{3\times 3}\;|\enskip \GA^T(s_3)=-\GA(s_3)\enskip\hbox{for a.e. }\; s_3\in (0,L)\Big\}\cr
{\cal H}{\cal S}=&\Big\{\GR\in (H^1(0,L))^{3\times 3}\;|\enskip  \GR(0)=\GI_3\enskip\hbox{and for any }\; s_3\in [0,L],\quad \GR(s_3)\in SO(3)\Big\}.\cr}$$ 

\noindent Let $\GA$ be a  matrix  belonging to ${\cal A}_3$ and let $\GR_\GA$ be the solution of the Cauchy's problem
$$\left\{\eqalign{
&\GR_\GA\in (H^1(0,L))^{3\times 3},\cr
&  {d\GR_\GA\over ds_3}(s_3)=\GR_\GA(s_3)\GA(s_3),\quad \hbox{for a.e.} \; s_3\in (0,L),\cr
&\GR_\GA(0)=\GI_3.\cr}\right.\leqno(V.1)$$  The map $\GA\longmapsto \GR_\GA$ is one to one from ${\cal A}_3$ onto ${\cal H}{\cal S}$. An element $\GR\in {\cal H}{\cal S}$ is associated to the element $\ds \GA=\GR^T{d\GR\over ds_3}$ of ${\cal A}_3$.

\noindent Taking into account  the definition of $\GG$,   the minimum $m_2$ is in fact the minimum of the functional 
$${\cal F}_{nl}(\GR)={E\over 2}\int_0^L\sum_{\alpha=1}^2I_\alpha\Big({d\GR\over ds_3}\Gt\cdot \GR \Gn_\alpha\Big)^2+{\mu K\over 4}  \int_0^L\Big({d\GR\over ds_3}\Gn_1\cdot \GR \Gn_2\Big)^2-\int_0^L<\GG , \GR-\GI_3>\leqno(V.2)$$ over the closed set ${\cal H}{\cal S}$.  In terms of $\GA$, $m_2$ is also the minimum of the functional 
$${\cal G}(\GA)={E\over 2}\int_0^L\sum_{\alpha=1}^2I_\alpha\Big({d\GR_\GA\over ds_3}\Gt\cdot \GR_\GA \Gn_\alpha\Big)^2+{\mu K\over 4}  \int_0^L\Big({d\GR_\GA\over ds_3}\Gn_1\cdot \GR_\GA \Gn_2\Big)^2-\int_0^L<\GG , \GR_\GA-\GI_3>\leqno(V.3)$$ over the space ${\cal A}_3$. In view of (V.1),
we have
$${\cal G}(\GA)={E\over 2}\int_0^L\sum_{\alpha=1}^2I_\alpha\Big( \GA 
\Gt\cdot  \Gn_\alpha\Big)^2+{\mu K\over 4}  \int_0^L\Big( \GA \Gn_1\cdot \Gn_2\Big)^2
-\int_0^L<  \GG\,, \GR_\GA-\GI_3 >.\leqno(V.4)$$ In what follows we derive the first and the second derivatives of the last term in (V.4). By a standard  calculation we show  that for any matrices $\GA$ and $\GB$ in ${\cal A}_3$ we have
$$\eqalign{
\GR_{\GA+\GB}(s_3)=&\GR_\GA(s_3)+\Big(\int_0^{s_3}\GR_\GA(s)\GB(s)\GR_\GA^T(s)ds \Big)\GR_\GA(s_3)\cr
+&\Big(\int_0^{s_3}\int_0^s\GR_\GA(t)\GB(t)\GR_\GA^T(t)\GR_\GA(s)\GB(s)\GR_\GA^T(s)dtds\Big)\GR_\GA(s_3) 
+ O\big(|||\GB|||^3_{(L^2(0,L))^{3\times 3}}\big),\cr} $$ as the consequence we obtain
$${\cal G}(\GA+\GB )={\cal G}(\GA)+{\cal G}^{'}(\GA)(\GB)+{1\over 2}{\cal G}^{''}(\GA)(\GB,\GB)+O\big(||\GB||^3_{(L^2(0,L))^{3\times 3}}\big),$$ where
$$\left\{\eqalign{
&{\cal G}^{'}(\GA)(\GB)= E\int_0^L\sum_{\alpha=1}^2I_\alpha\big(\GA\Gt\cdot  \Gn_\alpha\big)\big(\GB\Gt\cdot  \Gn_\alpha\big)+{\mu K\over 2}  \int_0^L\big( \GA\Gn_1\cdot \Gn_2\big)\big( \GB\Gn_1\cdot \Gn_2\big)\cr
&\qquad - \int_0^L<\GG(s_3) \, , \Big(\int^{s_3}_0\GR_\GA(s)\GB(s)\GR_\GA^T(s)ds \Big)\GR_\GA(s_3) >ds_3\cr
&{\cal G}^{''}(\GA)(\GB,\GB)= E\int_0^L\sum_{\alpha=1}^2I_\alpha\big(\GB\Gt\cdot  \Gn_\alpha\big)^2+{\mu K\over 2}  \int_0^L\big(\GB\Gn_1\cdot \Gn_2\big)^2\cr
&\qquad- 2\int_0^L<\GG(s_3) \, , \Big(\int_0^{s_3}\int_0^s\GR_\GA(t)\GB(t)\GR_\GA^T(t)\GR_\GA(s)\GB(s)\GR_\GA^T(s)dtds\Big)\GR_\GA(s_3) >ds_3.\cr}\right.\leqno(V.5)$$ In order to explicit the minimum of ${\cal G}$, we simplify the term  involving the forces in ${\cal G}^{'}(\GA)(\GB)$. We have
$$\eqalign{
&\int_0^L<\GG(s_3) \, , \Big(\int^{s_3}_0\GR_\GA(s)\GB(s)\GR_\GA^T(s)ds \Big)\GR_\GA(s_3) >ds_3\cr
 =&\int_0^L<\GG(s_3)\GR_\GA^T(s_3)  \, , \Big(\int^{s_3}_0\GR_\GA(s)\GB(s)\GR_\GA^T(s)ds \Big)>ds_3.}$$ We integrate by parts the right hand side term in the above equality. This gives
$$\eqalign{
&\int_0^L<\GG(s_3) \, , \Big(\int^{s_3}_0\GR_\GA(s)\GB(s)\GR_\GA^T(s)ds \Big)\GR_\GA(s_3) >ds_3\cr
=& \int_0^L<\Big(\int_{s_3}^L\GG (s)\GR_\GA^T(s) ds\Big)\,,  \GR_\GA(s_3)\GB(s_3)\GR_\GA^T(s_3) >ds_3\cr
=& \int_0^L<\GR_\GA^T(s_3)\Big(\int_{s_3}^L\GG (s)\GR_\GA^T(s) ds\Big)\GR_\GA(s_3)\,,  \GB(s_3) >ds_3.}$$ Using the fact that symmetric and antisymmetric matrices are orthogonal for the scalar product $<\cdot,\cdot>$, we finally  get for any matrix $\GB\in {\cal A}_3$ 
$$\left\{\eqalign{
{\cal G}^{'}(\GA)(\GB)=&E\int_0^L\sum_{\alpha=1}^2I_\alpha\Big(\GA\Gt\cdot  \Gn_\alpha\Big)\Big(\GB
\Gt\cdot  \Gn_\alpha\Big)+{\mu K\over 2}  \int_0^L\Big( \GA\Gn_1\cdot \Gn_2\Big)\Big( \GB\Gn_1\cdot \Gn_2\Big)\cr
-&\int_0^L<\GR_\GA^T(s_3)\Big(\int_{s_3}^L{1\over 2}\big[\GG (s)\GR_\GA^T(s)-\GR_\GA(s)\GG^T(s)\big]ds
\Big)\GR_\GA(s_3)\,,  \GB(s_3) >ds_3.\cr}\right.\leqno(V.6)$$  
The above derivations allow to prove the following theorem.
\smallskip
\noindent{\ggras Theorem V.1. }{\it   Let $({\cal V}_0,\GR_0)$  be in  $\GV_{nlin}$ and set 
$\GA_0=\ds\GR^T_0{d\GR_0\over ds_3} $. Then $({\cal V}_0,\GR_0)$ is a solution of 
the minimization problem (IV.2.21) if and only if 
$\GR_0$ is a solution of the following integro-differential problem
$$\left\{\eqalign{
\GA_0(s_3)\Gn_1(s_3)\cdot  \Gn_2(s_3)&={2\over \mu K}\GR^T_0(s_3)\Big(\int_{s_3}^L \big[\GG\GR^T_0 -\GR_0 \GG^T \big] \Big) \GR_0(s_3)\Gn_1(s_3)\cdot \Gn_2(s_3)\cr
 \GA_0(s_3)\Gt(s_3)\cdot  \Gn_1(s_3)&={1\over EI_1}\GR^T_0(s_3)\Big(\int_{s_3}^L \big[\GG\GR^T_0-\GR_0\GG^T\big] \Big) \GR_0(s_3)\Gt(s_3)\cdot  \Gn_1(s_3)\cr
 \GA_0(s_3)\Gt(s_3)\cdot  \Gn_2(s_3)&={1\over EI_2}\GR^T_0(s_3)\Big(\int_{s_3}^L \big[\GG\GR^T_0-\GR_0\GG^T\big]\Big) \GR_0(s_3)\Gt(s_3)\cdot  \Gn_2(s_3).\cr}\right.\leqno(V.7)$$
Moreover, if
$$||\GG||_{(L^2(0,L))^{3\times 3}}<{1\over  L^{3/2}} \inf\Big(EI_1,EI_2,{\mu K\over 2}\Big)\leqno(V.8)$$  the solution of the minimization problem (IV.2.21) is unique.}
\smallskip
\noindent{\ggras Proof. } An element  $({\cal V}_0,\GR_0)$ of $\GV_{nlin}$  is a minimizer of (IV.2.21)   only if $\GA_0$ is a minimizer of  the functional ${\cal G}$ given by (V.3).
Hence,   we have ${\cal G}^{'}(\GA_0)(\GB)=0$ for any $\GB\in {\cal A}_3$. In view of (V.5) and (V.6), the antisymmetric matrix  $\GA_0$ satisfies
$$\eqalign{
&E\int_0^L\sum_{\alpha=1}^2I_\alpha\Big(\GA_0\Gt\cdot  \Gn_\alpha\Big)\Big(\GB
\Gt\cdot  \Gn_\alpha\Big)+{\mu K\over 2}  \int_0^L\Big( \GA_0\Gn_1\cdot \Gn_2\Big)\Big( \GB\Gn_1\cdot \Gn_2\Big)\cr
=&\int_0^L<\GR^T_0(s_3)\Big(\int_{s_3}^L{1\over 2}\big[\GG (s)\GR^T_0-\GR\GG^T_0\big]
\Big)\GR_0(s_3)\,,  \GB(s_3) >ds_3,\qquad \forall \GB\in {\cal A}_3.}$$ This immediately gives (V.7).

\noindent Now we prove that the functional ${\cal G}$ admits a unique minimizer, under the assumption (V.8). For any $ \GA\in{\cal A}_3$ we get
$$\big\|\GA\big\|^2_{(L^2(0,L)^{3\times 3}}= 2\Big\{\big\|\GA\Gt\cdot  \Gn_1\big\|^2_{L^2(0,L)}+\big\|\GA\Gt\cdot  \Gn_2\big\|^2_{L^2(0,L)}+\big\|\GA\Gn_1\cdot \Gn_2\big\|^2_{L^2(0,L)}\Big\}.$$ 
\noindent From the expression (V.5) of ${\cal G}^{''}(\GA)(\GB,\GB)$ and the above equality we have
$${\cal G}^{''}(\GA)(\GB,\GB)\ge {1\over 2}\Big\{\inf\Big(EI_1,EI_2,{\mu K\over 2}\Big)-L^{3/2}||\GG||_{(L^2(0,L)^{3\times 3}}\Big\}||\GB||^2_{(L^2(0,L))^{3\times 3}}.\leqno(V.9)$$ As a consequence of the above inequality, if $\GG$ satisfies (V.8) the functional ${\cal G}$ is strictly convex, which insures the uniqueness of the minimizer $\GA_0$.\fin

\bigskip

\noindent{\ggras  Appendix. A few recalls on rotations }
\medskip
\noindent Let  $\GV$ be a matrix belonging to $SO(3)$. The matrix $\GV$ is the matrix of a rotation  ${\cal R}_{\Ga,\theta}$ in $\R^3$ where $\Ga$ is a unit vector belonging to the axis of the rotation and where $\theta$ belonging to $[0,\pi]$  is the angle of rotation about this axis. The rotation is written  as
$$\forall \Gx\in \R^3,\qquad {\cal R}_{\Ga,\theta}(\Gx)=\cos(\theta)\Gx+(1-\cos(\theta))<\Gx,\Ga>\Ga+\sin(\theta)\,\Ga\land\Gx.$$ We have 
 $$ |||\GI_3-\GV|||=2\sqrt 2\sin\Big({\theta\over 2}\Big)\ge {2\sqrt 2\over \pi}\theta.$$ For all $t\in [0,1]$, we denote by  $\GW(t)$ the matrix of the rotation ${\cal R}_{\Ga, t\theta}$. The function $t\to \GW(t)$   belongs to  $\big({\cal C}^1([0,1])\big)^{3\times 3}$  and satisfies
$$\GW(0)=\GI_3,\quad \GW(1)=\GV,\qquad\quad\GW(t)\in SO(3),\qquad \Big| \Big| \Big|{d\GW\over dt}(t)\Big| \Big| \Big|=\sqrt 2 \theta\le 2|||\GI_3-\GV|||,\qquad t\in [0,1].$$
 \noindent Now, let $\GU_0$ and  $\GU_1$ be two elements in  $SO(3)$. We set
$$\GV=\GU_0^{-1}\GU_1 $$ and we consider  the map
$$\GU(t)=\GU_0 \GW(t)\qquad\hbox{for any }\enskip t\in [0,1],$$ where $\GW(t)$ is defined above. We have built a path $\GU\in\big({\cal C}^1([0,1])\big)^{3\times 3}$ such that 
$$\eqalign{
&\GU(0) =\GU_0,\quad \GU(1)=\GU_1,\qquad  \GU(t)\in SO(3),\cr
& \Big| \Big| \Big|{d\GU\over dt}(t)\Big| \Big| \Big|\le 2|||\GU_1-\GU_0|||,\qquad t\in [0,1].\cr}$$
\noindent {\ggras  Lemma A. }{\it Let $\GR$ be in $(H^1(0,L))^{3\times 3}$ such that $\GR(0)=\GI_3$ and such that for any $s_3\in [0,L]$ the matrix $\GR(s_3)$ belongs to $SO(3)$. There exists a sequence of matrices $(\GR_N)_{N\in \N^*}$ satisfying $\GR_N\in  (W^{1,\infty}(0,L))^{3\times 3}$, $\GR_N(0)=\GI_3$ and for any $s_3\in [0,L]$ the matrix $\GR_N(s_3)$ belongs to $SO(3)$ and moreover
$$\GR_N\longrightarrow \GR\quad \hbox{strongly in }\quad (H^1(0,L))^{3\times 3}.$$}
\noindent{\ggras Proof.}  The matrix $\GA=\ds\GR^T{d\GR\over ds_3}$ is antisymmetric and belongs to $(L^2(0,L))^{3\times 3}$. Let  $\big(\GA_N\big)_{n\in \N}$ be a sequence of antisymmetric matrices such that
$$\GA_N\in ({\cal C}([0,L]))^{3\times 3}\qquad \hbox{and}\qquad \GA_N\longrightarrow \GA\quad \hbox{strongly in }\quad (L^2(0,L))^{3\times 3}.$$ Let $\GR_N$ ($N\in \N$) be the solution of the Cauchy's problem
$$\left\{\eqalign{
&{d\GR_N\over ds_3}=\GR_N\GA_N\cr
&\GR_N(0)=\GI_3\cr}\right.$$ We have $\GR_N\in ({\cal C}^1([0,L]))^{3\times 3}$ and for any $s_3\in [0,L]$ the matrix $\GR_N(s_3)$ belongs to $SO(3)$. From the above strong convergence we deduce that
$$\GR_{N}\longrightarrow \GR\quad \hbox{strongly in }\quad (H^1(0,L))^{3\times 3}.$$\fin
\bigskip
\noindent [1] E. Acerbi, G. Buttazzo and D. Percivale, A variational definition for the strain energy of an elastic string. J. Elast. 25 (1991) 137-148.

\noindent [2] S. S. Antman, The theory of rods. In: Truesdell, C.A. (ed.) Handbuch der Physik, vol. VIa, pp. 641-703, Springer, Berlin Heidelberg New York (1972).

\noindent [3] J.M. Ball, Convexity conditions and existence theorems in nonlinear elasticity. Arch. Ration. Mech. Anal. 63 (1976) 337-403.

\noindent [4] D. Blanchard, A. Gaudiello, G. Griso.  Junction of a periodic
family of elastic rods with a $3d$ plate. I. J. Math. Pures Appl. (9) 88 (2007), no 1, 149-190.

\noindent [5]   D. Blanchard, A. Gaudiello, G. Griso. Junction of a periodic
family of elastic rods with a thin plate. II. J. Math. Pures Appl. (9) 88 (2007), no 2, 1-33.

\noindent [6]  D. Blanchard,  G. Griso. Microscopic effects in the homogenization of the junction of rods and a thin plate.  Asympt. Anal. 56 (2008), no 1, 1-36.

\noindent [7] P.G. Ciarlet, Mathematical Elasticity, Vol.1, North-Holland, Amsterdam (1988).

\noindent [8] P.G. Ciarlet and P. Destuynder, A justification of a nonlinear model in plate theory. Comput. Methods Appl. Mech. Eng. 17/18 (1979) 227-258.

\noindent [9] G. Dal Maso: An Introduction to $\Gamma$-convergence. BirkhŠuser, Boston, 1993.

\noindent [10] L. Evans and R.F. Gariepy, Measure theory and finepropertiesof functions. Studies in Advanced Mathematics. CRC, Boca Raton, Fla., 1992.

\noindent [11] G. Friesecke, R. D. James and S. MŸller. A theorem on geometric rigidity and the derivation of nonlinear plate theory from the three-dimensional elasticity. Communications on Pure and Applied Mathematics, Vol. LX, 1461-1506 (II.1002).

\noindent [12] G. Friesecke, R. D. James and S. MŸller, A hierarchy of plate models derived from nonlinear elasticity by $\Gamma$-convergence. (2005)

\noindent [13] G. Griso. Decomposition of displacements of thin structures.  J. Math. Pures Appl.  89 (2008) 199-233.

\noindent [14] G. Griso. Asymptotic behavior of curved rods by the unfolding method. Math. Meth. Appl. Sci. 2004; 27: 2081-2110.

\noindent [15] J.J. Marigo, N. Meunier, Hierarchy of one-dimensional models in nonlinear elasticity. Journal of elasticity (2006) 83: 1-28.

\noindent [16] H. Le Dret, Convergence of displacements and stresses in linearly
elastic slender rods as the thickness goes to zero, Asymptot. Anal. 10  (1995) 367-402.

\noindent [17] H. Le Dret and A. Raoult, The nonlinear membrane model as variational limit of nonlinear three-dimensional elasticity. J. Math. Pures Appl. 75 (1995) 551-580.

\noindent [18] H. Le Dret and A. Raoult, The quasiconvex envelope of the Saint Venant-Kirchhoff  stored energy function. Proc. R. Soc. Edin., A 125 (1995) 1179-1192.

\noindent [19] M.G. Mora and S. MŸller, Derivation of the nonlinear bending-torsion theory for inextensible rods by Gamma-convergence. Calc. Var. 18 (2002) 287-305.

\noindent [20] M.G. Mora and S. MŸller, A nonlinear model for inextensible rods as a low energy $\Gamma$-limit of three-dimensional nonlinear elasticity. Ann. Inst. Henri PoincarŽ, Anal. non linŽaire 21 (2004) 271-299.

\noindent [21] O. Pantz,  On the justification of the nonlinear inextensional plate model. C. R. Acad. Sci. Paris SŽr. I Math. 332 (2001), no. 6, 587--592.

\noindent [22] O. Pantz, Le modle de poutre inextensionnelle comme limite de l'ŽlasticitŽ
non-linŽaire tridimensionnelle. Preprint,  2002.

\noindent [23] O. Pantz,  On the justification of the nonlinear inextensional plate model. Arch. Ration. Mech. Anal. 167 (2003), no. 3, 179--209. 

\noindent [24] L. Trabucho and J.M. Via$\tilde {\rm n} $o, Mathematical modelling of rods. In P.G. Ciarlet and J.L. Lions (eds), Handbook of Numerical Analysis, Vol. IV. North-Holland, Amsterdam (1996).

\bye
\vfill\eject
\noindent{\ggras ComplŽments } La poutre est  fixŽe aux deux extrŽmitŽs. 
\medskip
\noindent Dans la dŽcomposition d'une dŽformation on peut toujours choisir $\GR=\GI_3$  au voisinage des extrŽmitŽs de la poutre (voir Section II.2.4) les estimations (II.2.2) du ThŽorme II.2.2 restent inchangŽes. La dŽformation "moyenne" ${\cal V}$ est dŽfinie (II.2.5)  par la moyenne de $v$ sur la section droite de la poutre. Si la poutre est  fixŽe en ses deux extrŽmitŽs, on a donc toujours  ${\cal V}(0)=M(0)$ et ${\cal V}(L)=M(L)$.

\noindent On dŽfinit maintenant ${\cal V}_B$ et ${\cal V}_S$ par (voir (II.2.16) )
$$\forall s_3\in [0,L],\qquad {\cal V}_B(s_3)=M(0)+\int_0^{s_3}\GR(z)\Gt(z) dz-{s_3\over L}\int_0^{L}(\GR(z)-\GI_3)\Gt(z) dz\qquad {\cal V}_S(s_3)={\cal V}(s_3)-{\cal V}_B(s_3)$$ La dŽformation ${\cal V}_B$ appartient ˆ $(H^2(0,L))^{3}$ et elle vŽrifie ${\cal V}_B(0)=M(0)$, ${\cal V}_B(L)=M(L)$. Le dŽplacement ${\cal V}_S$ appartient ˆ $(H^1(0,L))^3$ et il vŽrifie ${\cal V}_S(0)={\cal V}_S(L)=0$. D'aprs (II.2.2) on a
$$\Big\|\int_0^{L}(\GR(z)-\GI_3)\Gt(z) dz\Big\|_2\le {C\over \delta}||dist\big(\nabla_x v, SO(3)\big)||_{L^2({\cal P}_\delta)}.$$ D'o
$$ \Big\|{d{\cal V}_S\over ds_3}\Big\|_{(L^2(0,L))^3}\le  {C\over \delta}||dist\big(\nabla_x v, SO(3)\big)||_{L^2({\cal P}_\delta)}.$$ La dŽformation ${\cal V}_B$ satisfait l'estimation
$$ \Big\|{d{\cal V}_B\over ds_3}-\GR\Gt\Big\|_{(L^2(0,L))^3}\le  {C\over \delta}||dist\big(\nabla_x v, SO(3)\big)||_{L^2({\cal P}_\delta)}.$$
Dans le lemme III.2.1, on a maintenant
$${\cal V}(0)=M(0),\quad {\cal V}(L)=M(L),\qquad {\cal V}_S\in (H^1_0(0,L))^3,\qquad {d{\cal V}\over ds_3}=\GR\Gt$$ et de plus
$$ \GR\in (H^1(0,L))^{3\times 3},\qquad \forall s_3\in [0,L]\quad \GR(s_3)\in SO(3),\qquad \GR(0)=\GR(L)=\GI_3,\qquad \int_0^{L}(\GR(z)-\GI_3)\Gt(z) dz=0.$$
\vfill\eject
\noindent{\Ggras III.2.* Limit behavior of the deformation for  $||\hbox{dist}(\nabla_x v,SO(3))||_{L^2({\cal P}_\delta)} \sim \delta^{1+\kappa/2}$, $0<\kappa\le 2 $}
  \medskip
Let us consider a sequence of deformations $v_\delta$ of $\big(H^1({\cal P}_\delta)\big)^3 $ such that 
$$||\hbox{dist}(\nabla_x v_\delta,SO(3))||_{L^2({\cal P}_\delta)} \le C \delta^{1+\kappa/2}.\leqno(III.2.1*)$$
We denote by ${\cal V}_\delta$, $\GR_\delta$ and $\overline{v}_\delta$ the three terms of the decomposition of $v_\delta$ given by Theorem II.2.2 and by ${\cal V}_{B,\delta}$ and ${\cal V}_{S,\delta}$  the two terms given by (II.2.16).  The estimates (II.2.2), (II.2.17), (II.3.5), (II.3.6), (III.11), (III.13) lead to the following lemma:

\noindent {\ggras  Lemma III.2.1.* }
{\it There exists a subsequence still indexed by $\delta$ such that 
$$\left\{\eqalign{
&\GR_\delta \rightharpoonup \GR \quad \hbox{weakly in}\quad \big(L^2(0,L)\big)^{3\times 3}\cr
&\delta^{1-\kappa/2}\GR^T_\delta{d\GR_\delta\over ds_3} \rightharpoonup \GA \quad \hbox{weakly in}\quad \big(L^2(0,L)\big)^{3\times 3}\cr
&{\cal V}_\delta\rightharpoonup  {\cal V} \quad \hbox{weakly in}\quad \big(H^1(0,L)\big)^3\cr
&{1\over \delta^{\kappa/2}}\GR^T_\delta\Big({d{\cal V}_{\delta}\over ds_3}-\GR_\delta\Gt\Big)\rightharpoonup {\cal Z}\quad \hbox{weakly in}\quad \big(L^2(0,L)\big)^3\cr
&{1\over \delta^{1+\kappa/2}}\GR^T_\delta\big(\Pi _\delta\overline {v}_\delta\big)\rightharpoonup  \overline {w} \quad \hbox{weakly in}\quad \big(L^2(0,L;H^1(\omega))\big)^3\cr}\right.\leqno(III.2.2*)$$ For a.e. $s_3\in ]0,L[$ the matrix $\GA(s_3)$ is antisymmetric and  
$$||\GR(s_3)\Gt(s_3)||_2\le 1,\qquad ||\GR(s_3)\Gn_\alpha(s_3)||_2\le 1.$$
 Moreover we have
$${\cal V}(0)=M(0), \qquad  {d{\cal V}\over ds_3}=\GR \Gt .\leqno (III.2.3*)$$ 
\noindent Furthermore, we also have
$$\left\{\eqalign{
\Pi _\delta v_\delta&\rightharpoonup  {\cal V} \quad \hbox{weakly in}\quad \big(H^1(\Omega)\big)^3,\cr
\Pi _\delta (\nabla_x v_\delta)& \rightharpoonup  \GR\quad \hbox{weakly in}\quad \big(L^2(\Omega)\big)^{3\times 3}.\cr}\right.\leqno(III.2.4*)$$}
\smallskip
\noindent{\ggras  Corollary III.2.2.* } {\it For the same subsequence we  have
$$\left\{\eqalign{
{1\over \delta^{\kappa/2}}\GR^T_\delta\big(\Pi _\delta (\nabla_x v_\delta)-\GR_\delta\big)\Gn_\alpha& \rightharpoonup  {\partial\overline{w}\over \partial S_\alpha}\quad \hbox{weakly in}\quad (L^2(\Omega))^3,\cr
{1\over \delta^{\kappa/2}}\GR^T_\delta\big(\Pi _\delta (\nabla_x v_\delta)-\GR_\delta\big)\Gt &\rightharpoonup  \GA\big(S_1\Gn_1+S_2\Gn_2)+{\cal Z} \quad \hbox{weakly in}\quad (L^2(\Omega))^3.\cr}\right.\leqno(III.2.5*)$$}
\noindent{\ggras  Proof. } The first convergence in (III.2.5*) is a direct consequence of (II.2.13) and (III.2.2*). In order to obtain the second convergence, remark first that, thanks to  estimates (III.11)  and (III.2.1*) the sequence $\ds{1\over \delta^{\kappa/2}}\GR^T_\delta\Pi_\delta\overline{v}_\delta$  is bounded in $H^1(\Omega)$. Due to (III.2.2*), its weak limit must be equal to $0$. Using now (II.2.15) and the convergences  (III.2.2*) leads to the result.\fin
\noindent 
To end this section, let us notice that  the relation 
$$(\nabla_xv_\delta)^T\nabla_x v_\delta-\GI_3=(\nabla_xv_\delta-\GR_\delta)^T(\nabla_x v_\delta-\GR_\delta)+(\nabla_xv_\delta-\GR_\delta)^T\GR_\delta+(\GR_\delta)^T(\nabla_x v_\delta-\GR_\delta)$$ permits to obtain the limit of the Green-St Venant's tensor in the rescaled domain $\Omega$
$${1\over 2\delta^{\kappa/2}}\Pi _\delta \big((\nabla_xv_\delta)^T\nabla_x v_\delta-\GI_3\big)\rightharpoonup   \GE \qquad\hbox{weakly in}\quad (L^1(\Omega))^{3\times 3},\leqno(III.2.6*)$$ where 
$$\eqalign{
\GE=&{1\over 2}\Big\{(\Gn_1\,|\, \Gn_2\,|\, \Gt) \Bigl({\partial\overline{w}\over \partial S_1}\,|\, {\partial\overline{w}\over \partial S_2}\,|\, \GA\big(S_1\Gn_1+S_2\Gn_2)+ {\cal Z}\Big)^T\cr
+& \Bigl({\partial\overline{w}\over \partial S_1}\,|\, {\partial\overline{w}\over \partial S_2}\,|\,\GA\big(S_1\Gn_1+S_2\Gn_2)+ {\cal Z}\Big) (\Gn_1\,|\, \Gn_2\,|\, \Gt)^T\Big\}\cr}\leqno(III.2.7*)$$ We can write $\GE$ in the form
$$\GE=(\Gn_1\,|\, \Gn_2\,|\, \Gt)\, \widehat{\GE}\,(\Gn_1\,|\, \Gn_2\,|\, \Gt)^T,$$ where, using the fact that the matrix $\GA$ is antisymmetric  
$$\widehat{\GE}= \pmatrix{
\ds{\partial\overline{w}\over \partial S_1}\cdot\Gn_1 & \ds{1\over 2}\Big\{{\partial\overline{w}\over \partial S_1}\cdot\Gn_2+{\partial\overline{w}\over \partial S_2}\cdot\Gn_1\Big\} & \ds {1\over 2}\Big\{{\partial\overline{w}\over \partial S_1}\cdot\Gt-S_2 \GA\Gn_1\cdot \Gn_2+{\cal Z}\cdot \Gn_1 \Big\}\cr
* & \ds{\partial\overline{w}\over \partial S_2}\cdot\Gn_2   &  \ds {1\over 2}\Big\{{\partial\overline{w}\over \partial S_2}\cdot\Gt+S_1\GA\Gn_1\cdot \Gn_2+ {\cal Z}\cdot \Gn_2\Big\}\cr
 * & *  & \ds \phantom{\ds\int} -S_1\GA\Gt\cdot  \Gn_1-S_2\GA\Gt\cdot  \Gn_2+ {\cal Z} \cdot \Gt}.\leqno(III.2.8*)$$ The matrix $\widehat{\GE}$ is symmetric.
\bye